\documentclass[11pt]{amsart}
\usepackage[T1]{fontenc}
\IfFileExists{lmodern.sty}{\usepackage{lmodern}}{}
\usepackage[margin=1.05in]{geometry}
\usepackage{amsmath,amssymb,amsthm,mathrsfs,bm}
\usepackage{longtable,booktabs,array,enumitem}
\usepackage[expansion=false]{microtype}
\usepackage{tikz}
\usepackage[bookmarks=false,colorlinks=true,linkcolor=blue,citecolor=blue,urlcolor=blue]{hyperref}
\setlist[enumerate]{leftmargin=*,itemsep=3pt}
\newtheorem{theorem}{Theorem}[section]
\newtheorem{proposition}[theorem]{Proposition}
\newtheorem{corollary}[theorem]{Corollary}
\newtheorem{lemma}[theorem]{Lemma}
\theoremstyle{definition}
\newtheorem{definition}[theorem]{Definition}

\theoremstyle{remark}
\newtheorem{remark}[theorem]{Remark}
\newtheorem*{ac}{Acknowledgements}
\newcommand{\bR}{\mathbb R}

\newcommand{\bZ}{\mathbb Z}
\newcommand{\T}{\mathbb T}
\newcommand{\bC}{\mathbb C}
\newcommand{\bG}{\mathbb G}
\newcommand{\hbG}{\widehat{\mathbb G}}
\newcommand{\cB}{\mathcal B}
\newcommand{\cM}{\mathcal M}
\newcommand{\hcM}{\widehat{\mathcal M}}
\newcommand{\cH}{\mathcal H}
\newcommand{\cR}{\mathcal R}
\newcommand{\cS}{\mathcal S}
\newcommand{\hcS}{\widehat{\mathcal S}}
\newcommand{\cF}{\mathcal F}
\newcommand{\tcF}{\widetilde{\mathcal F}}
\newcommand{\hcF}{\widehat{\mathcal F}}
\newcommand{\hnabla}{\widehat\nabla}
\newcommand{\Lamv}{\Lambda_\varphi}
\newcommand{\hLam}{\widehat\Lambda}
\newcommand{\hd}{\hat d_{\hat\varphi}}

\DeclareMathOperator{\Pol}{Pol}

\DeclareMathOperator{\Tr}{Tr}

\newcommand{\id}{\iota}
\newcommand{\Irr}{\operatorname{Irr}}

\title{Scaling Fourier Inequalities for Locally Compact Quantum Groups}

\author{Jinsong Wu}
\address{Jinsong Wu, Beijing Institute of Mathematical Sciences and Applications, Beijing, 101408, China}
\email{wjs@bimsa.cn}

\author{Fulin Yang}
\address{Fulin Yang, Beijing Institute of Mathematical Sciences and Applications, Beijing, 101408, China}
\address{Yau Mathematical Sciences Center, Tsinghua University, Beijing 100084, China}
\email{fulinyang@bimsa.cn}
\date{}

\begin{document}

\begin{abstract}
In this paper, we study Fourier and convolution inequalities on locally compact quantum groups in the presence of the scaling group. 
For $1\le p\le 2$, we investigate a family of scaling Fourier transforms and the associated norms $H_{p,\theta}(\mathbb G)$, and we show that $H_{p,\theta}(\mathbb G)<1$ for the real quantum $az+b$ group and $0<\theta<1$, $1<p<2$.
We also investigate scaling convolution for locally compact quantum groups and, by extending the ordinary convolution, obtain the corresponding scaling Young inequality in four different regions. 
In Young Region~II, under $1<p,s<2$ and $2\leq r<\infty$, we construct sequences proving norm one at every effective vertex and prove strict inequality at every other effective point.
\end{abstract}

\maketitle

\section{Introduction}

The Hausdorff--Young inequality and Young's inequality are fundamental inequalities in classical Fourier analysis. 
In 1958, Kunze \cite{Kunze1958} proved the Hausdorff--Young inequality for unimodular locally compact groups. 
However, it was challenging to determine the best constant for the Hausdorff--Young inequality when a locally compact group is given. 
In 1961, Babenko \cite{Bab61} studied the Hausdorff--Young inequality for the real line and found that the best constant should be strictly less than one. 
In 1975, Beckner \cite{Bec75} proved the sharp Hausdorff--Young inequality for the real line and found that Gaussian functions are extremizers. 
Lieb \cite{Lie90} showed that the extremizers are Gaussian functions. For unimodular locally compact groups, Russo \cite{Rus74} showed that the extremizers are modulations and translations of compact open subgroups for the Hausdorff-Young inequality with optimal constant one. 
Terp \cite{Terp2017} obtained the Hausdorff--Young inequality for locally compact groups. 
Klein and Russo \cite{KR} proved the sharp Hausdorff--Young inequality for Heisenberg groups when the conjugate exponent \(p'\) is an even integer. 
Cowling, Martini, M\"uller, and Parcet \cite{CMMP2019} studied the Hausdorff--Young inequality for semisimple Lie groups. 
Recently, Qiu and Xiong \cite{QX} proved the sharp Hausdorff--Young inequality for noncommutative Euclidean spaces and Heisenberg groups. Inspired by the work of Qiu and Xiong, the authors \cite{WY26a} obtained sharp Hausdorff--Young inequalities for solvable and semisimple finite-center Lie groups.

A fundamental result for the Fourier transform on locally compact groups is that it interchanges multiplication and convolution. 
Convolution for locally compact groups is very helpful in studying the Fourier transform. 
The inequality related to convolution is called Young's inequality, and it has been well studied for locally compact groups in \cite{Bec75,Nei63,Yap69,You12}.

Kac algebras \cite{EnoSch92,VaiKac74} generalized locally compact groups and their duals. 
Locally compact quantum groups, introduced by Kustermans and Vaes \cite{KusVae00}, are generalizations of Kac algebras. 
The Hausdorff--Young inequality for locally compact quantum groups was obtained by Cooney \cite{Coo10} and Caspers \cite{Cas13}, with the Hausdorff--Young constant bounded by $1$. 
It is a natural question to find non-Kac quantum groups whose Hausdorff--Young constant is strictly less than $1$ for interior parameters. 
Recently, the authors \cite{WY26b} obtained the sharp Hausdorff--Young inequality for some locally compact quantum groups, showing that the Hausdorff--Young constant for Kahng's quantum group \cite{Kahng} is strictly less than $1$ for interior parameters. 
One consequence of the Hausdorff--Young inequality is the uncertainty principle \cite{JLW18,LiuWu17}. 
To characterize the extremizers of the uncertainty principle in tracial cases, Jiang, Liu, and the first author suggested the notion of bi-shifts of group-like projections, which generalize modulations and translations of compact open subgroups, and showed that the minimizers of the uncertainty principle are bi-shifts of group-like projections. 
For unimodular Kac algebras \cite{EnoSch92,VaiKac74}, Liu and the first author \cite{LiuWu18} proved that the extremizers of the Hausdorff--Young inequality with optimal constant one are bi-shifts of group-like projections.

Convolution for locally compact quantum groups was first defined for linear functionals (viewed as measures) and operators (viewed as functions) via the comultiplication. 
With the help of the noncommutative $L^p$ space $L^p(\bG)$ for a locally compact quantum group $\bG$, the convolution of $x\in L^p(\bG)$ and $y\in L^s(\bG)$ is well-defined for $1\le p,s\le2$ \cite{JLW18,LWW17}. 
However, it was shown that the convolution defined in \cite{LWW17} is not suitable for all $1\le p,s\le\infty$. 
Moreover, they showed that the Fourier transform interchanges convolution and multiplication for locally compact quantum groups, and they established Young's inequality for convolution. 
Convolution and Young's inequality have also been studied for other quantum symmetries. 
Jiang, Liu, and the first author \cite{JLW16,JLW19} systematically investigated convolution and Young's inequality for subfactor planar algebras and showed that the extremizers of the Hausdorff--Young inequality and Young's inequality and the minimizers of the uncertainty principle are bi-shifts of biprojections. 
In category theory, Liu, Palcoux, and the first author \cite{LPW19} studied Fourier analysis for fusion categories and found that some inequalities in Fourier analysis become obstructions to unitary categorification of fusion rings.
A program of quantum Fourier analysis was suggested in \cite{JJL20}.

In this paper, we show that there is a family of Fourier transforms in the presence of scaling automorphism groups, which we call scaling Fourier transforms.
These Fourier transforms are equivalent to the twisted Fourier transforms introduced by Youn \cite{You26} for compact quantum groups. 
The Hausdorff--Young inequality for the twisted Fourier transforms on compact quantum groups is sharp with optimal constant $1$. 
The scaling Fourier transform exhibits much more interesting phenomena for locally compact quantum groups. 
We show that the Hausdorff--Young constant of the scaling Fourier transform is strictly less than $1$ for interior parameters for the real quantum $az+b$ group.
This gives a different perspective of nontriviality of the Hausdorff--Young inequality for real quantum $az+b$, where its usual Hausdorff-Young constant is one \cite{WY26b}.
We also characterize the extremizers of the scaling Hausdorff--Young inequality with optimal constant $1$ for interior parameters for unimodular locally compact quantum groups; 
they are scaling-invariant bi-shifts of group-like projections. 
We also extend ordinary convolution for locally compact quantum groups to scaling convolution and extend it to the full range.
Young's inequalities for scaling convolution naturally follow. 
We also characterize the pairs of extremizers of the scaling Young inequality with optimal constant $1$ for interior parameters for unimodular locally compact quantum groups. 
We show that the Young constant of scaling convolution is strictly less than $1$ for interior parameters for the real quantum $az+b$ group.
We also generalize the Nikolskii inequality to locally compact quantum groups.

The paper is organized as follows. 
In Section 2, we recall some basic facts about noncommutative $L^p$ spaces, locally compact quantum groups in the sense of Kustermans--Vaes, their ordinary Fourier transform, and bi-shifts of group-like projections. 
In Section 3, we introduce the scaling Fourier transform and its Hausdorff--Young inequality. 
In Section 4, we study equality in the Hausdorff--Young inequality for the ordinary Fourier transform and scaling Fourier transforms.
In Section 5, we study the Hausdorff--Young constant for locally compact quantum groups with periodic scaling groups. 
In Section 6, we extend ordinary convolution to convolution in the full range. 
In Section 7, we introduce scaling convolution and characterize the equality cases of Young's inequalities for unimodular locally compact quantum groups. 
In Section 8, we investigate the Nikolskii inequality for locally compact quantum groups. 
In Section 9, we study the scaling Hausdorff--Young inequality and scaling Young inequalities for the real quantum $az+b$ group.

\begin{ac}
J. Wu was supported by grants from Beijing Institute of Mathematical Sciences and Applications.
J.~Wu was supported by NSFC (Grant no. 12371124).  
Fulin Yang is supported by NSFC (Grant no. 12601240)
\end{ac}

\section{Preliminaries}

\subsection{The spatial $L^p$ spaces}
Let $\cM$ be a von Neumann algebra with a normal semifinite faithful
weight $\varphi$.  Put
\[
\mathfrak N_\varphi=\{a\in\cM:\varphi(a^*a)<\infty\},\qquad
\mathfrak M_\varphi=\operatorname{span}(\mathfrak N_\varphi^*\mathfrak N_\varphi).
\]
We use the GNS representation $(\cH_\varphi,\Lambda_\varphi)$, with inner products $\langle \cdot, \cdot\rangle_\varphi$ linear in the first variable and anti-linear in the second variable.  
If there is no confusion, we shall denote $\langle \cdot, \cdot\rangle_\varphi$ by $\langle \cdot, \cdot\rangle$.
Its modular data are modular conjugation $J_\varphi$, modular operator $\nabla_\varphi$ and modular automorphism $\sigma^\varphi$.  
In particular, whenever the right side is defined,
\begin{equation}\label{eq:correct-right-GNS}
 \Lambda_\varphi(xa)
 =J_\varphi\sigma^\varphi_{i/2}(a)^*J_\varphi\Lambda_\varphi(x).
\end{equation}
For example, this holds for $x\in\mathfrak N_\varphi$ and
$a\in\mathcal D(\sigma^\varphi_{i/2})$.
The Tomita algebra $\mathcal T_\varphi$ consists of the entire analytic
elements $a$ for which $\sigma_z^\varphi(a)\in
\mathfrak N_\varphi\cap\mathfrak N_\varphi^*$ for every $z\in\bC$.

Let $\cM_*$ be the predual of $\cM$, i.e. the space of all normal bounded linear functionals on $\cM$ and $\cM_*^+$ the space of all positive normal linear functionals in $\cM_*$. 
The spaces $\mathcal L_\varphi$ and $\mathcal R_\varphi$ are
\begin{align*}
\mathcal L_\varphi
 &=\{x\in\mathfrak N_\varphi:
 y^*\mapsto\varphi(y^*x)\text{ extends to an element of }\cM_*\},\\
\mathcal R_\varphi
 &=\{x\in\mathfrak N_\varphi^*:
 y\mapsto\varphi(xy)\text{ extends to an element of }\cM_*\}.
\end{align*}
The extensions are denoted by $x\varphi$ and $\varphi x$, respectively.
Then $\mathcal R_\varphi=\mathcal L_\varphi^*$ and $\mathcal T_\varphi^2\subset\mathcal L_\varphi\cap\mathcal R_\varphi$.
Here and below, $\mathcal T_\varphi^2$ may be replaced by its linear span.
The compatible structure is Izumi's structure at the parameter $\displaystyle -\frac{1}{2}$ \cite{Izu97}:
\[
 L^p(\cM)=(\cM,\cM_*)_{[1/p]},\qquad 1\le p\le\infty.
\]
Its common ambient space is $\mathcal R_\varphi^*$, using the pairings $\langle\cdot, \cdot\rangle$:
\[
 \langle\omega,y\rangle=\omega(y),\qquad
 \langle x,y\rangle=(\varphi y)(x).
\]
The embedding of $\mathcal L_\varphi$ into $L^p(\cM)$ is $\iota^p$.
All intersections of different interpolated $L^p$ spaces are understood in this fixed ambient space.

Let $\phi$ be a normal semifinite faithful weight on the commutant $\cM'$ and set $d_\varphi=d\varphi/d\phi$ the spatial derivative. (\cite{Con80,Hil81, Ter81} for more details.)  
We write $[x]$ for the closure of a closable operator $x$.  
The spatial space $L^p(\phi)$ and its integral $\displaystyle \int \cdot d\phi$ are in the sense of Connes and Hilsum \cite{Hil81}.  
If $x=w_x|x|\in L^p(\phi)$ is the polar decomposition, then
$w_x\in\cM$ and
\[
 |x|^p=\frac{d\omega}{d\phi},\quad \omega\in\cM_*^+,
 \qquad\|x\|_{p,\phi}^p=\omega(1),\quad 1\le p<\infty.
\]
In particular the polar decompositions used below are the spatial polar decompositions. 
We shall denote the norm $\|\cdot\|_{p, \phi}$ by $\|\cdot\|_p$ briefly if there is no confusion.
Products of spatial $L^p$ space operators mean their closed strong products.  
For $1/r=1/p+1/s\le1$ they satisfy the H\"{o}lder inequality:
\[
 \|xy\|_{r,\phi}\le\|x\|_{p,\phi}\|y\|_{s,\phi}.
\]
For a vector $\xi\in\cH_\varphi$, define initially $R^\phi(\xi)\Lambda_\phi(a')=a'\xi$ for $a'\in\mathfrak N_\phi\subset\cM'$. 
The vector is $\phi$-bounded when this map has a bounded extension.
Their space is denoted by $\mathcal D(\cH_\varphi,\phi)$.
The spatial integral $\displaystyle \int \cdot d\phi$ is cyclic on products whose reciprocal exponents sum to one.

We use the following precise comparison facts, rather than a formal identification of the same unbounded operator in different spaces.
By \cite[Proposition 2.21]{Cas13}, there are isometric isomorphisms
\begin{equation}\label{eq:Phi-correct-core}
 \Phi_p:L^p(\phi)\longrightarrow L^p(\cM),\qquad
 \Phi_p([a d_\varphi^{1/p}])=\iota^p(a),
 \quad a\in\mathcal T_\varphi^2.
\end{equation}
For $a\in\mathfrak N_\varphi$ and $2\le p\le\infty$, the operator $a d_\varphi^{1/p}$ is closable and its closure lies in $L^p(\phi)$;
in particular, $\Phi_2([a d_\varphi^{1/2}])=\Lambda_\varphi(a).$
By \cite[Theorems 3.3 and 3.5]{Cas13}, we have that 
\begin{equation}\label{eq:intersection-correct-endpoints}
 \mathcal I_\varphi=L^1(\cM)\cap L^2(\cM),\qquad
 \mathfrak N_\varphi=L^2(\cM)\cap L^\infty(\cM),
\end{equation}
where
\[
 \mathcal I_\varphi=
 \{\omega\in\cM_*:\Lambda_\varphi(a)\mapsto\omega(a^*)
       \text{ is bounded on }\Lambda_\varphi(\mathfrak N_\varphi)\}.
\]
The identifications in \eqref{eq:intersection-correct-endpoints} are
$\omega\mapsto\xi(\omega)$ and $a\mapsto\Lambda_\varphi(a)$.
The associated embeddings into intermediate spaces are denoted by
$\xi_p$ and $\iota^s$ in the rest of the paper.

\begin{lemma}\label{lem:graph-identification}
Suppose $1\le p,s<\infty$ and $a_\alpha\in\mathcal T_\varphi^2$ satisfy
\[
 [a_\alpha d_\varphi^{1/p}]\longrightarrow x \quad\text{in }L^p(\phi),
 \qquad
 [a_\alpha d_\varphi^{1/s}]\longrightarrow y\quad\text{in }L^s(\phi).
\]
Then $\Phi_p(x)=\Phi_s(y)$ in $\mathcal R_\varphi^*$, so this common
element belongs to $L^p(\cM)\cap L^s(\cM)$.
\end{lemma}
\begin{proof}
The continuous embeddings into $\mathcal R_\varphi^*$ give the same
limit to the net $a_\alpha$.  The ambient space is Hausdorff, hence its
two limits coincide.
\end{proof}

Suppose that $a\in \cM$ is analytic with respect to $\sigma^\varphi$.
Then we have that $a d_\varphi^z=d_\varphi^z\sigma^\varphi_{iz}(a).$
Beyond that analyticity, an identity involving $x d_\varphi^z$ requires a specified closable product and a graph argument, such as
Lemma~\ref{lem:graph-identification}.  

In the following, we shall recall the Stein interpolation theorem \cite{Stein56}.
\begin{proposition}[Stein interpolation]\label{prop:steininter}
Let $(\mathcal X_0,\mathcal  X_1)$ and $(\mathcal Y_0,\mathcal Y_1)$ be compatible Banach couples, and let
$T_z$ be an admissible analytic family on a common dense subspace for the strip $0\leq\Re z\leq1$. 
Suppose that the boundary operators satisfy, uniformly in $t\in\mathbb R$,
\[
 \|T_{it}\|_{\mathcal  X_0\to \mathcal  Y_0}\le M_0,\qquad
 \|T_{1+it}\|_{\mathcal  X_1\to \mathcal  Y_1}\le M_1.
\]
Assume the usual admissible growth bound on scalar pairings. 
Then, for $0<\theta<1$, $T_\theta$ extends boundedly from $(\mathcal  X_0, \mathcal X_1)_{[\theta]}$ to $(\mathcal Y_0,\mathcal Y_1)_{[\theta]}$, with norm at most
$M_0^{1-\theta}M_1^\theta$. 
The corresponding multilinear form follows by the multilinear Stein interpolation theorem.
\end{proposition}

We will also use the Clarkson--McCarthy inequality \cite{Clarkson36, Mc67}, valid for these noncommutative $L^s$ spaces when $s\ge2$:
\begin{equation}\label{near:Clarkson}
 \left\|\frac{x+y}{2}\right\|_s^s +\left\|\frac{x-y}{2}\right\|_s^s
\le\frac{\|x\|_s^s+\|y\|_s^s}{2}.
\end{equation}
Its extension to general spatial or Haagerup $L^s$ spaces follows from the tracial reduction theorem; see \cite[Theorem 3.1 and
Remark 3.3]{HJX}.

For $1\leq r\leq \infty$, we denote by $r'$ the H\"{o}lder conjugate of $r$.
For a nonzero spatial element $x=w_x|x|\in L^r(\phi)$, $1<r<\infty$, put
\[
 N_r(x)=\frac{w_x|x|^{r-1}}{\|x\|_r^{r-1}}\in L^{r'}(\phi).
\]
Then $\|N_r(x)\|_{r'}=1$ and $\displaystyle \int x N_r(x)^* d\phi=\|x\|_r$.
Moreover, we have that $N_{r'}(N_{r}(x))=x$ if $\|x\|_r=1$.

We shall use two facts about spatial $L^r$ spaces explicitly.
First, if $2\leq r\leq\infty$, $x\in L^r(\phi)$ and
$\xi\in\mathcal D(\mathcal H_\varphi,\phi)$ is a $\phi$-bounded vector, then
\begin{equation}\label{repair:bounded-vector}
 \|x\xi\|\leq \|x\|_{r,\phi}
 \|\xi\|^{1-2/r}\|R^\phi(\xi)\|^{2/r}.
\end{equation}
Second, if $x\in L^1(\phi)^+$ and $a\in\mathcal T_\varphi^2$, then
\begin{equation}\label{repair:L2path}
 z\longmapsto x^{z/2}a d_\varphi^{(1-z)/2}
\end{equation}
is an $L^2(\phi)$-valued analytic function for $0<\Re z<1$,
continuous on the closed strip, where $x^0$ is treated as the support projection of $x$.
The corresponding assertion with the
order reversed follows by taking adjoints.  
These are the spatial bounded-vector estimate and the analytic-section lemma, respectively \cite{Hil81,Ter82}.

\begin{lemma}\label{lemma:simultaneous-core}
Let $1\le p\le2$ and let $\omega$ be a common element of $L^p(\bG)\cap L^2(\bG)$.
There is a net $\omega_\alpha \in\mathcal I_\varphi$ such that
\[
 \xi_p(\omega_\alpha)\longrightarrow\xi_p(\omega) \text{ in }L^p(\bG),\qquad
 \xi(\omega_\alpha)\longrightarrow\xi(\omega) \text{ in }L^2(\bG).
\]
\end{lemma}
\begin{proof}
Choose a bounded Tomita approximate identity $a_\alpha$ such that
$\sigma_z^\varphi(a_\alpha)\to1$ strongly$^*$, with a uniform
operator-norm bound for each fixed $z$. 
For an analytic element $x$, the maps
\[
 R_x^{(1)}(\eta)=\eta\sigma_{i/2}^\varphi(x) \quad(\eta\in\cM_*),\qquad
 R_x^{(\infty)}(y)=y\sigma_{-i/2}^\varphi(x) \quad(y \in\cM)
\]
are compatible bounded maps. 
Indeed, for a Tomita-square element $b$ the KMS identity gives
\[
 ((b\varphi)\sigma_{i/2}^\varphi(x))(y)
 =\varphi(\sigma_{i/2}^\varphi(x)yb)
 =\varphi(yb\sigma_{-i/2}^\varphi(x)).
\]
The equality first holds on the common analytic domain and then for
the normal extensions. It proves agreement in the compatible ambient
space. 
Interpolation gives a map $R_x^{(s)}$ at every intermediate exponent $1\leq s<\infty$.
Its spatial formula, checked on the Tomita-square core, is
\begin{equation}\label{eq:common-core-cutoff}
 \Phi_s^{-1}(R_x^{(s)} y)
 =\Phi_s^{-1}(y)\sigma_{i(1/s-1/2)}^\varphi(x),
 \qquad1\le s<\infty.
\end{equation}
For example, if $y=\iota^s(b)$, the right side is $b\sigma_{-i/2}^\varphi(x)d_\varphi^{1/s}$.
The bounded module actions and core density extend this formula.

Put $x_2=\Phi_2^{-1}(\xi(\omega))$.
By the spatial H\"older inequality, we see that $x_2a_\alpha d_\varphi^{1/2}\in L^1(\phi).$
Let $\omega_\alpha$ be the corresponding normal functional. 
For $b\in\mathfrak N_\varphi$, the spatial pairing gives
\[
 \omega_\alpha(b^*)
 =\int x_2a_\alpha d_\varphi^{1/2}b^*\,d\phi
 =\langle\Phi_2(x_2a_\alpha),\Lamv(b)\rangle.
\]
Thus $\omega_\alpha\in\mathcal I_\varphi$, and its $L^2$ representative is $\Phi_2(x_2a_\alpha)=R_{a_\alpha}\xi(\omega)$.
The common-space identifications then show that its $L^p$ representative is $R_{a_\alpha}\xi_p(\omega)$ as well.
Formula~\eqref{eq:common-core-cutoff} and normality of the finite spatial module actions give convergence in both norms: 
the right multipliers are $a_\alpha$ for $s=2$ and $\sigma_{i(1/p-1/2)}^\varphi(a_\alpha)$ for $s=p$; 
both are uniformly bounded and converge strongly$^*$ to $1$.
Explicitly, if $1\le s\le2$, $x=w_x|x|\in L^s(\phi)$ and $b_\alpha\to1$ strongly$^*$ with a uniform norm bound, then
\[
 \|x(b_\alpha-1)\|_{s,\phi}
 \le\bigl\||x|^{1-s/2}\bigr\|_{2s/(2-s),\phi} \bigl\||x|^{s/2}(b_\alpha-1)\bigr\|_{2,\phi}
 \longrightarrow0.
\]
At $s=2$ the first factor is the support projection in operator norm.
The convergence of the second factor is the normal right module action on spatial $L^2$. 
This proves the needed norm convergence.
\end{proof}

\subsection{Equality in spatial H\"older inequalities}

We recall one elementary consequence of the strict convexity of
noncommutative $L^p$ spaces; this holds for arbitrary von Neumann
algebras, via the usual spatial/Haagerup identifications
\cite{Ter81,Hil81}.
\begin{lemma}\label{lem:norming-spatial}
Suppose that $1<p<\infty$, $x=w_x|x|\in L^p(\phi)$, and $\|x\|_{p,\phi}=1$.
If $y\in L^{p'}(\phi)$ satisfies $\|y\|_{p',\phi}\le1$ and $\displaystyle \int yx\,d\phi=1$, then $y=N_p(x)^*=|x|^{p-1}w_x^*$.
\end{lemma}
\begin{proof}
It follows from the strict convexity of $L^{p'}(\phi)$.
Indeed, both $y$ and $N_p(x)^*$ have norm one. If they were different, strict convexity would imply $\|y+N_p(x)^*\|_{p'}<2$.
However, 
\begin{align*}
\|y+N_p(x)^*\|_{p'}\geq \int yx d\phi +\int N_p(x)^* x d\phi=2.
\end{align*}
This is a contradiction.
\end{proof}

\begin{proposition}\label{prop:holder-repaired}
Let $1<p,s<\infty$, $1\le r<\infty$, and $1/r=1/p+1/s$.  
For nonzero $a\in L^p(\phi)$ and $b\in L^s(\phi)$,
\begin{equation}\label{eq:holder-equality-repaired}
 \|ab\|_{r,\phi}=\|a\|_{p,\phi}\|b\|_{s,\phi}
 \quad \text{ if and only if }\quad
 \frac{|a|^p}{\|a\|_{p,\phi}^p}
 =\frac{|b^*|^s}{\|b\|_{s,\phi}^s}.
\end{equation}
\end{proposition}
\begin{proof}
Assume that $\|a\|_{p,\phi}=\|b\|_{s,\phi}=1$ and write $a=w_a|a|$, $b=w_b|b|=|b^*|w_b$, $ab=w_{ab}|ab|$.
Suppose first $\|ab\|_{r,\phi}=1$.  
Set \[
 x=\begin{cases}
 |ab|^{r-1}w_{ab}^*,&r>1,\\
 w_{ab}^*,&r=1.
 \end{cases}
\]
For $r>1$ this belongs to $L^{r'}(\phi)$ with norm one; at $r=1$ it belongs to $\cM$ with norm one.  
In both cases, we have $\displaystyle \int x ab\,d\phi=1$.
Since $\displaystyle \frac1{p'}=\frac1s+\frac1{r'},\quad
 \frac1{s'}=\frac1{r'}+\frac1p,$ H\"older inequality gives $\|bx\|_{p',\phi}\le 1$ and $\|x a\|_{s',\phi}\le1$.  
Lemma~\ref{lem:norming-spatial} therefore yields
\[
 bx=|a|^{p-1}w_a^*,\qquad xa=|b|^{s-1}w_b^*.
\]
Multiplying the first identity on the right by $a$ and the second on
the left by $b$, in the associative spatial product, gives
\[
 |a|^p=b x a=w_b|b|^sw_b^*=|b^*|^s.
\]

Conversely, suppose $|a|^p=|b^*|^s$.  
Let $p_a$ be the support of $|a|$.  
Then $w_a^*w_a=p_a=w_bw_b^*$ and
\[
 ab=w_a |a|^{p/r}w_b,\qquad |ab|^r=w_b^*|a|^p w_b.
\]
Consequently
\[
 \|ab\|_{r,\phi}^r=\int w_b^*|a|^p w_b\,d\phi
 =\int |a|^p w_bw_b^*\,d\phi=\int |a|^p\,d\phi=1.
\]
This proves \eqref{eq:holder-equality-repaired}.
\end{proof}

\begin{proposition}\label{prop:holder-endpoint-repaired}
For $a\in\cM$ and $b\in L^1(\phi)$, $\|ab\|_{1,\phi}=\|a\|_\infty\|b\|_{1,\phi}$ if and only if $\mathcal R(b)$ is dominated by the spectral projection of $|a|$ at $\|a\|_\infty$.
\end{proposition}
\begin{proof}
If one of $a, b$ is zero, we see the proposition holds.
We shall assume that $\|a\|_\infty=\|b\|_{1, \phi}=1$.
Write $b=|b^*|w_b$, and put $h=|b^*|^{1/2}$.  
H\"older inequality and Cauchy-Schwarz inequality give
\[
 \|ab\|_{1,\phi}\le\|ah\|_{2,\phi}\|hw_b\|_{2,\phi}
 =\left(\int h a^*ah\,d\phi\right)^{1/2}\le1.
\]
If the equality holds, then $\displaystyle \int h(1-a^*a)h\,d\phi=0$.  
Thus $(1-a^*a)^{1/2}h=0$, and the support of $h$, namely
$\mathcal R(b)$, is contained in the spectral projection of $|a|$ at one. 

Conversely that support condition gives $b^*a^*ab=b^*b$ and therefore $|ab|=|b|$.
\end{proof}

\subsection{Locally compact quantum groups}

A locally compact quantum group is a quadruple $\bG=(\cM,\Delta,\varphi,\psi)$, where $\cM$ is a von Neumann algebra, 
$\Delta:\cM\to\cM\overline\otimes\cM$ is a normal unital $*$-homomorphism satisfying coassociativity, and $\varphi,\psi$ are normal semifinite faithful left and right Haar weights. 
Thus, for positive normal $\omega$ and positive Haar-integrable $a$,
\[
 \varphi((\omega\otimes\iota)\Delta(a))=\omega(1)\varphi(a),\qquad
 \psi((\iota\otimes\omega)\Delta(a))=\omega(1)\psi(a).
\]
We use the framework of \cite{KusVae00,KusVae03}. 
The multiplicative unitary $W_\varphi$ is characterized by
\[
 W_\varphi^*(\Lamv(a)\otimes\Lamv(b))=(\Lamv\otimes\Lamv)(\Delta(b)(a\otimes1)),
 \qquad a,b\in\mathfrak N_\varphi.
\]
In particular $\Delta(a)=W_\varphi^*(1\otimes a)W_\varphi$.
The antipode is the generally unbounded closed map $S=R\tau_{-i/2}$, where $R$ is the unitary antipode and $\tau$ the scaling automorphism group. 
We normalize the right Haar weight by $\psi=\varphi R$. 
The modular element $\delta$ is positive and nonsingular, affiliated with $\cM$, and satisfies
\[
 \psi=\varphi_\delta,\quad\Delta(\delta)=\delta\otimes\delta,
 \quad R(\delta)=\delta^{-1},\quad\tau_t(\delta)=\delta.
\]
The notation $\varphi_\delta$ denotes the standard weight perturbation; on its natural domain it is $\varphi(\delta^{1/2}\,\cdot\,\delta^{1/2})$.

The dual von Neumann algebra $\hcM$ is generated by the slices $\lambda(\omega)=(\omega\otimes\iota)(W_\varphi)$, and its comultiplication is implemented by $\widehat W=\Sigma W_\varphi^*\Sigma$, where $\Sigma$ is the flip of the tensor.
The dual Haar weight $\hat\varphi$ is normalized by $\hLam(\lambda(\omega))=\xi(\omega)$ for $\omega\in\mathcal I_\varphi$.
We write $\hat J,\hnabla,\hat\sigma$ for its modular data and $\hat\lambda$ for the dual Fourier representation.
From now on the spatial reference weights are
\[
 \phi=J_\varphi\varphi J_\varphi,\qquad
 \hat\phi=\hat J\hat\varphi\hat J,
 \qquad d_\varphi=\nabla_\varphi,\quad\hd=\hnabla.
\]

For a closed operator $x$, $\mathcal R(x)$ denotes its left range projection, the projection onto $\overline{\operatorname{Ran}x}$.
For a spatial element $x$, its right support is $\mathcal R(x^*)=\mathcal R(|x|)$. 
These projections belong to $\cM$.
We denote
\[
 \cS(x)=\varphi(\mathcal R(|x|)),\qquad x\in L^p(\phi).
\]
For $x\in L^p(\bG)$ this means $\cS(x)=\cS(\Phi_p^{-1}(x))$.
The modular centralizers are
\[
 \mathcal C_\varphi=\{x\in\cM:\sigma_t^\varphi(x)=x\ (t\in\bR)\},
 \qquad
 \mathcal C_{\hat\varphi}
 =\{y\in\hcM:\hat\sigma_t(y)=y\ (t\in\bR)\}.
\]
A projection in the first centralizer strongly commutes with $d_\varphi$, and the dual statement holds for $\hd$.

On normal functionals the isometric strongly continuous group $\rho$ is
\[
 \rho_t(\omega)(x)=\omega(\delta^{-it}\tau_{-t}(x)),
 \qquad \hat\sigma_t(\lambda(\omega))=\lambda(\rho_t(\omega)).
\]
We shall use the same symbol for the spatial spaces.

The scaling constant is denoted by $\nu$ and we have that
\[
 \sigma_t^\varphi(\delta)=\nu^t\delta,\quad
 \varphi\tau_t=\nu^{-t}\varphi,
 \quad\sigma_s^\varphi\tau_t=\tau_t\sigma_s^\varphi.
\]
The standard right-Haar cocycle reads
\[
 (D\psi:D\varphi)_t=\nu^{it^2/2}\delta^{it}.
\]
We have that 
\[
 d_\psi^{1/p}=\nu^{-ip^{-2}/2}[\delta^{1/p} d_\varphi^{1/p}].
\]

\subsection{The Fourier transform and Hausdorff-Young inequality}
The Fourier transform $\cF_p: L^p(\bG)\longrightarrow L^{p'}(\hbG)$, $1\le p\le2$, is obtained by interpolating the maps $\cF_1=\lambda$ and unitary transform $\cF_2:\cH_\varphi\to\cH_\varphi$ given by $\Lambda_\varphi(x)\to \widehat\Lambda(\lambda(x\varphi))$ with inverse the dual $L^2$ Fourier transform \cite[Corollary 5.5 and Theorem 5.6]{Cas13}.  
In the spatial form, we denote $\tcF_p=\widehat\Phi_{p'}^{-1}\cF_p\Phi_p.$
The formula used for computations is precisely
\begin{equation}\label{eq:Fourier-core-repaired}
 \tcF_p([x d_\varphi^{1/p}])=\left[\lambda(x\varphi)\hd^{1/p'}\right],\qquad
 x\in\mathcal T_\varphi^2,
\end{equation}
as in \cite[Theorem 5.7]{Cas13}.  
If an argument needs both $L^p$ and $L^2$ Fourier representatives, their compatibility is invoked on $L^p(\bG)\cap L^2(\bG)$.
The dual Fourier transforms are denoted by $\widehat\cF_p$ and $\widehat\tcF_p$ for the spatial spaces.

The Hausdorff--Young inequalities are
\[
 \left\|\tcF_p(x)\right\|_{p',\hat\phi}\le\|x\|_{p,\phi},\qquad
 \left\|\widehat{\tcF}_p(y)\right\|_{p',\phi}\le\|y\|_{p,\hat\phi},
 \qquad1\le p\le2.
\]
We denote by $\displaystyle H_p(\bG)=\sup_{\|x\|_{p,\phi}=1}\left\|\tcF_p(x)\right\|_{p',\hat\phi}=\|\tcF_p\|_{p\to p'}$.
We shall use the sesquilinear adjointness identity
\begin{equation}\label{eq:rewritten-Fourier-adjointness}
 \int y\tcF_p(x)^*\,d\hat\phi=\int\widehat{\tcF}_p(y)x^*\,d\phi,
 \qquad x\in L^p(\phi),\quad y\in L^p(\hat\phi).
\end{equation}
H\"older inequality and the Hausdorff-Young inequality make both pairings continuous in the displayed $L^p$ norms. 
Density of the Tomita-square core therefore proves \eqref{eq:rewritten-Fourier-adjointness} for all displayed inputs.
For a bilinear pairing the resulting formula is
\[
 \int\tcF_p(x)b\,d\hat\phi=\int x\bigl(\widehat{\tcF}_p(b^*)\bigr)^*\,d\phi.
\]

The Fourier transform $\tcF_p$ is injective for $1\le p\le2$. 
The range for $1<p\le2$ is norm dense in $L^{p'}(\hat\phi)$; 
at $p=1$ it is $\sigma$-weakly dense in $\hcM$.
These are the usual Fourier uniqueness and density properties of the compatible construction \cite{KusVae03,Cas13}.

\begin{proposition} \label{prop:HY-norming-rewritten}
Let $1<p<2$ and let $0\ne x\in L^p(\phi)$ with $y=\tcF_p(x)\ne0$.
Then 
\[
 \|y\|_{p',\hat\phi}=\|x\|_{p,\phi},
\quad  \text{ if and only if }\quad
 \widehat{\tcF}_p(N_{p'}(y))=N_p(x).
\]
\end{proposition}
\begin{proof}
If the norms agree, the Fourier adjointness gives
\[
 \int\widehat{\tcF}_p(N_{p'}(y))x^*\,d\phi
  =\int N_{p'}(y)y^*\,d\hat\phi
  =\|y\|_{p',\hat\phi}=\|x\|_{p,\phi}.
\]
The inverse Fourier image on the left has $L^{p'}$-norm at most one.
Uniqueness of the norming element, as in Lemma~\ref{lem:norming-spatial}, identifies it with $N_p(x)$.
Conversely the proposed identity and the same pairing equate the two norms. 
\end{proof}

We write $\hcS(x)=\cS(\tcF_p(\Phi_p^{-1}(x)))$ with respect to the dual Haar weight $\hat\varphi$, the index $p$ is fixed by the declared space of $x$.

\subsection{A joint analytic core with graph control}
Put
\[
 \mathcal E_\varphi=\mathfrak N_\varphi\cap\mathfrak N_\varphi^*,
 \qquad
 \|x\|_{\mathcal E_\varphi}
 =\|x\|_\infty+\|\Lambda_\varphi(x)\|
                  +\|\Lambda_\varphi(x^*)\|.
\]
This is a Banach $*$-algebra.  
Define
$\mathcal T_{\varphi,\delta,\tau}$ to consist of those $x$ for which
\begin{equation}\label{eq:joint-translate-graph}
 (z_1,z_2,z_3,z_4)\longmapsto
 \delta^{z_4}\tau_{z_2}(\sigma_{z_1}^\varphi(x))\delta^{z_3}
\end{equation}
extends to an entire $\mathcal E_\varphi$-valued function.  On the
imaginary $\delta$ axes its meaning is ordinary bounded multiplication;
the other values are the unique analytic extensions.  In particular,
boundedness and both GNS domain conditions are part of this definition.

\begin{lemma}\label{lem:joint-graph-core}
For $x\in\mathcal E_\varphi$ and $n\ge1$, set
\begin{align}\label{eq:gaussian-graph-core}
 x_n=\frac{n^4}{\pi^2}\int_{\bR^4}
 e^{-n^2(t_1^2+t_2^2+t_3^2+t_4^2)}
 \delta^{it_4}\tau_{t_2}(\sigma_{t_1}^\varphi(x))\delta^{it_3}\,dt_1dt_2dt_3dt_4.
\end{align}
The integral is interpreted weakly in $\cM$ and in norm after applying $\Lambda_\varphi$.  
Then $x_n\in\mathcal T_{\varphi,\delta,\tau}$,
$\|x_n\|_\infty\le\|x\|_\infty$, and
\[
 x_n\longrightarrow x\text{ strongly}^*,\qquad
 \Lambda_\varphi(x_n)\longrightarrow\Lambda_\varphi(x),\quad
 \Lambda_\varphi(x_n^*)\longrightarrow\Lambda_\varphi(x^*).
\]
If $x\in\mathcal T_\varphi$, these convergences also hold after
applying any fixed $\sigma_z^\varphi$, and the resulting operator
norms are bounded uniformly in $n\ge1$.
\end{lemma}
\begin{proof}
The four real actions in \eqref{eq:gaussian-graph-core} preserve
$\mathcal E_\varphi$.  
We have that 
\begin{align*}
 \|\Lambda_\varphi(\sigma_t^\varphi(x))\|=\|\Lambda_\varphi(x)\|, & \quad
 \|\Lambda_\varphi(\tau_t(x))\|=\nu^{-t/2}\|\Lambda_\varphi(x)\|,\\
 \|\Lambda_\varphi(\delta^{it}x)\|=\|\Lambda_\varphi(x)\|,& \quad 
 \|\Lambda_\varphi(x\delta^{it})\|=\nu^{-t/2}\|\Lambda_\varphi(x)\|.
\end{align*}
The last equality follows from \eqref{eq:correct-right-GNS} and $\sigma_{i/2}^\varphi(\delta^{it})=\nu^{-t/2}\delta^{it}$.
The corresponding estimates for adjoints have the same exponential
growth bound.  
Hence all graph integrals are absolutely convergent.
The real actions are continuous in both norms in $\cH_\varphi$ and in the strong$^*$ topology on bounded sets.  
Gaussian approximate identities give the three stated convergences.  
Closedness of the maps $\Lambda_\varphi$ identifies the weak operator integral with these graph integrals.

For completeness, analyticity can be checked one action at a time.
Conjugating the integrand by $\sigma_s^\varphi$ translates its $t_1$ variable and multiplies its coefficient by $\nu^{is(t_3+t_4)}$; 
applying $\tau_s$ translates $t_2$; 
the left and right multiplication by $\delta^{is}$ translate $t_4$ and $t_3$.
Thus any mixed complex action has an integral whose scalar coefficient is a translated Gaussian multiplied by an exponential linear in the real variables.  
On every compact set of complex parameters this coefficient, and each of its derivatives, is bounded by $\displaystyle C'\exp\left(-C\sum_jt_j^2+C'\sum_j|t_j|\right)$, with $C, C'>0$.
Together with the preceding graph estimates this proves entire analyticity in $\mathcal E_\varphi$. 
More explicitly, translated Gaussian kernel differences and their derivatives have the stated integrable majorants. 
Their integrals bound the operator-norm differences of the weak operator integrals and the norm differences of both images in $\cH_\varphi$. 
Differentiation is consequently in all three components of the graph norm; 
no Bochner integral in the unsmoothed operator-norm orbit is being assumed.
This argument also verifies that the analytic values are the closures of the indicated unbounded products.

If $x\in\mathcal T_\varphi$, move a fixed $\sigma_z^\varphi$ through the real actions, applying it to $x$ and retaining the scalar $\nu^{iz(t_3+t_4)}$. 
The resulting Gaussian approximate identities converge to $\sigma_z^\varphi(x)$ in the same three senses.
Completing the square bounds their operator norms uniformly for $n\ge1$.
\end{proof}

\begin{corollary}\label{cor:graph-Lp-density}
The space $\mathcal T_{\varphi,\delta,\tau}$ is a $*$-algebra contained in $\mathcal T_\varphi\cap\mathfrak N_\psi\cap\mathfrak N_\psi^*$.
For $1\le p<\infty$,
\[
 \iota^p(\operatorname{span}\mathcal T_{\varphi,\delta,\tau}^2)
 \text{ is norm dense in }L^p(\bG).
\]
\end{corollary}
\begin{proof}
Products and adjoints preserve the graph core: apply each mixed
translate, distribute the bounded factors between the two terms, and use $\|xy\|_{\mathcal E_\varphi} \le\|x\|_{\mathcal E_\varphi}\|y\|_{\mathcal E_\varphi}$.
The membership in $\mathfrak N_\psi$ follows from $x\delta^{1/2}\in\mathfrak N_\varphi$ and $\psi(x^*x)=\|\Lambda_\varphi(x\delta^{1/2})\|^2$;
apply this also to $x^*$.

The standard Tomita square core is norm dense after embedding into
each finite $L^p$ space.  
To approximate $ab$ with $a,b\in
\mathcal T_\varphi$, apply Lemma~\ref{lem:joint-graph-core} to each
factor, giving $a_n,b_n$ in the joint graph core.  The products are
uniformly bounded and converge strongly$^*$ to $ab$.  Moreover
\[
 [a_n b_n d_\varphi]
 =[a_n d_\varphi^{1/2}]
   [\sigma_{i/2}^\varphi(b_n)d_\varphi^{1/2}]
 \longrightarrow
 [a d_\varphi^{1/2}]
   [\sigma_{i/2}^\varphi(b)d_\varphi^{1/2}]
 =[ab d_\varphi]
\]
in $L^1(\phi)$, by the H\"older inequality.
For $z_n=a_nb_n-ab$, interpolation theorem gives
\[
 \|\iota^p(z_n)\|_p\le\|z_n\|_\infty^{1-1/p}\|z_n\varphi\|^{1/p}
 \longrightarrow 0.
\]
This proves the asserted density.
\end{proof}

\begin{remark}
The preceding lemma is an elementwise approximation, and suffices for the norm-density arguments. 
Independently, the usual Tomita approximate identity provides contractions $a_\alpha\in\mathcal T_\varphi$ with $\sigma_z^\varphi(a_\alpha)\to1$ strongly$^*$ and uniform bounds at each fixed $z$ \cite[Proposition 2.16]{Str81}.
Whenever the argument uses only $\sigma^\varphi$ it is this net that is used.  
No approximate identity is asserted to satisfy $\delta^{z_4}b_\beta\delta^{z_3}\to1$ for arbitrary $z_3,z_4$.
For bounded imaginary powers the correct limit is $\delta^{i(t_3+t_4)}$ if $b_\beta\to1$ strongly$^*$.
Positive or negative real powers require a specified graph topology.
\end{remark}

\subsection{Group-like projections}

A group-like projection is a nonzero projection $B\in\cM$ such that $\Delta(B)(1\otimes B)=B\otimes B.$
When $\varphi(B)<\infty$, it is an extremizer of the Hausdorff-Young inequality. 
Its shifts and their Fourier dual shifts lead to the
notion of a bi-shift. We fix the following equivalent analytic
convention from \cite{JLW18} so that this term has one meaning
throughout the paper.

\begin{definition}\label{def:bishift-rewritten}
A nonzero $\omega\in L^1(\bG)\cap L^2(\bG)$ is a Donoho--Stark minimizer if
\[
 \cS(\xi_p(\omega))\,\cS(\iota^s(\lambda(\omega)))=1
 \quad\text{for every }1\le p\le2,\ 2\le s\le\infty.
\]
We call these functionals bi-shifts of group-like projections, using
the structural equivalence proved in \cite{JLW18}. The same term for
a spatial element means a nonzero scalar multiple of its appropriate
spatial embedding.
\end{definition}
A minimizer can be written $\omega=\mu w\varphi$, where $\mu>0$, $w$ is a partial isometry, $p_{\omega}=|w|\in\mathcal C_\varphi$, and
$\varphi(p_{\omega})<\infty$. 
The Fourier transform $\lambda(\omega)$ is a multiple of a partial isometry and $|\lambda(\omega)|\in \mathcal{C}_{\hat{\varphi}}$.
Precisely, $\mu^{-1}\varphi(p_{\omega})^{-1} \lambda(\omega)$ is the partial isometry.
Moreover, we have that $\tau_t(p_\omega)=p_{\omega}$ and $\hat{\tau}_t(|\lambda(\omega)|)=|\lambda(\omega)|$ for all $t\in \bR.$
A bishift $\omega$ of group-like projection is scaling-invariant if $\tau_t(w)=w$ for all $t\in \bR$.

\begin{proposition}\label{prop:bishift-adjoint-criterion}
Let $\omega=w\varphi$ be a bi-shift of a group-like projection, with $w$ a nonzero partial isometry.  
For $1\le p<\infty$, the following are equivalent:
\begin{enumerate}
\item There is a scalar $\tilde\kappa$ such that
\begin{equation}\label{eq:bishift-adjoint-criterion}
 [w d_\varphi^{1/p}]^*=\tilde\kappa [w^*d_\varphi^{1/p}].
\end{equation}
\item There is a number $\kappa>0$ such that
\begin{equation}\label{eq:bishift-modular-eigen}
 \sigma_t^\varphi(w)=\kappa^{it}w\qquad(t\in\mathbb R).
\end{equation}
\end{enumerate}
In this case $\tilde\kappa=\kappa^{-1/p}>0$, $ww^*$ belongs to the centralizer of $\varphi$, and $\displaystyle \kappa=\frac{\varphi(ww^*)}{\varphi(w^*w)}$.
At $p=\infty$ it holds trivially with $\tilde\kappa=1$.

For $1<p<\infty$, we have that
\begin{equation}\label{eq:bishift-norming-correct}
 N_p([w d_\varphi^{1/p}])
 =\varphi(w^*w)^{-1/p'}[w d_\varphi^{1/p'}].
\end{equation}
\end{proposition}
\begin{proof}
We have that  $w^*w \in\mathcal C_\varphi$ and $0<\varphi(w^*w)<\infty$; see \cite[Proposition 4.7]{JLW18}.
Put $x=[w d_\varphi^{1/p}]$. 
Since $w^*w$ strongly commutes with $d_\varphi$, $|x|= w^*w d_\varphi^{1/p}$ and $\|x\|_p^p=\varphi(w^*w)$. 
Consequently \eqref{eq:bishift-norming-correct} follows from the polar decomposition directly.

Suppose \eqref{eq:bishift-adjoint-criterion} holds. 
In particular the two closed products are equal and $\tilde\kappa\ne0$. 
Let $E$ be the spectral measure of $d_\varphi$. 
Cutting the equality \eqref{eq:bishift-adjoint-criterion} with the spectral projections on the left and right yields, for Borel sets $A,B\subset(0,\infty)$ bounded away from zero and infinity,
\[
 E(A)d_\varphi^{1/p} w^*E(B)=\tilde\kappa E(A)w^*d_\varphi^{1/p} E(B).
\]
The spectral theorem for the self-adjoint (unbounded) operator $d_\varphi^{1/p}$ implies that $\tilde\kappa>0$ and that $E(A)w^*E(B)=0$ whenever $A\cap \tilde\kappa^{p}B=\emptyset$.  
It follows by functional calculus that
\[
 d_\varphi^{it}w^*d_\varphi^{-it}=\tilde\kappa^{itp}w^*,\qquad t\in\mathbb R.
\]
Taking adjoints gives \eqref{eq:bishift-modular-eigen} with $\kappa=\tilde\kappa^{-p}$.

Conversely, \eqref{eq:bishift-modular-eigen} makes $w,w^*$ entire analytic and
\[
 d_\varphi^{1/p} w^*=\sigma_{-i/p}^\varphi(w^*)d_\varphi^{1/p}=\kappa^{-1/p}w^*d_\varphi^{1/p}
\]
on their natural closed-product domains. 
Since $[w d_\varphi^{1/p}]^*=[d_\varphi^{1/p} w^*]$, this proves \eqref{eq:bishift-adjoint-criterion} with $\tilde\kappa=\kappa^{-1/p}$.

Under the condition \eqref{eq:bishift-modular-eigen}, we see that $ww^*$ is fixed by $\sigma^\varphi$.
The KMS relation implies that 
\[
  \varphi(ww^*)=\varphi(w^*\sigma_{-i}^\varphi(w))
 =\kappa\varphi(w^*w),
\]
where finiteness follows from the same identity. 
This shows $\displaystyle \kappa=\frac{\varphi(ww^*)}{\varphi(w^*w)}$.
\end{proof}

\section{Scaling Fourier Transform}

In this section, we shall extend the ordinary Fourier transform to scaling Fourier transforms.
For $1\leq p\leq 2$, let $\displaystyle \kappa_p=\frac{1}{p}-\frac{1}{2}$.
On the analytic core, we define $T_p(z)$ and $\widehat T_r(z)$ as follows:
\[
 T_p(z)[x d_\varphi^{1/p}]=\nu^{z/p} [\tau_z(x)d_\varphi^{1/p}],\qquad
 \widehat T_r(z)[y \hd^{1/r}]=\nu^{-z/r}[\hat\tau_z(y)\hd^{1/r}].
\]
The operators $T_p(z)$, $\widehat T_r(z)$ extend to operators on $L^p(\phi)$, $L^r(\hat{\phi})$ whenever they are bounded.

\begin{definition}[Scaling Fourier Transform]
Suppose that $\bG$ is a locally compact quantum group, $1\leq p\leq 2$ and $0\leq \theta\leq 1$.
The scaling Fourier transform on $\bG$ is defined as follows
\[
 \tcF_p^{(\theta)}([x d_\varphi^{1/p}])
 =\nu^{-i\kappa_p\theta/p'} [\hat\tau_{i\kappa_p\theta}(\lambda(x\varphi))\hd^{1/p'}],
\]
where $x$ is in the corresponding analytic core.
Thus $\tcF_p^{(\theta)}=\widehat T_{p'}(i\kappa_p\theta)\tcF_p$ on the analytic core.
Denote the norm of the scaling Fourier transform by
\begin{align*}
    H_{p,\theta}(\bG) =\left\|\tcF_p^{(\theta)}\right\|_{p\to p'}.
\end{align*}
\end{definition}

\begin{remark}
When $\theta=0$, we have that $\tcF_p^{(0)}=\tcF_p$ and $H_{p, 0}(\bG)=H_p(\bG)$.
When the scaling group $\tau$ is trivial, the scaling Fourier transform $\tcF_p^{(\theta)}=\tcF_p$ is the usual Fourier transform and $H_{p, \theta}(\bG)=H_p(\bG)$.
\end{remark}

\begin{remark}
We shall denote the scaling Fourier transform $\cF_p^{(\theta)}=\widehat\Phi_{p'}\tcF_p^{(\theta)}\Phi_p^{-1}$.
\end{remark}

\begin{proposition}
We have that $H_{p, \theta}(\bG)\leq 1$.    
\end{proposition}
\begin{proof}
For $0\le\Re z\le1$ and $\omega\in \mathcal I_\varphi$ such that $\lambda(\omega)$ is analytic with respect to the scaling group $\hat\tau$, we define that 
\[
 S_z( \omega)=\exp\!\left(\frac{i\theta\log\nu}{8}z^2\right)
       \hat\tau_{i\theta(1-z)/2}(\lambda(\omega)).
\]
We have that  $\hat\tau_{i/2}(\lambda(\omega))=\lambda(\overline\omega R)^*$ with norm at most $\|\omega\|$. 
Hence
\begin{align*}
    \|S_z( \omega)\|\leq \exp\!\left(-\frac{\theta\log\nu}{8}\Im(z^2)\right) \|\omega\|_1.
\end{align*}
On $z=iy$, we see that $\|S_{iy}(\omega)\|\leq \|\omega\|_1$.
On $z=1+iy$, we have that 
\begin{align*}
\|S_{1+iy}(\omega)\|_2=\left\|\exp\!\left(-\frac{\theta y\log\nu}{4}\right)\hat\tau_{\theta y/2}(\lambda(\omega)) \right\|_2 
=\|\lambda(\omega)\|_2=\|\omega\|_2.
\end{align*}
Now we can apply the Stein interpolation theorem \ref{prop:steininter} at the real interpolation parameter $z=2/p'$ and obtain that 
\begin{align*}
    \|S_{2/p'}(\omega)\|_{p'} \leq \|\omega\|_{p}.
\end{align*}
By taking approximation, we have that  
\begin{align*}
    \|\tcF_p^{(\theta)}(x)\|_{p'}\leq \|x\|_p, \quad x\in L^p(\phi).
\end{align*}
This shows that $H_{p, \theta}(\bG)\leq 1$.
\end{proof}

We define an anti-isometry $J_p: L^p(\phi)\to L^{p}(\phi)$ for $1\leq p \leq \infty$ as follows:
\[
 J_p([x d_\varphi^{1/p}])
     =\nu^{-i/2p^{2}}[R(x^*)\delta^{1/p} d_\varphi^{1/p}].
\]
A core may be obtained using Gaussian functions of $\log\delta$ and joint Gaussian smoothing for the real actions.
We see that $J_p$ is a surjective conjugate-linear isometry and an involution.
When $p=1$, $J_1([x d_\varphi])=\nu^{-i/2}[R(x^*)\delta  d_\varphi].$
When $p=\infty$, $J_\infty(x)= R(x^*).$

\begin{proposition}\label{prop:symmetry}
For $x\in L^p(\phi)$, we have that 
\begin{equation}\label{nu:reflection-phase}
 \tcF_p^{(\theta)}(J_p(x))
    =\nu^{-i/(2pp')}\bigl(\tcF_p^{(1-\theta)}(x)\bigr)^*.
\end{equation}
Consequently,
\[
 H_{p,\theta}(\bG)=H_{p, 1-\theta}(\bG),\qquad
 H_{p,0}(\bG)=H_{p, 1}(\bG)=H_p(\bG).
\]
There is no shift of the minimizing center caused by $\nu$.
\end{proposition}
\begin{proof}
Let $\omega=y\varphi$ and $\eta=\overline\omega R$.
We have that $\lambda(\eta)=\hat\tau_{-i/2}(\lambda(\omega)^*)$ and
\[
 \eta=\nu^{-i/2}(R(y^*)\delta)\varphi.
\]
Note that $(\eta\delta^{-1/p'})(\cdot)=\eta(\delta^{-1/p'}\cdot )$ and $\sigma_{-i}^\varphi(\delta^{-1/p'}) =\nu^{i/p'}\delta^{-1/p'}$.
By the KMS identity, we see that 
\[
 \eta\delta^{-1/p'}
  =\nu^{i(1/p'-1/2)}(R(y^*)\delta^{1/p})\varphi.
\]
Let $\eta_p=\nu^{-i/2{p'}^2}\eta\delta^{-1/p'}.$
We have that 
\begin{align*}
\lambda(\eta_p)
=\nu^{-i/2{p'}^2}\hat\tau_{i/p'}\hat\sigma_{-i/p'}\lambda(\eta)
=\nu^{-i/2{p'}^2}\hat\sigma_{-i/p'}\hat\tau_{-i\kappa_p}\!(\lambda(\omega)^*).
\end{align*}
Then 
\[
 \tcF_p(J_p(x))
   =\nu^{-i/2{p'}^2}[\hd^{1/p'}\hat\tau_{-i\kappa_p}(\lambda(\omega)^*)]
   =\nu^{-i/2{p'}^2}
\bigl([\hat\tau_{i\kappa_p}(\lambda(\omega))\hd^{1/p'}]\bigr)^*.
\]
Since $\tcF_p^{(\theta)}(x)=\nu^{-i\kappa_p\theta/p'} [\hat\tau_{i\kappa_p \theta}(\lambda(\omega))\hd^{1/p'}]$, this is
\[
 \tcF_p^{(\theta)}(J_p (x))
   =\nu^{-i(1/(2{p'}^2)+\kappa_p \theta /p')}\bigl([\hat\tau_{i\kappa_p(1-\theta)}(\lambda(\omega))\hd^{1/p'}]\bigr)^*
   =\nu^{-i/(2pp')}(\tcF_p^{(1-\theta)}(x))^*.
\]
This proves \eqref{nu:reflection-phase}.
The remaining equalities follow from \eqref{nu:reflection-phase} directly.
\end{proof}

\begin{proposition}\label{profile:shape}
The function $\theta\mapsto \log H_{p, \theta}(\bG)$ is convex on $[0,1]$.
Moreover, either $H_{p, \theta}(\bG)\equiv H_p(\bG)$ on $[0,1]$, or
$H_{p,\theta}(\bG)<H_p(\bG)$ for $\theta\in (0,1)$.
\end{proposition}
\begin{proof}
Suppose that $0\leq \theta_1<\theta_2\leq 1$.
Apply the three-lines theorem to the function $\nu^{-i\kappa_p z /p'} [\hat\tau_{i\kappa_p z }(\lambda(x\varphi))\hd^{1/p'}]$ in the
strip between $z=\theta_1, \theta_2$, we have that
\[
 \|\tcF_p^{((1-t)\theta_1+t\theta_2)}([x d_\varphi^{1/p}])\|_{p'}
 \le\|\tcF_p^{(\theta_1)}([x d_\varphi^{1/p}])\|_{p'}^{1-t}\|\tcF_p^{(\theta_2)}([x d_\varphi^{1/p}])\|_{p'}^t.
\]
By taking suprema, this proves convexity. 
Proposition \ref{prop:symmetry} gives
\begin{align*}
\log H_{p,1/2}(\bG)\le\frac12\left(\log H_{p,\theta}(\bG)+\log H_{p, 1-\theta}(\bG)\right),
\end{align*}
which is the global minimum assertion. 
Finally a convex function bounded above by its equal endpoint values cannot reach that upper value at an interior point unless it is constant.
\end{proof}

\begin{remark}
The convexity implies the following upper bound:
\begin{equation}\label{profile:chord-bound}
 H_{p,\theta}\le H_p^{\,|2\theta-1|}
 H_{p,1/2}^{\,2\min\{\theta,1-\theta\}},\qquad0\le\theta\le1.
\end{equation}
This turns a single midpoint upper bound into bounds at every interior parameter.
\end{remark}

Define the following contraction: for $1\leq p\leq 2$,
\[
 \mathcal{G}_p^{(\theta)}:L^p(\hat\phi)\longrightarrow L^{p'}(\phi)
\]
by the sesquilinear adjoint identity
\begin{equation}\label{twYoung:inverse-adjoint}
 \int \mathcal{G}_p^{(\theta)}(\tilde{y})x^*\,d\phi
   =\int \tilde{y} \tcF_p^{(\theta)}(x)^*\,d\hat\phi,
 \qquad x\in L^p(\phi).
\end{equation}
Here $\mathcal{G}_p^{(0)}=\widehat{\tcF}_p$ and $\mathcal{G}_p^{(\theta)}$ is the inverse scaling Fourier transform on $\widehat{\bG}$. 
On the joint analytic core its two equivalent formulas are
\begin{equation}\label{twYoung:inverse-core}
 \mathcal{G}_p^{(\theta)}(\tilde{y})
  =\widehat{\tcF}_p\bigl(\widehat T_p(i\kappa_p\theta )\tilde{y}\bigr)
  =T_{p'}(i\kappa_p\theta )\widehat{\tcF}_p(\tilde{y}),
\end{equation}
which follows from the Fourier duality, scaling covariance and
\[
 \int \left(\widehat T_p (z) x \right) \,y^*\,d\hat\phi
   =\int x\bigl(\widehat T_{p'}(-\bar z)y\bigr)^*\,d\hat\phi.
\]

\begin{proposition}\label{prop:dualityconstant}
Suppose that $\bG$ is a locally compact quantum group.
For every $1\le p\le2$ and $0\le\theta\le1$,
\begin{equation}\label{duality:profile}
H_{p, \theta}(\mathbb G)=H_{p, \theta}(\widehat{\mathbb G}).
\end{equation}
\end{proposition}
\begin{proof}
It follows from the scaling Fourier duality.
\end{proof}

\begin{corollary}
Suppose that $\bG$ is a compact or discrete locally compact quantum group.
Then 
\[
H_{p, \theta}(\mathbb G)=1\qquad(1\le p\le2,\ 0\le\theta\le1).
\]
\end{corollary}
\begin{proof}
Suppose that $\mathbb G$ is compact and $\varphi(1)=1$.
The trivial group-like projection $1$ implies that $H_{p, \theta}(\bG)=1$ for $\theta\in [0,1]$.  
Applying \eqref{duality:profile}, we see that $H_{p, \theta}(\mathbb G)=1$ for a discrete locally compact quantum group $\bG$.
\end{proof}

For $0<\theta<1$ and $1<p<2$, the duality also transports extremizers.
If $\|x\|_p=1$ and
$\|\tcF_{p}^{(\theta)}(x)\|_{p'}=H_{p, \theta}(\mathbb G)$,
put
\[
 y=N_{p'}(\tcF_{p}^{(\theta)}(x)).
\]
Then $\|y\|_p=1$ and the sharp norming certificate gives
\[
 \hcF_{p}^{(\theta)}(y)=H_{p, \theta} N_p(x).
\]
Hence $y$ is a unit extremizer for the dual family at the same parameter.
Applying the same operation on the dual side returns $x$, since $N_{p'}(N_p(x))=x$ for unit $x$. 
Thus these norming maps $N_s(\cdot)$ are inverse bijections between the unit extremizer sets.

\begin{proposition}
Suppose that $\theta$ is fixed. 
The function
\[
 \alpha \longmapsto\log H_{\frac{2}{2-\alpha},\theta}(\mathbb G), \quad 0\leq \alpha \leq 1,
\]
is convex.
\end{proposition}
\begin{proof}
Recall that on the strip $0\leq \Re z\leq 1$, 
\[
 S_z(\omega)=\exp\!\left(\frac{i\theta\log\nu}{8}z^2\right)
           \hat\tau_{i\theta(1-z)/2}(\lambda(\omega)).
\]
On a vertical line $z=\alpha+iy$, we have that 
\begin{align*}
    \|S_{\alpha+i y}(\omega)\|_{2/\alpha} \leq H_{\frac{2}{2-\alpha}, \theta}(\bG)\|\omega\|_{\frac{2}{2-\alpha}}.
\end{align*}
Applying Stein interpolation theorem \ref{prop:steininter} between any two such vertical lines $\Re z=\alpha_1$ and $\Re z=\alpha_2$, we see that the logarithmic convexity for the function holds.
\end{proof}

In the end of this section, we shall compare the scaling Fourier transform with the twisted Fourier transform introduced by Youn \cite{You26}.
In the following, we shall assume that $\bG$ is a compact quantum group and the irreducible unitary representations are
\[
 u^\alpha=(u^\alpha_{ij})_{1\le i,j\le n_\alpha},\qquad \alpha\in\Irr(\bG).
\]
The comultiplication $\Delta$ satisfies
\begin{align*}
 \Delta(u^\alpha_{ij})=\sum_k u^\alpha_{ik}\otimes u^\alpha_{kj}.
\end{align*}
We diagonalize each positive Woronowicz matrix:
\begin{equation}\label{eq:Q}
 Q_\alpha=\operatorname{diag}(q_{\alpha,1},\ldots,q_{\alpha,n_\alpha}),
 \qquad q_{\alpha,i}>0,\qquad
 d_\alpha=\Tr(Q_\alpha)=\Tr(Q_\alpha^{-1}).
\end{equation}
When a single representation is being considered we shall suppress $\alpha$ and write $q_i$ and $d$.
The normalization is the one in which
$(Q_\alpha^{1/2}\otimes1)(u^\alpha)^c
(Q_\alpha^{-1/2}\otimes1)$ is unitary, where
$(u^\alpha)^c=( (u^\alpha_{ij})^*)_{ij}$.
Let $h$ be the Haar state on $\bG$.
The Schur orthogonality formulas are
\begin{align}
 h((u^\beta_{kl})^*u^\alpha_{ij})
 &=\delta_{\alpha\beta}\delta_{ik}\delta_{jl}
       \frac{q_{\alpha,i}^{-1}}{d_\alpha},\label{eq:schur1}\\
 h(u^\beta_{kl}(u^\alpha_{ij})^*)
 &=\delta_{\alpha\beta}\delta_{ik}\delta_{jl}
       \frac{q_{\alpha,j}}{d_\alpha}.\label{eq:schur2}
\end{align}
The polynomial algebra $\Pol(\bG)$ is the span of these matrix coefficients. 
We have that $ h(xy)=h(y\sigma^h_{-i}(x))$ and
\begin{align}
 \sigma^h_t(u^\alpha_{ij})=(q_{\alpha,i}q_{\alpha,j})^{it}u^\alpha_{ij}, \quad
 \tau_t(u^\alpha_{ij})=\left(\frac{q_{\alpha,i}}{q_{\alpha,j}}\right)^{it}u^\alpha_{ij}.
 \label{eq:modular-source}
\end{align}
The dual scaling action is
\begin{equation}\label{eq:dual-scaling}
 \widehat\tau_t(x)_\alpha
 =Q_\alpha^{-it}x_\alpha Q_\alpha^{it},\qquad t\in\mathbb R.
\end{equation}

\begin{equation}\label{eq:canonical-family-intro}
 \cF^{(\theta)}(x)=Q^{\kappa_p\theta}\lambda(x h)Q^{-\kappa_p\theta}.
\end{equation}
The twisted Fourier transform $\cF_{Y, p}^{(\theta)}$ is defined by
\begin{equation}\label{eq:youn-family-intro}
 \cF_{Y, p}^{(\theta)}(x)=Q^{-\kappa_p \theta}\lambda^*(x h)Q^{\kappa_p \theta},
\end{equation}
where $\left(\lambda^*(\omega)\right)_{\alpha}=(\iota\otimes \omega)((u^\alpha)^*)$, $\omega\in \cM_*$.

\begin{lemma}\label{lem:KMS-coefficient}
For $x\in\Pol(\bG)$ and $1\le p\le2$,
\begin{equation}\label{eq:KMS-coefficient}
\lambda(\sigma_{-i/p}^h(x^*)h)
=Q^{1/p'}\lambda^*(xh)^*Q^{1/p'}.
\end{equation}
\end{lemma}
\begin{proof}
Fix a block and write $u=u^\alpha$, $Q=\operatorname{diag}(q_1,\ldots,q_n)$.
By modular invariance and the KMS identity,
\begin{align*}
\lambda(\sigma_{-i/p}^h(x^*)h)_{kj}
 &=h(u_{kj}\sigma_{-i/p}^h(x^*))
 =(q_kq_j)^{-1/p}h(u_{kj}x^*)\\
 &=(q_kq_j)^{1/p'}h(x^*u_{kj})
 =(q_kq_j)^{1/p'}\lambda^*(xh)^*_{kj}.
\end{align*}
This proves the identity blockwise.
\end{proof}

For $1\le s<\infty$ let $\ell^s_\pm$ be the completions of finitely
supported block families for the norms
\begin{equation}\label{eq:norm-mp}
 \|a\|_{\pm,s}=\Big(\sum_\alpha d_\alpha
 \|a_\alpha Q_\alpha^{\pm1/s}\|_s^s\Big)^{1/s}.
\end{equation}
At $s=\infty$ use the block supremum norm.

\begin{proposition}\label{prop:equivalence}
Let $\bG$ be a compact quantum group, $1\le p\le2$, $0\le\theta\le1$, and $x\in\Pol(\bG)$. 
Then
\begin{equation}\label{eq:equivalence}
 Q^{-1/p'}\cF_p^{(\theta)}(\sigma_{-i/p}^h(x^*))Q^{-1/p'}
 =\bigl(\cF_{Y,p}^{(\theta)}(x)\bigr)^*.
\end{equation}
The two Fourier maps into $\ell^{p'}_-$ and $\ell^{p'}_+$,
respectively, have bounded extensions simultaneously, and their norms agree.
\end{proposition}
\begin{proof}
The coefficient identity and the definitions give
\[
 Q^{-1/p'}\cF_p^{(\theta)}(\sigma_{-i/p}^h(x^*))Q^{-1/p'}
 =Q^{\kappa_p\theta}\lambda^*(xh)^*Q^{-\kappa_p\theta},
\]
which is the adjoint of $\cF_{Y,p}^{(\theta)}(x)$.
The map $x\mapsto\sigma_{-i/p}^h(x^*)$ is an involutive isometry
on the embedded polynomial $L^p$ core: its spatial representative is
the adjoint of $[x d_h^{1/p}]$.
Likewise $A\mapsto(Q^{-1/p'}AQ^{-1/p'})^*$ is a bijective
conjugate-linear isometry from $\ell_-^{p'}$ to $\ell_+^{p'}$, since
\[
 \|(Q^{-1/p'}AQ^{-1/p'})^*Q^{1/p'}\|_{p'}
 =\|AQ^{-1/p'}\|_{p'}.
\]
The same conclusion holds for $p'=\infty$ with $Q^0=1$.
Density now gives simultaneous boundedness and equality of norms.
\end{proof}

\begin{remark}
In \cite{You26}, Youn shows that $\cF_{Y, p}^{(\theta)}$ could be unbounded outside the interval $[0,1]$.
One could also investigate the scaling Fourier transform for $\theta\in \bR$, in which the boundedness must first be established.
\end{remark}

\section{Equality of the Hausdorff-Young inequality}

In this section, we shall investigate the equality of the ordinary and scaling Hausdorff-Young inequality for locally compact quantum groups.

The space
\[
 \mathcal D_a(\mathcal H_\varphi,\phi)
 =\mathcal D(\mathcal H_\varphi,\phi)\cap\bigcap_{t\in\mathbb R}\mathcal D(d_\varphi^t)
\]
is a core for every spatial $L^r$-operator with $r\geq2$ and for every
real power of $d_\varphi$; 
the same assertions hold on the dual.
Here is the reduction to the standard bounded-vector core theorem
of spatial integration \cite{Hil81}. 
Put $D=d_\varphi$.
The unitaries $d_\varphi^{iu}$ preserve $\phi$-bounded vectors and their $R^\phi$ norms. 
For $2\le r<\infty$, the operators
$A_t=d_\varphi^{-it}Ad_\varphi^{it}$ vary continuously in $L^r(\phi)$.
For a bounded vector $\xi$, the identity
$Ad_\varphi^{it}\xi=d_\varphi^{it}A_t\xi$ and
\eqref{repair:bounded-vector} show continuity in the graph norm of $A$.
The Gaussian smoothings
\[
 \xi_n=\frac n{\sqrt\pi}\int_{\bR}e^{-n^2t^2}d_\varphi^{it}\xi\,dt
\]
therefore converge to $\xi$ in that graph norm. 
They remain $\phi$-bounded and lie in $\mathcal D(d_\varphi^t)$ for every real $t$, by spectral calculus. 
Since the bounded vectors form a core for $A$, this proves the assertion for $A$; 
the bounded case $r=\infty$ needs only Hilbert-space density. 
For a power $d_\varphi^t$ itself, first smooth by
$S_n=\exp(- (\log d_\varphi)^2/(4n^2))$. 
Both $S_n$ and $d_\varphi^t S_n$ are bounded. 
Approximate an arbitrary vector by bounded vectors before applying $S_n$, and then let $n$ tend to infinity in the graph norm
of $d_\varphi^t$. 
This proves the power-core assertion as well.
These common-core assertions specify which closed operator is meant when a bounded Fourier coefficient is multiplied by a spatial power.

\begin{lemma}\label{repair:Fouriercutoff}
Suppose that $x\in L^p(\phi)$ with $1\leq p \leq 2$ and $a\in\mathcal T_\varphi^2$.
Then $xa d_\varphi^{1/p'}\in L^1(\phi)$.
If $\omega_{x,a}\in\mathcal M_*$ is the normal functional with this spatial density, then
\begin{equation}\label{repair:cutoffformula}
 \tcF_p(xa)
   =\bigl[\lambda(\omega_{x,a})\hat d_{\hat\varphi}^{1/p'}\bigr].
\end{equation}
In particular the closed product on the right belongs to $L^{p'}(\hat\phi)$.
\end{lemma}
\begin{proof}
The analytic Tomita core has $a d_\varphi^{1/p'}\in L^{p'}(\phi)$;
H\"older's inequality implies that $xa d_\varphi^{1/p'}\in L^1(\phi)$.
Choose $b_n\in\mathcal T_\varphi^2$ such that
$x_n=b_n d_\varphi^{1/p}\to x$ in $L^p(\phi)$.
Because
\[
 x_n a=b_n\sigma^\varphi_{-i/p}(a)d_\varphi^{1/p}
\]
is again in the spatial Tomita-square core $\mathcal T_\varphi^2$, we have that
$\tcF_p(x_n a)=\bigl[\lambda(\omega_{x_n,a})\hat d_{\hat\varphi}^{1/p'}\bigr].$
The bounded module action and the Hausdorff--Young inequality imply
$\tcF_p(x_n a)\to \tcF_p(xa)$ in $L^{p'}(\hat\phi)$.
Furthermore
\[
 \|\omega_{x_n,a}-\omega_{x,a}\|
 \leq\|x_n-x\|_{p,\phi}\|a d_\varphi^{1/p'}\|_{p',\phi},
\]
so $\lambda(\omega_{x_n,a}) \to \lambda(\omega_{x,a})$ in operator norm.
For a vector
$\xi\in\mathcal D_a(\mathcal H_{\hat\varphi},\hat\phi)$,
Inequality \eqref{repair:bounded-vector} and the fact that $\xi\in\mathcal D(\hat d_{\hat\varphi}^{1/p'})$ give
\[
 \tcF_p(x a)\xi=\lim_n \tcF_p(x_na)\xi
      =\lim_n \lambda(\omega_{x_n,a})\hat d_{\hat\varphi}^{1/p'}\xi
      =\lambda(\omega_{x,a})\hat d_{\hat\varphi}^{1/p'}\xi.
\]
This common vector space is a core for both $ \tcF_p(x a)$ and $\hat d_{\hat\varphi}^{1/p'}$.
Since $\lambda(\omega_{x,a})$ is bounded, approximation in the graph norm of $\hat d_{\hat\varphi}^{1/p'}$ shows
$\lambda(\omega_{x,a})\hat d_{\hat\varphi}^{1/p'}\subseteq \tcF_p(x a)$.
Conversely, using that the same space is a core for $\tcF_p(x a) $ shows
$\tcF_p(x a) \subseteq[\lambda(\omega_{x,a})\hat d_{\hat\varphi}^{1/p'}]$.
Thus the closed operators are equal, proving Equation \eqref{repair:cutoffformula}.
\end{proof}

Put
\[
 y=\tcF_p(x),\qquad x=w_x|x|,\qquad y=w_{\hat x}|y|,
 \qquad p_x=w_x^*w_x,\qquad q_{\hat x}=w_{\hat x}^*w_{\hat x}.
\]

\begin{proposition}\label{repair:HYpaths}
Let $1<p<2$, $0\ne x\in L^p(\phi)$, and normalize $\|x\|_p=1$.
Then the following conditions are equivalent:
\begin{enumerate}[label=\textup{(\roman*)}]
\item $\|y\|_{p'}=1$.
\item $\widehat{\tcF}_p(N_{p'}(y))=N_p(x)$.
\item $\widehat{\tcF}_p(w_{\hat x}|y|^{p'-1})=w_x|x|^{p-1}$.
\item The single Plancherel identity holds:
\begin{equation}\label{eq:HY-Plancherel-certificate}
 \tcF_2(w_x|x|^{p/2})=w_{\hat x}|y|^{p'/2}.
\end{equation}
\item $\widehat{\tcF}_1(w_{\hat x}|y|^{p'})=w_x$.
\item $\|\widehat{\tcF}_1(w_{\hat x}|y|^{p'})\|_\infty=1$.
\end{enumerate}
Under these conditions the endpoint identities are
\begin{equation}\label{repair:endpoints}
 \tcF_1(w_x|x|^p)=w_{\hat x},\qquad
 \widehat{\tcF}_1(w_{\hat x}|y|^{p'})=w_x,
\end{equation}
and the polar parts satisfy
\begin{equation}\label{repair:Wrelations}
 W_\varphi(w_x\otimes q_{\hat x})=p_x\otimes w_{\hat x},
 \qquad
 (p_x\otimes w_{\hat x}^*)W_\varphi=w_x^*\otimes q_{\hat x}.
\end{equation}
\end{proposition}

\begin{proof}
Assume first that (i) holds, so both $x$ and $y$ are unit vectors. For $z$ in the closed strip $0\leq\Re z\leq1$, set
\[
 F(z)=\int
 \tcF_{2/(1+\Re z)}
 \bigl(w_x|x|^{p(1+z)/2}d_\varphi^{-i\Im z/2}\bigr)\hat d_{\hat\varphi}^{-i\Im z/2}
 |y|^{p'(1+z)/2}w_{\hat x}^*\,d\hat\phi.
\]
H\"older's inequality and Hausdorff-Young inequality give $|F(z)|\leq1$.
The relevant input exponent is $2/(1+\Re z)$.

Choose bounded nets $a_\alpha\in\mathcal T_\varphi^2$ and $b_\beta\in\mathcal T_{\hat\varphi}^2$ converging strongly* to $1$.
Insert $a_\alpha$ immediately before $d_\varphi^{-i\Im z/2}$ and $b_\beta$ immediately after $\hat d_{\hat\varphi}^{-i\Im z/2}$ in the displayed definition of the function $F$;
call the result $F_{\alpha,\beta}(z)$.
The Fourier slice formula and the spatial trace pairing express it as
\begin{align*}
 F_{\alpha,\beta}(z)=\big\langle
 W_\varphi\bigl(w_x|x|^{p/2}\otimes
 \hat d_{\hat\varphi}^{(1-z)/2}b_\beta|y|^{p'z/2}\bigr),\,
 d_\varphi^{(1-\bar z)/2}a_\alpha^*|x|^{p\bar z/2}
 \otimes w_{\hat x}|y|^{p'/2}
 \big\rangle_{L^2(\phi)\otimes L^2(\hat\phi)}.
\end{align*}
The inner product is linear in its first variable.
We see that this expression is analytic on the open strip and continuous on its closure (\cite{Hil81,Ter82} for more details). 
Furthermore H\"older's inequality gives the uniform estimate
\begin{align*}
 |F_{\alpha,\beta}(z)-F(z)|
 \leq{}&
 \bigl\||x|^{p/2}(a_\alpha-1)\bigr\|_{2,\phi}
       \sup_\beta\|b_\beta\|
 +\bigl\|(b_\beta-1)|y|^{p'/2}\bigr\|_{2,\hat\phi}.
\end{align*}
For the first error, use the spatially homogeneous factorization
\[
 |x|^{p(1+z)/2}(a_\alpha-1)d_\varphi^{-i\Im z/2}
 =\bigl(|x|^{pz/2}d_\varphi^{-i\Im z/2}\bigr)
   \bigl(d_\varphi^{i\Im z/2}|x|^{p/2}(a_\alpha-1)d_\varphi^{-i\Im z/2}\bigr).
\]
The first factor has $L^{2/\Re z}$-norm one; the imaginary powers occur in the spatial cocycle with the compensating power of $d_\varphi$.
The second factor has the same $L^2$-norm as $|x|^{p/2}(a_\alpha-1)$ by modular invariance of that norm.
For $\Re z=0$ the first factor is a contraction on its support.
The dual error is treated with the analogous hatted factorization.
Both terms tend to zero by the normal left and right module actions on spatial $L^2$.
Consequently $F$ is analytic on the open strip and continuous on the closure.
Since
\[
 F(2/p-1)=\int |y|^{p'}\,d\hat\phi=1,
\]
the maximum modulus principle yields $F\equiv1$.
At $z=0$ the scalar function is
\[
 F(0)=\left\langle\tcF_2(w_x|x|^{p/2}),
              w_{\hat x}|y|^{p'/2}\right\rangle=1.
\]
Both vectors have $L^2$ norm one. Equality in the Hilbert-space Cauchy--Schwarz inequality proves \eqref{eq:HY-Plancherel-certificate}. Fourier adjointness and uniqueness of norming elements give (ii), as in Proposition~\ref{prop:HY-norming-rewritten}. Since the two norms are one, (iii) follows as well.

We next justify \eqref{repair:endpoints}.
Put $u_s=w_x|x|^{ps}$, $1/2\leq s\leq1$.
The normalization and H\"older's inequality imply
\begin{equation}\label{repair:uniformcutoff}
 \|u_s(a_\alpha-1)\|_{1/s,\phi}
 \leq\bigl\||x|^{p/2}(a_\alpha-1)\bigr\|_{2,\phi},
 \qquad 1/2\leq s\leq1.
\end{equation}
For fixed $a=a_\alpha$, after multiplication on the left by $w_x|x|^{p/2}\in L^2(\phi)$, the analytic-section lemma gives
\[
 h_s:=w_x|x|^{ps}a d_\varphi^{1-s}
       \longrightarrow w_x|x|^p a=:h_1
       \quad\hbox{in }L^1(\phi),\quad s\uparrow1.
\]
Applying the Fourier transform, we obtain, on every $\hat\phi$-bounded vector $\xi\in\mathcal D(\hat d_{\hat\varphi}^{1/2})$,
\[
 \tcF_{1/s}(u_s a)\xi
   =\tcF_1(h_s)\hat d_{\hat\varphi}^{1-s}\xi
   \longrightarrow\tcF_1(h_1)\xi.
\]
Such vectors form a dense subspace; one may use
$\xi=\hat\Lambda(b)$ with $b\in\mathcal T_{\hat\varphi}$.
Indeed $\tcF_1(h_s)\to\tcF_1(h_1)$ in operator norm, while
$\hat d_{\hat\varphi}^{1-s}\xi\to\xi$ by spectral calculus.
By \eqref{repair:bounded-vector}, the Hausdorff--Young inequality and
\eqref{repair:uniformcutoff}, the error made by removing $a_\alpha$ is
bounded, uniformly in $s$, by
\[
 \bigl\||x|^{p/2}(a_\alpha-1)\bigr\|_{2,\phi}
 \max\{\|\xi\|,\|R^{\hat\phi}(\xi)\|\}.
\]
Also $w_x|x|^p a_\alpha\to w_x|x|^p$ in $L^1(\phi)$.
It follows by first fixing $\alpha$, then sending $s\uparrow1$, and finally sending $\alpha$ to infinity, that
\[
 \tcF_{1/s}(u_s)\xi\longrightarrow\tcF_1(w_x|x|^p)\xi.
\]
On the other hand, equality in the real-parameter pairing $F(2s-1)=1$ identifies $\tcF_{1/s}(u_s)$ with the unit norming vector for $w_{\hat x}|y|^{p's}$. 
Its value on $\xi$ tends to $w_{\hat x}\xi$ by spectral calculus.
For this last assertion, boundedness of $\xi$ ensures
$\xi\in\mathcal D(|y|^{p'/2})$ by
\eqref{repair:bounded-vector}, and dominated convergence applies to the spectral measure of $|y|$ for exponents between $0$ and $p'/2$.
Thus the first endpoint identity holds on a dense subspace and hence as a bounded-operator identity. 
Apply this same endpoint argument to the unit dual equality input $N_{p'}(y)$. 
By (ii) its Fourier image is $N_p(x)$, with polar part $w_x$, and its endpoint $L^1$ density is $w_{\hat x}|y|^{p'}$. This proves the second identity and (v), hence (vi).

Finally, let $\widetilde x=|x|^p$ and $\widetilde y=|y|^{p'}$.
These are positive spatial $L^1$-densities of normal states, with support
projections $p_x$ and $q_{\hat x}$.
Write $P=p_x\otimes q_{\hat x}$ and
\[
\widetilde W=(p_x\otimes w_{\hat x}^*)W_\varphi(w_x\otimes q_{\hat x}).
\]
Then $\widetilde W\in P(\mathcal M\overline\otimes\widehat{\mathcal M})P$ is a contraction.
The equality $F(1)=1$ says that the normal state with density $\widetilde x \otimes \widetilde y$ takes the value $1$ at $\widetilde W$.
That state is faithful on the corner with unit $P$.
Since $P-\Re \widetilde W\geq0$, it follows that $\Re \widetilde W=P$, and contractivity
then gives $\widetilde W=P$.
For $U=W_\varphi(w_x\otimes q_{\hat x})$ and $V=p_x\otimes w_{\hat x}$ we have $U^*U=V^*V=P$ and $V^*U=P$.
Thus $(U-V)^*(U-V)=0$ and $U=V$.
Multiplying the adjoint equality by $W_\varphi$ gives the second identity
in \eqref{repair:Wrelations}.

It remains to verify the converses, without any analytic continuation. In (ii), pairing with $x$ and applying Fourier adjointness gives $\|y\|_{p'}=1$. In (iii) the same pairing gives $\|y\|_{p'}^{p'}=1$. In (iv), Plancherel gives
\[
 1=\|w_x|x|^{p/2}\|_2
  =\|w_{\hat x}|y|^{p'/2}\|_2=\|y\|_{p'}^{p'/2}.
\]
Finally, for all unit inputs, the endpoint Fourier contraction gives
\[
 \|\widehat{\tcF}_1(w_{\hat x}|y|^{p'})\|_\infty\le\|w_{\hat x}|y|^{p'}\|_1
 =\|y\|_{p'}^{p'}\le1.
\]
Thus either (v) or (vi) forces (i). This proves all the equivalences.
\end{proof}

\begin{theorem} \label{repair:conditionalrigidity}
Let $\mathbb G$ be a locally compact quantum group and $0\ne x\in L^p(\phi)$ with $1< p<2$.
Set $y=\tcF_p(x)$ and let $p_x=\mathcal R(|x|)$ and $q_{\hat x}=\mathcal R(|y|)$.
Suppose $\|y\|_{p',\hat\phi}=\|x\|_{p,\phi}.$
If either
\begin{equation}\label{repair:modularhyp}
 p_x\in\mathcal C_\varphi
 \quad\hbox{or}\quad
 q_{\hat x}\in\mathcal C_{\hat\varphi},
\end{equation}
then $x$ is the spatial $L^p$-embedding of a bi-shift of a group-like
projection.
\end{theorem}

\begin{proof}
Suppose that $\|x\|_{p, \phi}=1$.
We first assume $q_{\hat x}\in\mathcal C_{\hat\varphi}$.
Thus $q_{\hat x}$ strongly commutes with $\hat d_{\hat\varphi}$.
Choose a bounded strongly$^*$ approximating net $a_\alpha\in\mathcal T_\varphi^2$.
For each $\alpha$, the spatial $L^1$-density
$|x|a_\alpha d_\varphi^{1/p'}$ defines a normal functional $\omega_{|x|, a_\alpha}$.
Slice the first relation in \eqref{repair:Wrelations} by the functional $\omega_{|x|, a_\alpha}$.
Since $w_x|x|=x$ and $p_x|x|=|x|$, the result is
\[
 \tcF_1(xa_\alpha d_\varphi^{1/p'})q_{\hat{x}}=\mu_\alpha w_{\hat x},
 \qquad
 \mu_\alpha=\int|x|a_\alpha d_\varphi^{1/p'}\,d\phi.
\]
Multiplication by $\hat d_{\hat\varphi}^{1/p'}$ and the commutation of $q_{\hat{x}}$ with that operator give
\begin{equation}\label{repair:keycutoff}
 \tcF_p(xa_\alpha)q_{\hat{x}} =\mu_\alpha\,[w_{\hat x}\hat d_{\hat\varphi}^{1/p'}].
\end{equation}
Here no assertion of finite Haar weight for $q_{\hat{x}}$ has yet been made.
The left side belongs to $L^{p'}(\hat\phi)$.
Moreover $xa_\alpha\to x$ in $L^p(\phi)$, so the left side converges to $yq_{\hat{x}}=y\ne0$ in $L^{p'}(\hat\phi)$.
Consequently $\mu_{\alpha_0}\ne0$ for some $\alpha_0$.
Using that index in \eqref{repair:keycutoff} proves $[w_{\hat{x}}\hat d_{\hat\varphi}^{1/p'}]\in L^{p'}(\hat\phi)$.
Since $w_{\hat{x}}^*w_{\hat{x}}=q_{\hat x}$ and $q_{\hat{x}}$ commutes with $\hat d_{\hat\varphi}$, we have $\displaystyle \bigl|[w_{\hat{x}}\hat d_{\hat\varphi}^{1/p'}]\bigr|^{p'}=q_{\hat{x}} \hat d_{\hat\varphi}.$
It follows immediately that
\[
 0< \hat\varphi(q_{\hat{x}})<\infty.
\]
The nonzero fixed vector $[w_{\hat{x}}\hat d_{\hat\varphi}^{1/p'}]$ in the Banach space $L^{p'}(\hat\phi)$ now shows that $\mu_\alpha$ converges to a scalar $\mu$ and that $y=\mu[w_{\hat{x}}\hat d_{\hat\varphi}^{1/p'}]$.
Its prescribed polar part $w_{\hat{x}}$ forces $\mu>0$; 
the norm normalization then gives $\mu=\hat\varphi(q_{\hat{x}})^{-1/p'}$.

The spatial density $w_{\hat{x}}\hat d_{\hat\varphi}$ belongs to $L^1(\hat\phi)$, because $w_{\hat{x}}q_{\hat{x}}=w_{\hat{x}}$, $q_{\hat{x}}$ commutes with $\hat d_{\hat\varphi}$ and $\hat\varphi(q_{\hat{x}})<\infty$.
Let $\eta=w_{\hat{x}}\hat\varphi$ denote the corresponding normal functional.
It satisfies $\eta q_{\hat{x}}=\eta$ and $\eta(w_{\hat{x}}^*)=\hat\varphi(q_{\hat{x}}).$
The spatial trace pairing allows $q_{\hat x}$ to be moved past the density $\hd$; 
its commutation with $\hat d_{\hat\varphi}$ and $w_{\hat{x}}q_{\hat{x}}=w_{\hat{x}}$ then give $y=\hat\varphi(q_{\hat{x}})^{-1/p'}w_{\hat x}\hat d_{\hat\varphi}^{1/p'}$.

Take the adjoint of the first relation in \eqref{repair:Wrelations}:
\[
 (w_x^*\otimes q_{\hat{x}})W_\varphi^*=p_x\otimes w_{\hat{x}}^*.
\]
Slice its second leg by $\eta$.
By the preceding identities and the definition $\hat\lambda(\eta)=(\iota\otimes\eta)(W_\varphi^*)$, one obtains
\begin{equation}\label{repair:correctslice}
 w_x^*\hat\lambda(w_{\hat{x}}\hat\varphi)=\hat\varphi(q_{\hat{x}}) p_x.
\end{equation}
Equation~\eqref{eq:HY-Plancherel-certificate} gives
\[
 \tcF_2(w_x|x|^{p/2})=\hat\varphi(q_{\hat{x}})^{-1/2}w_{\hat{x}}\hat d_{\hat\varphi}^{1/2}.
\]
The vector on the right is the spatial $L^2$-embedding of an element in the common Fourier core $L^1(\widehat{\mathbb G})\cap L^2(\widehat{\mathbb G})$.
Plancherel inversion identity and the inverse Fourier transform therefore yield
\[
w_x|x|^{p/2}=\hat\varphi(q_{\hat{x}})^{-1/2}\hat\lambda(w_{\hat{x}}\hat\varphi)d_\varphi^{1/2}.
\]
Multiply on the left by $w_x^*$ and use \eqref{repair:correctslice}.
This gives 
\begin{align}\label{repair:sourceflat}
|x|^{p/2}=\hat\varphi(q_{\hat{x}})^{1/2}p_x d_\varphi^{1/2}.
\end{align}
We show that the positivity in this equality implies the required modular invariance of $p_x$; no Gaussian uniqueness argument is needed.
The operator $[p_x d_\varphi^{1/2}]=\hat\varphi(q_{\hat{x}})^{-1/2}|x|^{p/2}$ is positive self-adjoint.
Hence
\[
 p_x d_\varphi^{1/2}\subseteq[p_x d_\varphi^{1/2}]=[p_x d_\varphi^{1/2}]^*=(p_x d_\varphi^{1/2})^*= d_\varphi^{1/2} p_x.
\]
Thus $p_x$ preserves $\mathcal D(d_\varphi^{1/2})$ and commutes with $d_\varphi^{1/2}$ there.
A bounded projection with this property reduces the self-adjoint operator $d_\varphi^{1/2}$, and therefore strongly commutes with $d_\varphi$.
Equivalently $p_x\in\mathcal C_\varphi$.
Taking the squared $L^2(\phi)$ norm in \eqref{repair:sourceflat} gives
\begin{align}\label{repair:DSproduct}
 1=\hat\varphi(q_{\hat{x}}) \varphi(p_x),
\end{align}
which proves both finiteness of $\varphi(p_x)$ and \eqref{repair:DSproduct}.
Functional calculus also gives
\[
 x=\hat\varphi(q_{\hat{x}})^{1/p}w_x d_\varphi^{1/p}.
\]
Since $p_x\in\mathcal C_\varphi$ and $\varphi(p_x)<\infty$,
$w_x\varphi$ belongs to the common $L^1$--$L^2$ Fourier domain.
The preceding Plancherel identity now yields 
\begin{align}
\lambda(w_x\varphi)&=\hat\varphi(q_{\hat{x}})^{-1}w_{\hat x}.\label{repair:polarFourier}
\end{align}
Thus $w_x$ has support $p_x$, its Fourier transform has support $q_{\hat{x}}$, and $\varphi(p_x)\hat\varphi(q_{\hat{x}})=1$.
The characterization of Donoho--Stark minimizers for locally compact quantum groups \cite{JLW18} identifies $w_x$ as a bi-shift of a group-like projection.  
The displayed formula for $x$ proves the assertion for its spatial $L^p$-embedding.

Suppose instead that $p_x\in\mathcal C_\varphi$ is the given hypothesis.
Put $z=N_{p'}(y)$. The dual norming certificate in Proposition~\ref{repair:HYpaths} gives $\widehat{\tcF}_p(z)=N_p(x)$. Both have norm one. The input has right support $q_{\hat x}$ and its Fourier image has right support $p_x$. Applying the first case to this dual equality input proves $q_{\hat x}\in\mathcal C_{\hat\varphi}$.
We may therefore apply the first case to the original $x$.
\end{proof}


\begin{corollary}\label{repair:oneHaartrace}
Suppose that $\bG$ is a locally compact quantum group and either $\varphi$ or $\hat\varphi$ is tracial.
Then the nonzero extremizers of the Hausdorff--Young inequality with optimal constant one, for $1<p<2$, are exactly the bi-shifts of group-like projections.
\end{corollary}
\begin{proof}
A tracial Haar weight has the whole von Neumann algebra as its
centralizer, so the condition \eqref{repair:modularhyp} always holds.
The forward implication follows from Theorem~\ref{repair:conditionalrigidity}; the converse is the known Fourier formula for bi-shifts of group-like projections \cite{JLW18}.
\end{proof}

\begin{proposition} \label{prop:scaling-family-equality}
Suppose $1<p<2$, $0\neq x\in L^p(\phi)$, and the scaling constant is $\nu=1$. 
The following conditions are equivalent:
\begin{enumerate}
 \item $\|\tcF_p ^{(\theta_0)}(x)\|_{p'}=\|x\|_p$ for some
       $0<\theta_0<1$;
 \item $\|\tcF_p (x)\|_{p'}=\|x\|_p$ and
       $\widehat T_{p'}(s)\tcF_p(x )=\tcF_p (x)$ for every real $s$;
 \item $\tcF_p^{(\theta)}(x)=\tcF_p (x)$ for every $\theta\in[0,1]$,
       and all these vectors have norm $\|x\|_p$.
\end{enumerate}
\end{proposition}
\begin{proof}
The strip function
\[
 F(z)=\widehat T_{p'}(i\kappa_p z)\tcF_p (x),
       \qquad0<\Re z<1,
\]
is an $L^{p'}$-valued bounded analytic function with $\|F(z)\|_{p'}\le \|x\|_p $.
This follows from the analytic-generator strip theorem, since the two real endpoints $\tcF_p (x)$ and $\widehat T_{p'}(i\kappa_p)\tcF_p (x)$ have norm at most $\|x\|_p$ and real scaling acts isometrically. 

Suppose that (1) holds.
We choose a norm-one linear functional $\ell$ on $L^{p'}(\hat\phi)$ for $F(\theta_0)$, with $\ell(F(\theta_0))=\|x\|_p$.
The scalar maximum principle gives $\ell(F(z))=\|x\|_p $ throughout the strip. 
Strict convexity of the $L^{p'}$-space implies that the intersection of its closed radius-$\|x\|_p $ ball with this supporting hyperplane consists of the single point $F(\theta_0)$. 
Thus $F$ is constant. 
In particular
\[
 \widehat T_{p'}(-\kappa_p t)F(\theta_0)
       =F(\theta_0+it)=F(\theta_0),\qquad t\in\mathbb R.
\]
Boundary continuity now gives $F(\theta_0)=\tcF_p (x)$, proving (2).
An invariant vector is fixed by every imaginary translate in this strip, which proves (2)$\Rightarrow$(3); 
(3)$\Rightarrow$(1) is immediate.
This completes the proof of the proposition.
\end{proof}

\begin{corollary} \label{cor:unimodular-scaling-equality}
Suppose that $\bG$ is a unimodular locally compact quantum group, $1<p<2$, and $0\ne x\in L^p(\phi)$.  The following conditions are equivalent.
\begin{enumerate}
 \item $\left\|\tcF_p ^{(\theta_0)}(x)\right\|_{p'}=\|x\|_p$ for some
       $0<\theta_0<1$;
\item $x$ is an $L^p$ embedding of a scaling-invariant bi-shift of a group-like projection.
\end{enumerate}
\end{corollary}
\begin{proof}
By a direct computation, we have that (2)$\Rightarrow$(1).
We shall prove (1)$\Rightarrow$(2) in the following.
Unimodularity implies that the scaling constant is $1$ and that $\hat\sigma_s=\hat\tau_s$ on $\hcM$. 
Thus the Fourier image in Proposition~\ref{prop:scaling-family-equality} is invariant under the canonical spatial modular action. 
Its polar decomposition is invariant as well, and its right support belongs to
$\mathcal C_{\hat\varphi}$. 
Theorem~\ref{repair:conditionalrigidity} applies to see the corollary holds.
\end{proof}

\section{Averaging scaling Hausdorff-Young constants}

In this section, we shall investigate the averaging scaling Hausdorff-Young constants for a locally compact quantum group with periodic scaling group and compare them with the ordinary Hausdorff-Young constants.

\begin{lemma}\label{near:analytic-lemma}
Let $\mathcal X$ be a complex Banach space satisfying the Clarkson inequality \eqref{near:Clarkson} with $s\geq 2$, and let $F$ be holomorphic in $S=\{z:0<\operatorname{Re}z<1\}$ with $\|F(z)\|\le1$.
Fix $\theta_0\in(0,1)$ and write $\varepsilon=1-\|F(\theta_0)\|.$

Suppose that $K\subset S$ is a compact set.
Let $c_K$ be a Harnack constant
from $\theta_0$ to $K$: every nonnegative harmonic function $u$ on $S$ satisfies $\displaystyle \sup_{z\in K}u(z)\le c_Ku(\theta_0).$
Then
\begin{equation}\label{near:compact-strip-estimate}
 \sup_{z\in K}\|F(z)-F(\theta_0)\|
 \le 2\left(\frac s2(c_K+1)\varepsilon\right)^{1/s}.
\end{equation}
\end{lemma}
\begin{proof}
Assume that $F(\theta_0)\neq 0$.
By the Hahn-Banach theorem, we take a norm-one complex linear functional $\phi$ on $\mathcal X$ with $\phi(F(\theta_0))=\|F(\theta_0)\|$. 
If this vector is zero, any norm-one functional suffices and the asserted bound is already at least the elementary bound $2$. 

Note that the function $u(z)=1-\operatorname{Re}\phi(F(z))$ is nonnegative and harmonic, and $u(\theta_0)=\varepsilon$.
Thus $\operatorname{Re}\phi(F(z))\ge1-c_K\varepsilon$ for
$z\in K$.
It follows that
\[
 \left\|\frac{F(z)+F(\theta_0)}2\right\|
     \ge \left(1-\frac{c_K+1}{2}\varepsilon\right)_+.
\]
Apply \eqref{near:Clarkson} and use the inequality $1-(1-a)_+^s\le s a$ for $a\ge0$. 
This proves \eqref{near:compact-strip-estimate}.
\end{proof}

\begin{theorem} \label{near:almost-invariance}
Fix $1<p<2$, $\theta_0\in(0,1)$ and $r>0$. 
There is a finite constant $K_{p,\theta_0,r}$ such that every $x\in L^p(\phi)$ with $\|x\|_p=1$ satisfies
\begin{equation}\label{near:orbit-bound}
 \sup_{|s|\le r}
 \left\|\widehat T_{p'}(s)\tcF_p^{(\theta_0)}(x)-\tcF_p^{(\theta_0)}(x)\right\|_{p'}
 \le K_{p,\theta_0,r } H_{p} \left(1-\frac{\|\tcF_p^{(\theta_0)}(x)\|_{p'}}{H_p}\right)^{1/p'}.
\end{equation}
\end{theorem}
\begin{proof}
Suppose that $x\in L^p(\phi)$ with $\|x\|_p=1$. 
We define an $L^{p'}(\hat\phi)$-valued function $F_x(z)$ as follows:
\begin{equation}\label{near:strip-function}
 F_x(z)=H_p^{-1}\widehat T_{p'}(i\kappa_p z)\tcF_p (x),
 \qquad0<\operatorname{Re}z<1,\qquad \|F_x(z)\|_{p'}\le1.
\end{equation}
It is analytic in the interior and continuous on the boundary for $x$ in the analytic core.
Then applying the argument in Proposition \ref{repair:HYpaths}, we see that it is analytic in the interior and continuous on the boundary.
Note that $F_x(\theta_0-is/\kappa_p)=H_p^{-1}\widehat T_{p'}(s)\tcF_p^{(\theta_0)}(x).$
Apply Lemma~\ref{near:analytic-lemma} to $F_x$  and to the compact vertical segment $K_r=\{\theta_0-is/\kappa_p:|s|\le r\}$,
we see that \eqref{near:orbit-bound} holds.
\end{proof}

\begin{remark}
For the displayed Harnack constant, we map $S=\{z: 0< \Re z <1\}$ to the upper half-plane by the map $z\mapsto e^{\pi iz}$ and then map $e^{\pi i\theta_0}$ to zero in the disk. 
The largest pseudo-hyperbolic distance on $K_r$ is exactly 
\begin{align*}
\rho_r=\frac{\sinh  \frac{\pi r}{2\kappa_p}}{\sqrt{\sinh^2 \frac{\pi r}{2\kappa_p}+\sin^2(\pi\theta_0)}}.
\end{align*}
The disk Harnack constant gives $\displaystyle c_{K_r}=\frac{1+\rho_r}{1-\rho_r}$. 
Finally apply the same lemma to a compact real interval inside the strip, we have that 
\begin{align}
 K_{p,\theta_0,r}=2\left(\frac {p'}2(c_{K_r}+1)\right)^{1/p'}.
 \label{near:explicit-constant}
\end{align}
\end{remark}

Assume now that the scaling actions have a common period $P>0$. Define their contractive averaging projections by
Bochner integration:
\[
 E_p (x)=\frac1P\int_0^P T_p(t)(x)\,dt,\qquad
 \widehat E_{p'}(y)=\frac1P\int_0^P\widehat T_{p'}(t)(y)\,dt.
\]
We have that 
\begin{equation}\label{near:averaging-Fourier}
 \widehat E_{p'}\tcF_p=\tcF_pE_p,\qquad
 \widehat E_{p'}\tcF_p^{(\theta)}=\tcF_pE_p \quad(0\le\theta\le1).
\end{equation}
For the second identity, the averaging projection kills every nonzero scaling frequency, while the imaginary translate is the identity on frequency zero. 
Put
\[
 H_p^{\mathrm{fix}}(\mathbb G)=\sup\left\{\|\tcF_p (x)\|_{p'}: E_p (x)=x,\ \|x\|_p\le1\right\}.
\]

\begin{theorem}\label{near:periodic-criterion}
Suppose that $\bG$ is a locally compact quantum group with periodic scaling group and $1<p<2$.
For every $0<\theta_0<1$, we have that 
\begin{equation}\label{near:gap-equivalence}
 H_{p, \theta_0}(\mathbb G)=H_p(\mathbb G)
\text{ if and only if }
 H_p^{\mathrm{fix}}(\mathbb G)=H_p(\mathbb G).
\end{equation}
More precisely, recall $K_{p,\theta_0,P}$ from \eqref{near:explicit-constant}. 
Then
\begin{equation}\label{near:quantitative-gap}
 \frac{H_{p,\theta_0}(\bG)}{H_p(\bG)}
  \le1-\left(\frac{1-H_p^{\mathrm{fix}}(\mathbb G)/H_p(\bG)}{K_{p,\theta_0, P}+1}\right)^{p'}.
\end{equation}
\end{theorem}
\begin{proof}
Let $\|x\|_p=1$, put $y=\tcF_p^{(\theta_0)}(x)$, $K=K_{p, \theta_0, P}$ and $d=1-\|y\|_{p'}/H_p\in[0,1]$. 
Averaging \eqref{near:orbit-bound} over $0\le s\le P$ gives
\[
 \|y-\widehat E_{p'}(y)\|_{p'}\le K H_p d^{1/p'}.
\]
By Equation \eqref{near:averaging-Fourier}, $\widehat E_{p'}(y)=\tcF_p(E_p (x))$, so its norm is at most $H_p^{\mathrm{fix}}(\mathbb G)$. 
Therefore
\[
 1-\frac{H_p^{\mathrm{fix}}(\mathbb G)}{H_p}
      \le d+K d^{1/p'}\le(K+1)d^{1/p'}.
\]
Rearranging and taking the supremum over $x$ proves \eqref{near:quantitative-gap}.

If $H_{p, \theta_0}(\mathbb G)=H_p$, then the same inequality forces $H_p^{\mathrm{fix}}(\mathbb G)=H_p$. 
Conversely, if $H_p^{\mathrm{fix}}(\mathbb G)=H_p$, by taking a fixed-subspace extremizing sequence, we see that $H_{p, \theta}(\mathbb G)=H_p$ for every $\theta$.
\end{proof}

\begin{proposition} \label{azgap:fixed-interpolation}
Suppose that $\bG$ is a locally compact quantum group with scaling constant one and periodic scaling group.
Then  $\log H_p^{\mathrm{fix}}(\bG)$ is convex as a function of $1/p$ on $1\le p\le2$. 

In particular, if $H_{p_0}^{\mathrm{fix}}(\bG)<1$ at one exponent $1<p_0<2$, then $H_p^{\mathrm{fix}}(\bG)<1$ at every $1<p<2$.
\end{proposition}
\begin{proof}
The period average defines compatible contractive projections \(E_1\) and \(E_2\) on the Fourier interpolation couple \(\bigl(L^1(\bG),L^2(\bG)\bigr)\). 
Their ranges form an interpolation retract: indeed, the inclusion of each range and the projection onto it are compatible contractions, so complex interpolation in both directions identifies the interpolated range with \(E_pL^p(\bG)\), equipped with its inherited norm.
The Fourier map restricted to these ranges is compatible at the endpoints. 
The interpolation theorem therefore yields log-convexity of its norm as a function of \(1/p\). 
Finally, the endpoint norms are at most one. 
If the norm is strictly below one at an interior exponent, then interpolation between that exponent and either endpoint gives strict inequality at every other interior exponent.
\end{proof}

\section{Extended convolution and Young's inequality}

Throughout this section, we shall take $\phi=J_\varphi\varphi J_\varphi$ and $\hat\phi=\hat J\hat\varphi\hat J$.
The pairing between conjugate spatial spaces is
\[
 (x,y)\longmapsto\int xy^*\,d\phi,
 \qquad x\in L^p(\phi),\quad y\in L^{p'}(\phi).
\]
It is linear in $x$ and conjugate-linear in $y$.
We use the isometric duality $(L^{p'}(\phi))^*=L^p(\phi)$ for $1\leq p'<\infty$,
including $(L^1(\phi))^*=\mathcal M$.

For $1\leq p\leq2$, let
\[
 L^p_\#(\phi)
 =\{x\in L^p(\phi):\text{there is }\tilde{x}\in L^p(\phi)
       \text{ with }\tcF_p(\tilde{x})=\tcF_p(x)^*\}.
\]
The injectivity of $\tcF_p$ makes $\tilde{x}$ unique; we write $\tilde{x}=x^\#$.
The map $\#$ is conjugate-linear.  
Its graph is closed by the Hausdorff--Young inequality, and it is an involution on its domain.
The usual simultaneous analytic core is contained in its domain and is norm dense in $L^p(\phi)$; see \cite{JLW18} for this standard densely defined dual involution.  
At $p=1$, under the identification $L^1(\phi)=\mathcal M_*$, its definition is equivalently
\begin{equation}\label{eq:repair-predual-sharp}
 \omega^\#=\overline\omega S,
 \qquad
 \overline\omega(a)=\overline{\omega(a^*)},
\end{equation}
whenever the functional on $\mathcal D(S)$ on the right has a normal bounded extension.  
The notation $\omega^\#$ denotes that extension.

For $2\leq p\leq\infty$, we use the conjugate adjoint definition.
Thus $x\in L^p_\#(\phi)$ if there is $\tilde{x} \in L^p(\phi)$ such that
\begin{equation}\label{eq:repair-adjoint-sharp}
 \int \tilde{x} y^*\,d\phi=\int y^\# x^*\,d\phi
 \quad\text{for all }y\in L^{p'}_\#(\phi),
\end{equation}
then we define $x^\#:=\tilde{x}.$
The density of $L^{p'}_\#(\phi)$ in $L^{p'}(\phi)$ implies uniqueness. 
This definition also applies at $p=2$ and agrees there with
$\tcF_2^{-1}\circ *\circ\tcF_2$.

\begin{proposition}\label{prop:repair-sharp}
For $2\leq p<\infty$, the operator $\#$ defined in \eqref{eq:repair-adjoint-sharp} is a densely defined closed conjugate-linear involution.  
At $p=\infty$, its graph is $\sigma$-weakly closed and its domain is $\sigma$-weakly dense.
For $1\leq s\leq2$ and $a\in L^s(\hat\phi)$,
\begin{equation}\label{eq:repair-sharp-Fourier-range}
 \widehat{\tcF}_s(a)\in L^{s'}_\#(\phi),\qquad
 \widehat{\tcF}_s(a)^\#=\widehat{\tcF}_s(a^*).
\end{equation}
Moreover,
\begin{equation}\label{eq:repair-sharp-infty}
 L^\infty_\#(\phi)
 =\{x\in\mathcal M:x^*\in\mathcal D(S)\},
 \qquad x^\#=S(x^*).
\end{equation}
\end{proposition}

\begin{proof}
Closedness follows by passing to the limit in each pairing in \eqref{eq:repair-adjoint-sharp}. 
The same argument gives $\sigma$-weak graph closedness at $p=\infty$, since the test elements $y,y^\#$ belong to $L^1(\phi)$.
To verify involutivity, suppose $x^\#=\tilde{x}$.  
Replacing $y$ by $y^\#$ in \eqref{eq:repair-adjoint-sharp} gives
\[
 \int \tilde{x}(y^\#)^*\,d\phi=\int yx^*\,d\phi.
\]
Conjugation and the cyclicity of the spatial $L^1$ pairing give
\[
 \int xy^*\,d\phi=\int y^\#\tilde{x}^*\,d\phi.
\]
Thus $\tilde{x}\in L^p_\#(\phi)$ and $\tilde{x}^\#=x$.

For Equalities \eqref{eq:repair-sharp-Fourier-range}, take $y\in L^s_\#(\phi)$.
The Fourier pairing identity and $\tcF_s(y^\#)=\tcF_s(y)^*$ give
\begin{align*}
 \int\widehat{\tcF}_s(a^*)y^*\,d\phi
 &=\int a^*\tcF_s(y)^*\,d\hat\phi\\
 &=\overline{\int a\tcF_s(y)\,d\hat\phi}\\
 &=\overline{\int\widehat{\tcF}_s(a)(y^\#)^*\,d\phi}
 =\int y^\#\widehat{\tcF}_s(a)^*\,d\phi.
\end{align*}
This proves \eqref{eq:repair-sharp-Fourier-range}.
The Fourier range is norm dense in $L^{s'}(\phi)$ for $1<s\leq2$
and $\sigma$-weakly dense in $\mathcal M$ for $s=1$.
These are the density consequences of the Fourier pairing identity
and injectivity of the Fourier transform, by Hahn--Banach separation.
This also proves the stated domain-density assertions.

Finally, if $x^*\in\mathcal D(S)$, then for
$y=d\omega/d\phi\in L^1_\#(\phi)$,
\[
 \int S(x^*)y^*\,d\phi
 =\overline\omega(S(x^*))
 =\omega^\#(x^*)
 =\int y^\#x^*\,d\phi.
\]
Conversely, \eqref{eq:repair-adjoint-sharp} says that $\overline\omega(z)=\overline\omega S(x^*)$ for every normal
functional for which $\overline\omega S$ is bounded and normal.
The graph of $S$ is $\sigma$-weakly closed.  
Its preannihilator consists precisely of the pairs of normal functionals $(-\eta S,\eta)$ for which $\eta S$ has a normal bounded extension.
The bipolar theorem therefore puts $(x^*,z)$ in the graph of $S$.
This proves \eqref{eq:repair-sharp-infty}.
\end{proof}


The region labels refer to the following closed portions of the Young triangle; formulas on an overlap are compared explicitly rather than identified directly:
\begin{center}
\begin{tabular}{cll}
\toprule
Region & Input exponents & Output exponent \\
\midrule
I & $1\le p,s\le2$ & $1\le r\le2$ \\
II & $1\le p,s\le2$ & $2\le r\le\infty$ \\
III & $1\le s\le2\le p\le\infty$ & $r\ge p$ \\
IV & $1\le p\le2\le s\le\infty$ & $r\ge s$ \\
\bottomrule
\end{tabular}
\end{center}
In every row the exponents satisfy
\begin{equation}\label{eq:repair-young-exponents}
 1+\frac1r=\frac1p+\frac1s.
\end{equation}
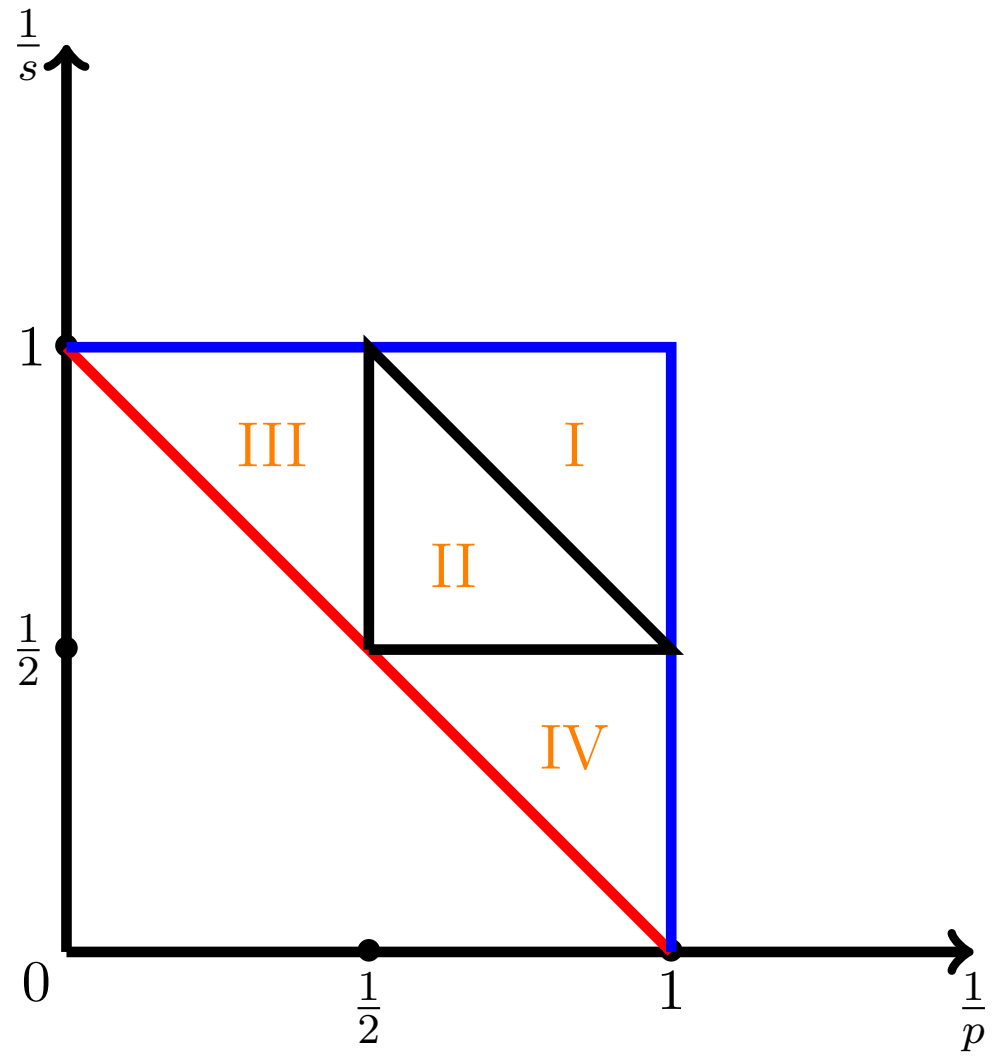
\begin{figure}
\begin{center}
\begin{tikzpicture}[scale=2]
\draw [line width=2] [->](0,0)--(0,3);
\node at (-0.1,-0.1) {0};
\node [left] at (0, 3) {$\frac{1}{s}$};
\draw [line width=2]  [->](0,0)--(3,0);
\node [below] at (1,0) {$\frac{1}{2}$};
\node at (1,0) {$\bullet$};
\node [below] at (2, 0) {1};
\node at (2,0) {$\bullet$};
\node [below] at (3,0) {$\frac{1}{p}$};
\node [left] at (0,1) {$\frac{1}{2}$};
\node at (0,1) {$\bullet$};
\node [left] at (0, 2) {1};
\node at (0,2) {$\bullet$};
\draw [red, line width=2] (0,2)--(2,0);
\draw [blue, line width=2] (0,2)--(2,2)--(2,0);
\draw [line width=2] (1,1)--(2,1)--(1,2)--(1,1);
\node at (1.28, 1.28) {\color{orange}\large $\mathrm{II}$};
\node at (1.68, 1.68) {\color{orange}\large $\mathrm{I}$};
\node at (0.68, 1.68) {\color{orange}\large $\mathrm{III}$};
\node at (1.68, 0.68) {\color{orange}\large $\mathrm{IV}$};
\end{tikzpicture}
\end{center}
\caption{The Regions of Young's inequality}
\label{Fig:young}
\end{figure}
Let $1\leq p,s\leq2$ and $\displaystyle 1+\frac1r=\frac1p+\frac1s$.
For $r\leq2$, denote by
\[
 \mathsf C^{\mathrm I}_{p,s}(x,y)
 :=x*\rho_{-i/p'}(y)
\]
the Young convolution of \cite{LWW17}, acting on the spatial spaces.
We recall the two properties of that construction which are used here:
\begin{align}
 \|\mathsf C^{\mathrm I}_{p,s}(x,y)\|_{r,\phi}&\leq\|x\|_{p,\phi}\|y\|_{s,\phi},
 \label{eq:repair-region-I-bound}\\
 \tcF_r(\mathsf C^{\mathrm I}_{p,s}(x,y))&=\tcF_p(x)\tcF_s(y).
 \label{eq:repair-region-I-Fourier}
\end{align}
The product in \eqref{eq:repair-region-I-Fourier} is the closed spatial $L^{r'}$ product.  
In particular, for $p=s=r=1$ this convolution is the usual predual convolution
\begin{equation}\label{eq:repair-predual-convolution}
 \mathsf C^{\mathrm I}_{1,1}
 \left(\frac{d\omega}{d\phi},\frac{d\eta}{d\phi}\right)
 =\frac{d((\omega\otimes\eta)\Delta)}{d\phi}.
\end{equation}
These are established Young inequality.

For $r\geq2$, define
\begin{equation}\label{eq:repair-region-II}
 \mathsf C^{\mathrm {II}}_{p,s}(x,y)
 :=\widehat{\tcF}_{r'}\bigl(\tcF_p(x)\tcF_s(y)\bigr).
\end{equation}
Indeed, $1/r'=1/p'+1/s'$, so H\"older's inequality puts the argument of $\widehat{\tcF}_{r'}$ in $L^{r'}(\hat\phi)$.
Since $1\leq r'\leq2$, the Hausdorff-Young inequality gives
\begin{equation}\label{eq:repair-region-II-bound}
 \|\mathsf C^{\mathrm {II}}_{p,s}(x,y)\|_{r,\phi}
 \leq\|\tcF_p(x)\|_{p',\hat\phi}
       \|\tcF_s(y)\|_{s',\hat\phi}
 \leq\|x\|_{p,\phi}\|y\|_{s,\phi}.
\end{equation}
At $r=2$, \eqref{eq:repair-region-I-Fourier} and
$\widehat{\tcF}_2=\tcF_2^{-1}$ show that the two definitions agree.
We consequently write $\mathsf C_{p,s}$ for this unambiguously
defined lower-square product.

The next definition gives bounded bilinear products directly.
It does not require either argument to belong to the domain of an unbounded operation such as $\#$ or $\rho_{-i/p'}$.

\begin{definition}\label{def:repair-transposed-products}
Suppose that $\displaystyle 1+\frac1r=\frac1p+\frac1s.$
If $1\leq p\leq2\leq s\leq\infty$, define
$\mathsf C^{\mathrm {IV}}_{p,s}(x,y)\in L^r(\phi)$ by
\begin{equation}\label{eq:repair-IV-duality}
 \int \mathsf C^{\mathrm {IV}}_{p,s}(x,y)a^*\,d\phi
 =\overline{\int
     \mathsf C^{\mathrm I}_{p,r'}(x^*,a)y^*\,d\phi},
 \qquad a\in L^{r'}(\phi).
\end{equation}
If $1\leq s\leq2\leq p\leq\infty$, define $\mathsf C^{\mathrm {III}}_{p,s}(x,y)\in L^r(\phi)$ by
\begin{equation}\label{eq:repair-III-duality}
 \int \mathsf C^{\mathrm {III}}_{p,s}(x,y)a^*\,d\phi
 =\overline{\int
     \mathsf C^{\mathrm I}_{r',s}(a,y^*)x^*\,d\phi},
 \qquad a\in L^{r'}(\phi).
\end{equation}
\end{definition}

\begin{theorem}[Young's inequality]
\label{thm:repair-all-Young}
Suppose that $\bG$ is a locally compact quantum group.
The products $\mathsf C^{\mathrm {III}}_{p,s}(x,y)$ and $\mathsf C^{\mathrm {IV}}_{p,s}(x,y)$ exist, are unique, and are bilinear.
We have that 
\begin{align}
 \|\mathsf C^{\mathrm {IV}}_{p,s}(x,y)\|_{r,\phi}
 &\leq\|x\|_{p,\phi}\|y\|_{s,\phi},
 \label{eq:repair-IV-bound}\\
 \|\mathsf C^{\mathrm {III}}_{p,s}(x,y)\|_{r,\phi}
 &\leq\|x\|_{p,\phi}\|y\|_{s,\phi}.
 \label{eq:repair-III-bound}
\end{align}
Together with \eqref{eq:repair-region-I-bound} and \eqref{eq:repair-region-II-bound}, these give contractive bilinear products for every exponent triple in \eqref{eq:repair-young-exponents}.
\end{theorem}

\begin{proof}
In Region IV, the exponent relation gives
\[
 r\geq s\geq2,\qquad
 1\leq r'\leq2,\qquad
 1+\frac1{s'}=\frac1p+\frac1{r'}.
\]
Thus $\mathsf C^{\mathrm I}_{p,r'}(x^*,a)$ is defined and lies in $L^{s'}(\phi)$.  
By H\"older inequality and \eqref{eq:repair-region-I-bound}, the modulus of the right-hand side of \eqref{eq:repair-IV-duality} is at most
\[
 \|x^*\|_{p,\phi}\|a\|_{r',\phi}\|y\|_{s,\phi}
 =\|x\|_{p,\phi}\|y\|_{s,\phi}\|a\|_{r',\phi}.
\]
It is a bounded conjugate-linear functional of $a$.
Here $r'\leq2<\infty$, so spatial duality gives a unique representing element of $L^r(\phi)$, including when $r=\infty$.
The same estimate gives \eqref{eq:repair-IV-bound}.
Conjugation of $x^*$ in the right-hand side makes the result
linear in $x$, and conjugation of $y^*$ makes it linear in $y$.

In Young Region III, the corresponding identities are
\[
 r\geq p\geq2,\qquad
 1\leq r'\leq2,\qquad
 1+\frac1{p'}=\frac1{r'}+\frac1s.
\]
Apply the same argument to $\mathsf C^{\mathrm I}_{r',s}(a,y^*)\in L^{p'}(\phi)$.
This proves existence, bilinearity, and \eqref{eq:repair-III-bound}.
\end{proof}

\begin{remark}[Relation with the notation involving $\#$]
The transposed products can be denoted by
\[
 \mathsf C^{\mathrm {IV}}_{p,s}(x,y)
   =x^{*\#}*\rho_{-i/p'}(y),\qquad
 \mathsf C^{\mathrm {III}}_{p,s}(x,y)
   =x*\rho_{-i/p'}(y^{*\#}),
\]
provided these are understood as names for
\eqref{eq:repair-IV-duality} and
\eqref{eq:repair-III-duality} on general $L^p$ inputs.
One must not first apply an unbounded $\#$ or an imaginary-time
map to an arbitrary input.  The precise Frobenius reciprocity
identities for these products are their defining pairings.
For example, for $1\leq p,s,r\leq2$,
\begin{align*}
 \int\mathsf C^{\mathrm I}_{p,s}(x,y)a^*\,d\phi
 &=\overline{\int
     \mathsf C^{\mathrm {IV}}_{p,r'}(x^*,a)y^*\,d\phi}\\
 &=\overline{\int
     \mathsf C^{\mathrm {III}}_{r',s}(a,y^*)x^*\,d\phi}.
\end{align*}
The stars on $x^*$ and $y^*$ make these transpose conventions bilinear. Their compatibility with the canonical products on common edges is examined in Section~\ref{sec:edge-compatibility}.
\end{remark}

\section{Scaling convolution and equality of Young inequality}

In this section, we shall extend the convolutions to the scaling convolution and investigate the equality of scaling Young inequality.

\begin{definition}
Suppose that $1\leq p,s\leq2$, $1\leq r\leq\infty$, $1/r+1=1/p+1/s$, and $x\in L^p(\phi)$, $y\in L^s(\phi)$.
We say the pair $(x, y)$ is a bishift $(p,s)$-pair of a group-like projection if 
\begin{align*}
    x=\mu_x w_x d_\varphi^{1/p}, \quad y=\mu_y w_y d_\varphi^{1/s},\quad \mu_x, \mu_y\neq 0, 
\end{align*}
where $w_x\varphi, w_y\varphi$ are bishift of a group-like projections, $w_x, w_y$ are polar parts of $x, y$ and
\begin{align*}
    w_{\hat x}^*w_{\hat x}=w_{\hat y}w_{\hat y}^*,
\end{align*}
where $\hat{x}=\tcF_p(x)$, $\hat{y}=\tcF_s(y)$ and $w_{\hat{x}}$, $w_{\hat{y}}$ are the polar parts respectively.
We say a bishift $(p,s)$-pair $(x, y)$ of a group-like projection is scaling-invariant if $w_x, w_y$ are invariant under real scaling actions. 
\end{definition}

\begin{proposition}
Suppose that $\bG$ is a unimodular locally compact quantum group.
Suppose that $(x, y)$ is a scaling-invariant bishift $(p, s)$-pair of a group-like projection. 
Then $\mathcal{S}(x)=\mathcal{S}(y)$.
\end{proposition}
\begin{proof}
Scaling covariance makes $w_{\hat x},w_{\hat y}$ invariant under $\widehat\tau$.
Unimodularity of the locally compact quantum group gives $\widehat\sigma=\widehat\tau$, so both polar parts lie in $\mathcal C_{\hat\varphi}$. 
The restriction of $\hat\varphi$ to this centralizer is tracial.
By the fact that $w_{\hat x}^*w_{\hat x}=w_{\hat y}w_{\hat y}^*$, we have that $\widehat{\mathcal S}(x)=\widehat{\mathcal S}(y)$.
By the Donoho-Stark uncertainty principles, we have that $\mathcal{S}(x)\widehat{\mathcal S}(x)=1$ and $\mathcal{S}(y)\widehat{\mathcal S}(y)=1$.
Hence $\mathcal{S}(x)=\mathcal{S}(y)$.
\end{proof}

\subsection{Young Region II}

For $\theta_1,\theta_2,\theta_3\in[0,1]$, set
\begin{equation}\label{twYoung:definition}
\mathsf C_{p,s}^{\boldsymbol\theta} (x, y)
 :=\mathcal{G}_{r'}^{(\theta_3)}\bigl(\tcF_p^{(\theta_1)}(x)\tcF_s^{(\theta_2)}(y)\bigr).
\end{equation}

\begin{theorem}\label{twYoung:family}
Suppose that $\displaystyle 1\le p,s\le2,\quad 2\le r\le\infty,\quad 1+\frac1r=\frac1p+\frac1s.$
The map $\mathsf C_{p,s}^{\boldsymbol\theta}$ is a bounded bilinear map from $L^p(\phi)\times L^s(\phi)$ to $L^r(\phi)$, and
\begin{equation}\label{twYoung:bound}
 \|\mathsf C_{p,s}^{\boldsymbol \theta}(x,y)\|_{r,\phi}
       \le\|x\|_{p,\phi}\|y\|_{s,\phi}.
\end{equation}
\end{theorem}
\begin{proof}
The exponent relation is $1/r'=1/p'+1/s'$. 
The spatial H\"older inequality and the Hausdorff-Young inequalities give
\begin{align*}
 \|\mathsf C_{p,s}^{\boldsymbol \theta}(x,y)\|_r
 &\le\|\tcF_p^{(\theta_1)}(x)\tcF_s^{(\theta_2)}(y)\|_{r'}\\
 &\le\|\tcF_p^{(\theta_1)}(x)\|_{p'}
       \|\tcF_s^{(\theta_2)}(y)\|_{s'}
 \le\|x\|_p\|y\|_s.
\end{align*}
This proves all the assertions, since $\kappa_2=0$.
\end{proof}

\begin{remark}
At $\boldsymbol{\theta}=(0,0,0)$, this is ordinary Young convolution in Region II. 
If $r=2$, the third parameter has no effect. 
If $r=\infty$, necessarily $p=s=2$, so the first two
parameters have no effect.
\end{remark}

\begin{remark}[Effective region of the parameter $\boldsymbol{\theta}$]

Reformulate the convolution $\mathsf C_{p,s}^{\boldsymbol \theta}$, we obtain that 
\begin{equation}\label{twYoung:effective-core}
 \mathsf C_{p,s}^{\boldsymbol \theta}(x,y)
 =\widehat{\tcF}_{r'}\left(
       \widehat T_{p'}(i(\kappa_p\theta_1 +\kappa_{r'}\theta_3)) \tcF_p(x) \,
       \widehat T_{s'}(i(\kappa_s\theta_2+\kappa_{r'}\theta_3))\tcF_s(y)\right).
\end{equation}
We see that there are two effective coordinates
\[
 (\alpha,\beta)=(\kappa_p\theta_1 +\kappa_{r'}\theta_3,\kappa_s\theta_2+\kappa_{r'}\theta_3).
\]
Hence the effective pair $(\alpha, \beta)$ lie in the compact convex set
\begin{equation}\label{twYoung:polygon}
 \Omega_{p,s}=
 \bigl([0,\kappa_p]\times[0,\kappa_s]\bigr)+\{(w,w):0\le w\le\kappa_{r'}\},
\end{equation}
with the relation
\[
 \kappa_p+\kappa_s=\frac1r,\qquad  \kappa_p+\kappa_s+\kappa_{r'}=\frac12.
\]
We may therefore write the convolution as $\mathsf C_{p,s}^{[\alpha,\beta]}$.
If $r>2$ and $1< p, s<2$, $\kappa_p, \kappa_s, \kappa_{r'}$ are positive.
The vertices of the hexagon $\Omega_{p, s}$ are
\[
 (0,0),\ (\kappa_p,0),\ (\kappa_p+\kappa_{r'},\kappa_{r'}),\ (\kappa_p+\kappa_{r'},\kappa_s+\kappa_{r'}),\ (\kappa_{r'},\kappa_s+\kappa_{r'}),\ (0,\kappa_s).
\]
At $r=2$,  the parameter set is the rectangle $[0,\kappa_p]\times[0,\kappa_s]$.
\end{remark}

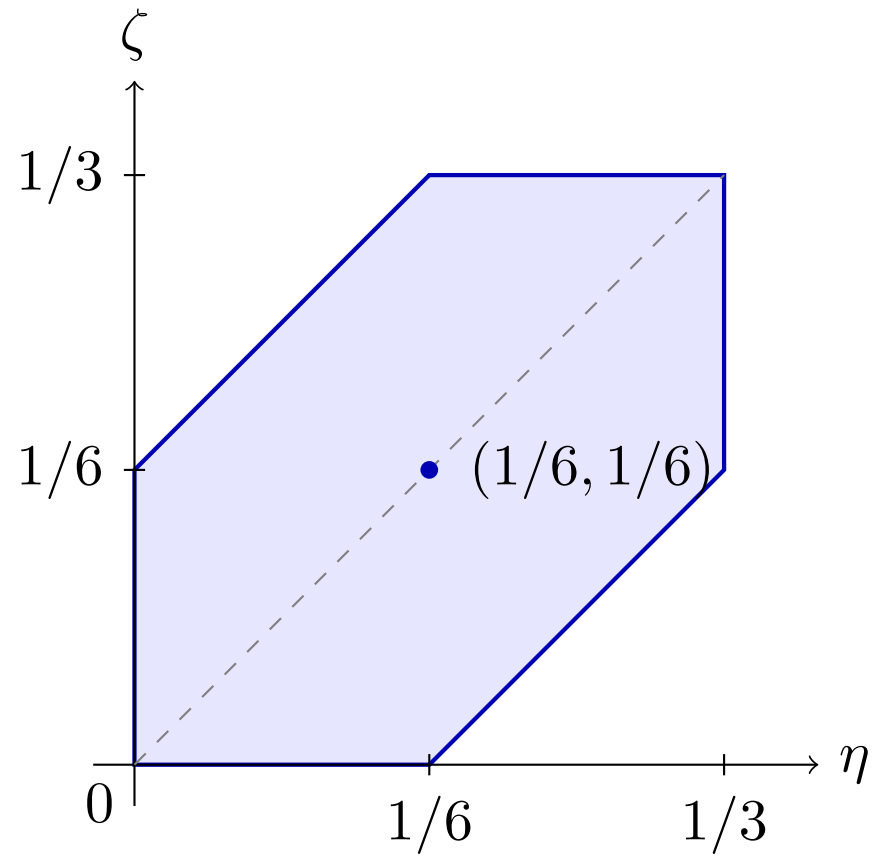
\begin{figure}[ht]
\centering
\begin{tikzpicture}[scale=3.9]
 \fill[blue!10] (0,0)--(.5,0)--(1,.5)--(1,1)--(.5,1)--(0,.5)--cycle;
 \draw[thick,blue!70!black] (0,0)--(.5,0)--(1,.5)--(1,1)--(.5,1)--(0,.5)--cycle;
 \draw[->] (-.07,0)--(1.16,0) node[right] {$\eta$};
 \draw[->] (0,-.07)--(0,1.16) node[above] {$\zeta$};
 \foreach \x/\lab in {.5/{1/6},1/{1/3}} {
  \draw (\x,.018)--(\x,-.018) node[below] {$\lab$};
  \draw (.018,\x)--(-.018,\x) node[left] {$\lab$};
 }
 \node[below left] at (0,0) {$0$};
 \draw[dashed,gray] (0,0)--(1,1);
 \fill[blue!70!black] (.5,.5) circle (.015);
 \node[right] at (.5,.5) {$\ (1/6,1/6)$};
\end{tikzpicture}
\caption{The effective Young parameters for $p=s=3/2$, $r=3$.
The marked interior point is already represented by the cube vertex $(\theta_1,\theta_2,\theta_3)=(0,0,1)$. 
Interior equality is therefore a property of the effective map, not of a particular cube representation.}
\end{figure}

\begin{theorem}\label{twYoung:certificate}
Suppose that $1<p,s\le2$, $2\le r<\infty$ and $1/r+1=1/p+1/s$,
and let $0\ne x\in L^p(\phi)$ and $0\ne y\in L^s(\phi)$. 
The equality in \eqref{twYoung:bound} holds if and only if
\begin{align}
 \|\tcF_p^{(\theta_1)}(x)\|_{p'}&=\|x\|_p,\qquad
 \|\tcF_s^{(\theta_2)}(y)\|_{s'}=\|y\|_s,\label{twYoung:forward-equality}\\
 \frac{|\tcF_p^{(\theta_1)}(x)|^{p'}}{\|\tcF_p^{(\theta_1)}(x)\|_{p'}^{p'}}
  &=\frac{|\tcF_s^{(\theta_2)}(y)^*|^{s'}}{\|\tcF_s^{(\theta_2)}(y)\|_{s'}^{s'}},
       \label{twYoung:Holder-equality}\\
 \left\|\mathcal{G}_{r'}^{(\theta_3)}(\tcF_p^{(\theta_1)}(x)\tcF_s^{(\theta_2)}(y))\right\|_r&=\left\|\tcF_p^{(\theta_1)}(x)\tcF_s^{(\theta_2)}(y)\right\|_{r'}.
       \label{twYoung:inverse-equality}
\end{align}
At $r=2$ the last condition is automatic.
\end{theorem}
\begin{proof}
Equality at the two ends of the norm chain in the proof of Theorem~\ref{twYoung:family} forces equality at each step. 
Since both inputs are nonzero, the product of the last two contraction bounds is sharp only when both individual bounds are sharp. 
The equality of the spatial H\"older inequality gives \eqref{twYoung:Holder-equality}. 

Conversely, the three conditions make every step an equality. 
Plancherel identity proves the last assertion.
\end{proof}

\begin{theorem}\label{thm:cube-any-two}
Suppose that $\bG$ is a locally compact quantum group with scaling constant $\nu=1$, $1<p,s<2$, $2\leq r<\infty$ and $1/r+1=1/p+1/s$.
Let $0\ne x\in L^p(\phi)$ and $0\ne y\in L^s(\phi)$.  
The following conditions are equivalent.
\begin{enumerate}
\item Equality in \eqref{twYoung:bound} holds for one $\boldsymbol\theta$ with $\theta_1,\theta_2\in(0,1)$.  
\item Equality in \eqref{twYoung:bound} holds for every $\boldsymbol\theta\in[0,1]^3$.
\item The equality of canonical Young inequality holds, namely
\[
 \left\|\widehat{\tcF}_{r'}(\tcF_p(x)\tcF_s(y))\right\|_r
     =\|x\|_p\|y\|_s,
\]
and $\tcF_p(x)$ and $\tcF_s(y)$ are invariant under the real dual-scaling actions on their spatial spaces.
\end{enumerate}
When these conditions hold, the output itself is independent of all three parameters:
\[
\mathsf C_{p,s}^{\boldsymbol\theta}(x,y)=\mathsf C_{p,s}^{\boldsymbol0}(x,y).
\]
\end{theorem}
\begin{proof}
Put 
\begin{align*}
    \widehat{x} =\tcF_p^{(\theta_1)}(x), \quad \widehat{y}=\tcF_s^{(\theta_2)}(y).
\end{align*}

Suppose (1) holds. Equality in the norm chain gives equality in both forward Hausdorff--Young inequalities. Proposition~\ref{prop:scaling-family-equality} implies that $\tcF_p(x)$ and $\tcF_s(y)$ are fixed under the real dual-scaling actions, and that their twisted transforms equal their untwisted transforms. Their product is therefore fixed under the dual-scaling action. By the core formula and closedness, the inverse transform of this product is independent of $\theta_3$. Thus the canonical output already has maximal norm and (2) and (3) follow. Conversely, (2) trivially implies (1). If (3) holds, the same covariance and closedness show that all three twists leave the output unchanged, proving (2).
\end{proof}

\begin{corollary}\label{cor:hex-equality}
Suppose that $\bG$ is a locally compact quantum group with scaling constant $\nu=1$, $1<p,s< 2$, $2\leq r<\infty$ and $1/r+1=1/p+1/s$. 
For $0\ne x\in L^p(\phi)$ and $0\ne y\in L^s(\phi)$, the following are equivalent:
\begin{enumerate}
\item Equality in \eqref{twYoung:bound} holds for one $(\alpha,\beta)$ in the interior of $\Omega_{p,s}$.
\item Equality in \eqref{twYoung:bound} holds at every point of $\Omega_{p,s}$.
\item The equality of canonical Young inequality for Region II holds and both inputs are invariant under their real scaling actions.
\end{enumerate}
When $r>2$, equality at either cube vertex $(0,0,1)$ or $(1,1,0)$ already implies these conclusions. 
At $r=2$, the same equivalence holds with an interior point of the rectangle $[0,\kappa_p]\times[0,\kappa_s]$.
\end{corollary}
\begin{proof}
At an interior point $(\alpha, \beta)$ of the nondegenerate hexagon, we have that 
\[
 \max\{0,\alpha-\kappa_p,\beta-\kappa_s\}
  <\min\{\kappa_{r'},\alpha,\beta\}.
\]
Choose $w$ strictly between these two quantities and set $u=\alpha-w$, $v=\beta-w$. 
Then $0<u<\kappa_p$, $0<v<\kappa_s$, $0<w<\kappa_{r'}$,
so the same map has a representation with all three cube parameters interior. 
Apply Theorem \ref{thm:cube-any-two}, we see that (2), (3) hold. 

Conversely, equality for all cube parameters gives equality for all points in the hexagon. 
This shows that (3) implies (2).

The points $(\kappa_{r'},\kappa_{r'})$ and $(\kappa_p,\kappa_s)$ are interior when $\kappa_p,\kappa_s,\kappa_{r'}>0$, proving the cube-vertex assertion. 
When $\kappa_{r'}=0$, interiority of the rectangle means interiority of both input twists, which is exactly the case used in the proof.
\end{proof}

\begin{theorem}\label{twYoung:structural}
Suppose that $\bG$ is a unimodular locally compact quantum group, $1<p,s<2$, $2\le r<\infty$, and $1+1/r=1/p+1/s$.  
For nonzero $x\in L^p(\phi)$ and
$y\in L^s(\phi)$, equality at one interior point of $\Omega_{p,s}$
is equivalent to $(x,y)$ being a scaling-invariant bi-shift
$(p,s)$-pair in the sense of the preceding definition.  
In that case equality holds at every point of $\Omega_{p,s}$.
\end{theorem}
\begin{proof} 
Suppose first that equality holds at an interior point.  
By Corollary~\ref{cor:hex-equality}, it holds at every cube point and the two inputs are scaling-invariant.  
Equality in the two forward Hausdorff--Young inequalities and
Corollary~\ref{cor:unimodular-scaling-equality} show that $x,y$ are
embeddings of scaling-invariant bi-shifts of a group-like projection. 
Write their Fourier images as
\[
 \hat{x}=\tcF_p(x)=\hat{\mu}_x w_{\hat{x}} \hd^{1/p'},\qquad
 \hat{y}=\tcF_s(y)=\hat{\mu}_y w_{\hat{y}} \hd^{1/s'},\qquad \hat{\mu}_x,\hat{\mu}_y>0,
\]
where $w_{\hat{x}},w_{\hat{y}}$ are partial isometries. 
Scaling covariance makes $\hat{x},\hat{y}$ invariant under the spatial dual-scaling actions.  
Since $\delta=1$, these are the spatial modular actions of $\hat\varphi$; 
hence $w_{\hat{x}},w_{\hat{y}}$ and their support projections belong to the centralizer of $\hat\varphi$.  
The H\"older equality criterion \eqref{twYoung:Holder-equality} therefore compares two scalar multiples of $(w_{\hat{x}}^*w_{\hat{x}})\hd$ and $(w_{\hat{y}}w_{\hat{y}}^*)\hd$. 
Since $\hd$ is nonsingular and both projections are nonzero, their support projections coincide:
$w_{\hat{x}}^*w_{\hat{x}}=w_{\hat{y}}w_{\hat{y}}^*$.  
This is precisely the matching condition in the definition of a bi-shift pair of a group-like projection.

Conversely, suppose that $(x,y)$ is such a pair.  
By the properties of the bishift of a group-like projection, we have that
\[
 \hat{x}=\mu_x \mathcal{S}(x) w_{\hat{x}} \hd^{1/p'},\quad
 \hat{y}=\mu_y\mathcal{S}(y) w_{\hat{y}} \hd^{1/s'},\quad
 \hat\varphi(w_{\hat{x}}^*w_{\hat{x}})=\mathcal{S}(x)^{-1},\quad
 \hat\varphi(w_{\hat{y}}^*w_{\hat{y}})=\mathcal{S}(y)^{-1}.
\]
Here the phases are included in $w_{\hat{x}},w_{\hat{y}}$, and $\mu_x,\mu_y>0$.
Scaling invariance and $\hat\sigma=\hat\tau$ put $w_{\hat{x}},w_{\hat{y}}$ in the dual Haar centralizer.  
Since the Haar weight is tracial on its centralizer, the matching $w_{\hat{x}}^*w_{\hat{x}}=w_{\hat{y}}w_{\hat{y}}^*$ gives
\[
 \mathcal{S}(x)^{-1}=\hat\varphi(w_{\hat{x}}^*w_{\hat{x}})=\hat\varphi(w_{\hat{y}}w_{\hat{y}}^*)
 =\hat\varphi(w_{\hat{y}}^*w_{\hat{y}})=\mathcal{S}(y)^{-1}.
\]
Put $\hat{w}=w_{\hat{x}}w_{\hat{y}}$. 
Then $\hat{w}$ is a partial isometry with $\hat{w}^*\hat{w}=w_{\hat{y}}^*w_{\hat{y}}$ and, for
$1/p'+1/s'=1/r'$,
\[
 \hat{x} \hat{y}=\mu_x\mu_y \mathcal{S}(x)^2 \hat{w} \hd^{1/r'},\qquad
 \|\hat{x} \hat{y}\|_{r'}=\mu_x\mu_y \mathcal{S}(x)^{2-1/r'}
 =\|x\|_p\|y\|_s.
\]

We next show that the inverse transform preserves the norm of this
product.  
We first embed the same bi-shift pair of a group-like projection at exponent $6/5$ and let $v=\mathsf C^{\mathrm I}_{6/5,6/5}(w_xd_\varphi^{5/6}, w_yd_\varphi^{5/6})$.  
The Region I Fourier-product identity \eqref{eq:repair-region-I-Fourier} and the preceding calculation yield
\[
 \tcF_{3/2}(v)=\kappa \hat{w}\hd^{1/3},\qquad
 \|\tcF_{3/2}(v)\|_3=\kappa \mathcal{S}(x)^{-1/3}
\]
for some $\kappa>0$.  
Region I Young inequality and the Hausdorff--Young inequality, we see that $\|v\|_{3/2}=\kappa \mathcal{S}(x)^{-1/3}$.  
Pairing against $\hat{w} \hd^{2/3}$ and using the Fourier duality and H\"older inequality gives
\[
 \kappa\mathcal{S}(x)^{-1}\le \kappa \mathcal{S}(x)^{-1/3}
 \bigl\|\widehat{\tcF}_{3/2}(\hat{w} \hd^{2/3})\bigr\|_3
 \le \kappa \mathcal{S}(x)^{-1}.
\]
Consequently the inverse Hausdorff--Young inequality is sharp at exponent $3/2$ for this input.

For $0\le t\le1/2$ set $f(t)=\|\hat\lambda(\hat{w}\hat\varphi)d_\varphi^t\|_{1/t}$, with
$f(0)$ the operator norm.  
Interpolation between $L^\infty$ and $L^2$ gives $\log f(t)$ is convex.  
The inverse Hausdorff--Young inequality gives $f(t)\le \mathcal{S}(x)^{t-1}$, and the preceding equality gives equality at $t=1/3$.  
A convex function bounded above by an affine function and touching it at an interior point agrees with it throughout the interval.  
Therefore
\[
 \|\widehat{\tcF}_{\tilde{p}}(\hat{w} \hd^{1/\tilde{p}})\|_{\tilde{p}'}=\mathcal{S}(x)^{-1/\tilde{p}},
 \qquad 1\le \tilde{p}\le2.
\]
Taking $\tilde{p}=r'$ proves equality for the canonical Region II convolution.
Both inputs are scaling-invariant, so Corollary~\ref{cor:hex-equality}
gives equality at every point of $\Omega_{p,s}$.
\end{proof}

\subsection{Young Region I}
We denote $T(t)\omega=\omega\tau_{-t}$ for real $t$ and $\omega\in L^1(\bG)$.
For $x\varphi\in L^1(\bG)$, we have that $T(t)(x\varphi)=\tau_t(x)\varphi$. 
The analytic continuations $T_p(it)$ below are first used only on a joint analytic core.

\begin{theorem}\label{thm:YF-low-output-family}
Suppose that $\bG$ is a locally compact quantum group with scaling constant $\nu=1$, $1\leq p,s,r\leq2$ and $1+1/r=1/p+1/s$.
For every $\alpha,\beta\in[0,1]$ there is a contractive bilinear map
\[
 D^{\alpha,\beta}_{p,s}:L^p(\phi)\times L^s(\phi) \longrightarrow L^r(\phi)
\]
defined as follows:
on the joint analytic core,
\begin{equation}\label{eq:YF-low-core-formula}
 D^{\alpha,\beta}_{p,s}(x,y)
 =T_p\bigl(i\alpha(1/p-1/r)\bigr)x *\rho_{-i/p'}\!\left(T_s\bigl(i\beta(1/s-1/r)\bigr)y\right).
\end{equation}
Explicitly,
\begin{align}\label{eq:scalingyoung2}
    \| D^{\alpha,\beta}_{p,s}(x,y)\|_{r} \leq \|x\|_p \|y\|_s.
\end{align}
This expression means the bounded extension of the whole bilinear
formula.

At $\alpha=\beta=0$ it is $\mathsf C^{\mathrm I}_{p,s}$.
At $r=2$ it is exactly $\mathsf C^{\alpha,\beta,0}_{p,s}$ from
\eqref{twYoung:definition}. 
At $r=1$ all scaling parameters $\alpha, \beta$ vanish.
\end{theorem}

\begin{proof}
The case $r=1$ is the predual convolution corner $p=s=1$.
For $1<r\leq2$ put
\[
 \vartheta=2(1-1/r),\qquad
 \frac1{p_0}=1-\frac{1-1/p}{\vartheta},\qquad
 \frac1{s_0}=1-\frac{1-1/s}{\vartheta}.
\]
Then $1\leq p_0,s_0\leq2$ and
$1/p_0+1/s_0=3/2$.  Set
\[
 \alpha_0=\alpha(1/p_0-1/2),\qquad
 \beta_0=\beta(1/s_0-1/2).
\]
We interpolate between $(1,1,1)$ and $(p_0,s_0,2)$.
We shall construct an analytic family in $\mathcal I_\varphi$, with norm
$\|\omega\|+\|\xi_2(\omega)\|_2$.
Both $T(t)$ and $\rho_t$ act isometrically on its two endpoint norms and commute. 
For $T$ this uses scaling constant one; 
for $\rho$ use $\lambda\rho_t=\hat\sigma_t\lambda$ and invariance of $\hat\varphi$ under its modular group.  
Joint Gaussian smoothing therefore supplies a common dense core $\mathcal I_{\varphi}^{a}$, entire in these two graph norms for both groups $\tau, \rho$.
The convolution of two elements of $\mathcal I_\varphi$ lies in $\mathcal I_\varphi$: 
besides the predual bound, the standard identity
\[
 \xi_2(\omega*\eta)=\lambda(\omega)\xi_2(\eta)
\]
gives its bound on the Hilbert space.  
In particular convolution is a continuous bilinear map $\mathcal I_\varphi\times \mathcal I_\varphi\to \mathcal I_\varphi$.

For $\omega,\eta\in \mathcal I_{\varphi}^{a}$ , we define
\begin{equation}\label{eq:YF-analytic-bilinear-family}
 \mathscr C_z(\omega,\eta)
   =(T(i\alpha_0 z)\omega) *\rho_{-iz/p_0'}\bigl(T(i\beta_0 z)\eta\bigr),
       \qquad 0\leq\Re z\leq1.
\end{equation}
It is analytic with values in $\mathcal I_\varphi$, hence in the compatible target sum.  
Its norm on $\mathcal I_\varphi$ is bounded on the strip for fixed $\omega,\eta$:
the imaginary parts of the group times stay in fixed compact intervals and the remaining real actions are isometries.

At $z=it$, we have that
\[
 \|\mathscr C_{it}(\omega,\eta)\|_1 \leq\|\omega\|_1\|\eta\|_1.
\]
At $z=1+it$, we put
\[
 \omega_t=T(-\alpha_0 t)\omega,\qquad \eta_t=\rho_{t/p_0'}T(-\beta_0 t)\eta.
\]
Their norms at $p_0,s_0$ are unchanged.  
The common $L^2$ representative of \eqref{eq:YF-analytic-bilinear-family} is precisely $\mathsf C^{\alpha,\beta,0}_{p_0,s_0}(\omega_t,\eta_t).$
Indeed the predual Fourier coefficient of \eqref{eq:YF-analytic-bilinear-family} is $\hat\tau_{i\alpha_0}\lambda(\omega_t)\,\hat\sigma_{-i/p_0'}\!\left(\hat\tau_{i\beta_0}\lambda(\eta_t)\right).$
 
This is exactly the coefficient of the closed spatial product $\tcF_{p_0}^{(\alpha)}(\omega_t)\,\tcF_{s_0}^{(\beta)}(\eta_t)$ after the two dual densities are commuted together.  
Plancherel identity and Young's inequality for Region II thus give
\[
 \|\mathscr C_{1+it}(\omega,\eta)\|_2\leq\|\omega\|_{p_0}\|\eta\|_{s_0}.
\]
All assertions first hold on $\mathcal I_{\varphi}^{a}$ and then on the stated boundary spaces by the density of $\mathcal I_{\varphi}^{a}$.

The bilinear Stein interpolation theorem (Proposition \ref{prop:steininter}), using the fixed compatible couples $(L^1,L^{p_0})$, $(L^1,L^{s_0})$ and $(L^1,L^2)$,
now gives the contraction at $z=\vartheta$.  
The exponent identities are
\[
 \frac1p=1-\vartheta+\frac{\vartheta}{p_0},\quad
 \frac1s=1-\vartheta+\frac{\vartheta}{s_0},\quad
 \frac1r=1-\frac{\vartheta}{2}.
\]
Finally
\[
 \vartheta \alpha_0=\alpha(1/p-1/r),\qquad
 \vartheta \beta_0=\beta(1/s-1/r),\qquad
 \frac{\vartheta}{p_0'}=\frac1{p'}.
\]
These give \eqref{eq:YF-low-core-formula} after spatial transport.
The dense-core formula determines the extension uniquely. 

The $r=2$ identification was already verified in the boundary calculation.
\end{proof}

The scaling convolutions in Regions III, IV are now defined, without any unbounded operation on arbitrary inputs, by
\begin{align}
 \left\langle D^{\mathrm {IV};\alpha,\beta}_{p,s}(x,y),a\right\rangle
 &=\overline{\langle D^{\alpha,\beta}_{p,r'}(x^*,a),y\rangle},
 \label{eq:YF-twisted-IV}\\
 \left\langle D^{\mathrm {III};\alpha,\beta}_{p,s}(x,y),a\right\rangle
 &=\overline{\langle D^{\alpha,\beta}_{r',s}(a,y^*),x\rangle}.
 \label{eq:YF-twisted-III}
\end{align}
By the spatial duality, we have that the following scaling Young inequalities:
\begin{align}
  \left \| D^{\mathrm {IV};\alpha,\beta}_{p,s}(x,y)\right\|_r \leq & \|x\|_p \|y\|_s, \label{eq:scalingyoung4}\\
   \left \| D^{\mathrm {III};\alpha,\beta}_{p,s}(x,y)\right\|_r \leq & \|x\|_p \|y\|_s. \label{eq:scalingyoung3}
\end{align}
For orientation, the formal analytic-core descriptions are
\begin{align}
 D^{\mathrm {IV};\alpha,\beta}_{p,s}(x,y)
 &=T_r\bigl(i\beta(1/s-1/r)\bigr)
   \mathsf C^{\mathrm {IV}}_{p,s}
                 \bigl(T_p(-i\alpha/r)x,y\bigr),
 \label{eq:YF-IV-core-twists}\\
 D^{\mathrm {III};\alpha,\beta}_{p,s}(x,y)
 &=T_r\bigl(i\alpha(1/p-1/r)\bigr)
   \mathsf C^{\mathrm {III}}_{p,s}
                 \bigl(x,T_s(-i\beta/r)y\bigr).
 \label{eq:YF-III-core-twists}
\end{align}
At $(p,s,r)=(1,\infty,\infty)$ and $(\infty,1,\infty)$ all these scaling twists vanish, retaining the untwisted regional slice actions.  
On other boundary overlaps, one should keep the regional definitions explicit instead of identifying products by their exponent triple alone.

In the following, we shall investigate the equality for scaling Young inequalities for Regions I, III, IV.
We begin with the following lemma.
\begin{lemma}\label{lem:YF-analytic-equality}
Let $\mathcal X,\mathcal Y, \mathcal Z$ be strictly convex Banach spaces, and let $U,V,W$ be strongly continuous groups of linear isometries on them.
Suppose that $0\leq\alpha,\beta\leq1$ and $B_{\alpha,\beta}:\mathcal X\times \mathcal Y\to \mathcal Z$ are bilinear contractions with joint real covariance
\[
 B_{\alpha,\beta}(U_t(x),V_t(y)) =W_t( B_{\alpha,\beta}(x,y)).
\]
Assume a dense common entire core supplies separate analytic
continuations to the two closed parameter strips, with one of the
following pairs of vertical covariance rules, where $\kappa,\tilde\kappa\ne0$:
\begin{align*}
 \text{(I)}\quad
 B_{\alpha+it,\beta}(x,y) &=B_{\alpha,\beta}(U_{t\kappa}(x),y),&
 B_{\alpha,\beta+it}(x,y) &=B_{\alpha,\beta}(x,V_{\tilde \kappa t}(y));\\
 \text{(IV)}\quad
 B_{\alpha+it,\beta}(x,y)&=B_{\alpha,\beta}(U_{t\kappa}(x),y),&
 B_{\alpha,\beta+it}(x,y) &=W_{\tilde \kappa t}B_{\alpha,\beta}(x,y);\\
 \text{(III)}\quad
 B_{\alpha+it,\beta}(x,y) &=W_{\kappa t}(B_{\alpha,\beta}(x,y)),&
 B_{\alpha,\beta+it}(x,y)&=B_{\alpha,\beta}(x,V_{\tilde\kappa t}(y)).
\end{align*}
For nonzero $x,y$, the following are equivalent:
\begin{enumerate}
\item $\|B_{\alpha_0,\beta_0}(x,y)\|
        =\|x\|\|y\|$ for some $0<\alpha_0,\beta_0<1$;
\item $\|B_{0,0}(x,y)\|=\|x\|\|y\|$, and
      $U_t(x)=x$, $V_t(y)=y$ for every real $t$;
\item Equality holds for every $(\alpha,\beta)\in[0,1]^2$.
\end{enumerate}
In this case, we have that $B_{\alpha, \beta}(x, y)=B_{0,0}(x,y)$, and it is invariant under the action of $W$.
\end{lemma}
\begin{proof}
We first justify the analytic continuations for arbitrary inputs.
Choose entire-core approximants $x_n\to x$, $y_n\to y$.
For a fixed real $\beta$ and $z=\alpha+it$, the vertical covariance
rule, isometry, and real-parameter contractivity give
\begin{align*}
 &\|B_{z,\beta}(x_n,y_n)-B_{z,\beta}(x_m,y_m)\|\\
 &\qquad\leq
   \|x_n-x_m\|\,\|y_n\|
       +\|x_m\|\,\|y_n-y_m\|.
\end{align*}
This estimate is uniform on the entire closed strip.
Hence the core functions converge uniformly there to a continuous function, holomorphic in its interior. The limit agrees at real parameters with the already defined bounded family members, obeys the same
vertical covariance, and has norm at most $\|x\|\|y\|$.
The identical argument applies to the second parameter strip.
No imaginary source action has thereby been separately defined on an arbitrary input.

Suppose (1) holds, and choose a norm-one functional $\ell\in \mathcal Z^*$ which norms $B_{\alpha_0,\beta_0}(x,y)$ with positive value $\|x\|\|y\|$.
The scalar maximum principle in the first strip makes $\ell(B_{z,\beta_0}(x,y))=\|x\|\|y\|$. 
Strict convexity of $\mathcal Z$ makes the corresponding supporting face of its radius-$\|x\|\|y\|$ ball a singleton,
so $B_{z,\beta_0}(x,y)$ is constant. 
In particular every real $\alpha$ gives the same output $Z$ of norm $\|x\|\|y\|$. 
Apply the same argument in the second strip at each such $\alpha$. 
It follows that every real-parameter output equals $Z$, and $Z=B_{0,0}(x,y)$.

In case (I), the first vertical rule now gives
\[
 B_{\alpha_0,\beta_0}(U_t(x),y)=Z\qquad(t\in\mathbb R),
\]
since $\kappa\ne0$. 
Bilinearity and contractivity imply
\[
\|x\|\|y\|
=\left\|B_{\alpha_0,\beta_0}\left(\frac{x+U_t(x)}{2},y\right)\right\|
   \leq\left\|\frac{x+U_t(x)}{2}\right\|\|y\|
   \leq \|x\|\|y\|.
\]
Strict convexity of $\mathcal X$ therefore gives $U_t(x)=x$. 
The second rule and strict convexity of $\mathcal Y$ gives $V_ty=y$.

In case (IV), the same midpoint argument gives $U_t(x)=x$, while the second vertical rule gives $W_tZ=Z$.
Canonical covariance consequently gives
\[
 B_{0,0}(x,V_t(y))
   =B_{0,0}(U_t(x),V_t(y))=W_t(Z)=Z.
\]
The midpoint argument, now for the canonical contraction, gives
$V_ty=y$. Case (III) is symmetric. This proves (1)$\Rightarrow$(2).

Conversely, suppose the statement (2) holds. 
Joint covariance fixes every output under $W$, since both inputs are fixed. 
Each of the relevant vertical rules consequently says that its analytic strip function is constant in the imaginary direction. 
A holomorphic Banach-valued function with that property has derivative zero and is constant.
Applying this to both coordinates proves that every output equals
$B_{0,0}(x,y)$, giving (3). 
The implication (3)$\Rightarrow$(1)
is immediate. Canonical covariance fixes the common output.
\end{proof}

\begin{theorem}\label{thm:YF-low-outer-equality}
Suppose that $\bG$ is a locally compact quantum group with scaling constant $\nu=1$.
Let $x,y$ be nonzero and let the exponents $p, s, r$ be strictly between one and infinity in the relevant admissible Young region.
In Region I, take $1<p,s<r\le2$ and put
\[
 a_0=\frac1p-\frac1r=\frac1{s'}>0,\qquad
 b_0=\frac1s-\frac1r=\frac1{p'}>0.
\]
For each of the proved families
\[
 D^{\alpha,\beta}_{p,s},\qquad
 D^{\mathrm{IV};\alpha,\beta}_{p,s},\qquad
 D^{\mathrm{III};\alpha,\beta}_{p,s},
\]
the following are equivalent:
\begin{enumerate}
\item The equality of the corresponding scaling Young inequality holds for the fixed pair $(x,y)$ at some
      $(\alpha,\beta)\in(0,1)^2$;
\item The corresponding untwisted regional Young inequality is an equality, and
      $T_p(t)x=x$, $T_s(t)y=y$ for every real $t$;
\item The equality of the corresponding scaling Young inequality holds for every
      $(\alpha,\beta)\in[0,1]^2$.
\end{enumerate}
All outputs in the specified regional family then coincide with
its canonical output and are fixed by $T_r$.
\end{theorem}
\begin{proof}
The analytic core formulas give the vertical covariance rules
\begin{align*}
 D^{\alpha+it,\beta}(x,y)&=D^{\alpha,\beta}(T_p(-a_0t)x,y),&
 D^{\alpha,\beta+it}(x,y) &=D^{\alpha,\beta}(x,T_s(-b_0t)y);\\
 D^{\mathrm{IV};\alpha+it,\beta}(x,y)&=D^{\mathrm{IV};\alpha,\beta}(T_p(t/r)x,y),&
 D^{\mathrm{IV};\alpha,\beta+it}(x,y) &=T_r(-t/p')D^{\mathrm{IV};\alpha,\beta}(x,y);\\
 D^{\mathrm{III};\alpha+it,\beta}(x,y)&=T_r(-t/s')D^{\mathrm{III};\alpha,\beta}(x,y),&
 D^{\mathrm{III};\alpha,\beta+it}(x,y)&=D^{\mathrm{III};\alpha,\beta}(x,T_s(t/r)y).
\end{align*}
For the outer regions these follow from their proved core formulas
\begin{align*}
 D^{\mathrm{IV};\alpha,\beta}(x,y)
  &=T_r(i\beta/p')\,
       \mathsf C^{\mathrm{IV}}
                    (T_p(-i\alpha/r)x,y),\\
 D^{\mathrm{III};\alpha,\beta}(x,y)
  &=T_r(i\alpha/s')\,
       \mathsf C^{\mathrm{III}}
                    (x,T_s(-i\beta/r)y).
\end{align*}
Spatial $L^p,L^s,L^r$ spaces are strictly convex in this range.
Lemma~\ref{lem:YF-analytic-equality} therefore applies.
\end{proof}

\begin{lemma}\label{lem:boundary}
Let $\Omega=\{z\in\mathbb C:0<\operatorname{Re}z<1\}$, $0<\vartheta<1$, and $\epsilon>0$.
Suppose $F_n$ are bounded holomorphic functions on $\Omega$ such that
\[
 \sup_{z\in \Omega}|F_n(z)|\leq1+\varepsilon_n,\qquad
 \varepsilon_n\geq0,\quad
 \varepsilon_n\longrightarrow0,\qquad
 F_n(\vartheta)\longrightarrow1.
\]
Assume that
\[
 F_n(1+it)\longrightarrow f(t)
 \quad\text{for almost every }t\in(-\epsilon,\epsilon),
\]
where $f$ has a representative continuous at $t=0$.
Then $f(t)=1$ almost everywhere on this interval and $f(0)=1$.
Moreover, $F_n\to1$ locally uniformly in $\Omega$.
\end{lemma}

\begin{proof}
Put $G_n=F_n/(1+\varepsilon_n)$.
The Poisson representation of the bounded harmonic function $1-\operatorname{Re}G_n$ gives
\[
 \int_{-\epsilon}^{\epsilon}
 \frac{\sin(\pi\vartheta)}
      {2\bigl(\cosh(\pi t)+\cos(\pi\vartheta)\bigr)}
 \bigl(1-\operatorname{Re}G_n(1+it)\bigr)\,dt
 \leq 1-\operatorname{Re}G_n(\vartheta)\longrightarrow0.
\]
Fatou's lemma and strict positivity of the kernel imply $\operatorname{Re}f(t)=1$ almost everywhere.
Since $|f(t)|\leq1$ almost everywhere, $f(t)=1$ there.
Continuity at zero implies $f(0)=1$.
Finally, every subsequential locally uniform limit of $G_n$ has modulus at most one and takes the value one at $\vartheta$.
Montel's theorem and the maximum-modulus principle therefore give $G_n\to1$, and hence $F_n\to1$, locally uniformly.
\end{proof}

\begin{proposition}\label{lift:kac}
Let $\mathbb G$ be a unimodular locally compact quantum group. Let $1<p,s<r<2$ satisfy $1+1/r=1/p+1/s$, and put
\[
 p_0'=\frac{2p'}{r'},\qquad s_0'=\frac{2s'}{r'},\qquad
 \vartheta=\frac2{r'}.
\]
If $0\ne x\in L^p(\phi)$ and $0\ne y\in L^s(\phi)$ satisfy
\[
 \|\mathsf C^{\mathrm I}_{p,s}(x,y)\|_r=\|x\|_p\|y\|_s,
\]
then, writing $x=w_x|x|$ and $y=w_y|y|$, one has
\[
 \left\|\mathsf C^{\mathrm I}_{p_0,s_0}
 \bigl(w_x|x|^{p/p_0},w_y|y|^{s/s_0}\bigr)\right\|_2
 =\|x\|_p^{p/p_0}\|y\|_s^{s/s_0}.
\]
\end{proposition}
\begin{proof}
Normalize $\|x\|_p=\|y\|_s=1$ and set $v=\mathsf C^{\mathrm I}_{p,s}(x,y)=w_v|v|$, so $\|v\|_r=1$.
The exponent identities are
\[
 \frac1{p_0}+\frac1{s_0}=\frac32,\qquad
 \frac1p=1-\vartheta+\frac\vartheta{p_0},\quad
 \frac1s=1-\vartheta+\frac\vartheta{s_0},\quad
 \frac1r=1-\frac\vartheta2.
\]
In particular $1<p_0,s_0<2$. 
Unimodularity implies $\nu=1$.

Use the analytic family already constructed in the proof of Theorem~\ref{thm:YF-low-output-family}, with the scaling twists set to zero:
\begin{equation}\label{lift:correct-family}
 \mathscr C_z(\omega,\eta)
 =\omega*\rho_{-iz/p_0'}(\eta),\qquad0\le\Re z\le1.
\end{equation}
All formulas in this paragraph refer to the common compatible data on the joint analytic core. 
Its lower boundary is a contraction from $L^1\times L^1$ to $L^1$. 
Its upper boundary, under the spatial identifications, is
\[
 \mathscr C_{1+it}(a,b)
 =\mathsf C^{\mathrm I}_{p_0,s_0}\bigl(a,\rho_{t/p_0'}b\bigr).
\]
The real action in this formula is isometric, so the upper boundary is contractive from $L^{p_0}\times L^{s_0}$ to $L^2$, uniformly in $t$.
At $z=\vartheta$, the interpolated map is exactly
$\mathsf C^{\mathrm I}_{p,s}$, because $\vartheta/p_0'=1/p'$.
Thus the modular correction is retained throughout the interpolation.

Put $a(z)=1-z/p_0'$ and $b(z)=1-z/s_0'$. Use the compatible polar interpolation sections whose ambient representatives are
\[
 x(z)=[w_x|x|^{p a(z)}d_\varphi^{1-a(z)}],\qquad
 y(z)=[w_y|y|^{s b(z)}d_\varphi^{1-b(z)}],
\]
and the dual section with ambient representative
\[
 v(z)=[w_v|v|^{rz/2}d_\varphi^{-z/2}].
\]
These expressions specify sections in the compatible ambient spaces, rather than operators belonging to one fixed homogeneous spatial degree for every $z$.
Their physical representatives at $\vartheta$ are $x,y,N_r(v)$; at $1$ they are
\[
 x_1=w_x|x|^{p/p_0},\qquad y_1=w_y|y|^{s/s_0},\qquad
 v_1=w_v|v|^{r/2}.
\]
The boundary norms of $x(z),y(z)$ are at most one in the input couples
$(L^1,L^{p_0})$, $(L^1,L^{s_0})$; those of the dual section are at most one in $(L^\infty,L^2)$. 
These are the standard compatible polar sections of spatial interpolation; see \cite{Ter81,Izu97,Cas13}. 
On the upper boundary the first physical representative is
\[
 w_x|x|^{p/p_0-ipt/p_0'}d_\varphi^{it/p_0'},
\]
with the analogous formula for $y$. 
The relative cocycle, with its support understood, is strongly continuous, so these boundary sections are norm-continuous in their finite endpoint spaces. 


There are $\epsilon>0$ and admissible analytic core sections $x_n,y_n,v_n$ with the following properties.
\begin{enumerate}
\item Their scalar pairings
\[
 F_n(z)=
 \bigl\langle \mathscr C_z(x_n(z),y_n(z)),v_n(\bar z)\bigr\rangle
\]
are bounded holomorphic functions on $\{z: 0< \Re z <1\}$, have boundary values, and satisfy the admissible growth condition for Hadamard three-lines theorem.
On both boundaries, the products of the three relevant endpoint norms
are at most $1$.
The endpoint triples are
\[
 (L^1,L^1,L^\infty),\qquad
 (L^{p_0},L^{s_0},L^2).
\]
Here the conjugate parameter in the sesquilinear pairing is understood in the compatible dual-section convention.

\item At $\vartheta$, 
\[
 x_n(\vartheta)\to x,\qquad
 y_n(\vartheta)\to y,\qquad
 v_n(\vartheta)\to N_r(v)
\]
in $L^p,L^s,L^{r'}$, respectively.

\item For almost every $|t|<\epsilon$,
\[
 x_n(1+it)\to x_1(t),\quad
 y_n(1+it)\to y_1(t),\quad
 v_n(1-it)\to v_1(t)
\]
in $L^{p_0},L^{s_0},L^2$, respectively.
These three limiting sections are norm-continuous at zero and satisfy
\[
 x_1(0)=w_x|x|^{p/p_0},\quad
 y_1(0)=w_y|y|^{s/s_0},\quad
 v_1(0)=w_v|v|^{r/2}.
\]
\end{enumerate}
In fact, we can take 
\begin{align*}
    x_n(z)= & \int_{\bR} g_n(\tilde{t}\,)T(\tilde{t}\,)\int_{\bR}g_n(t)x(z+it) dt d\tilde{t}, \\
y_n(z)= & \int_{\bR} g_n(\tilde{t}\,)T(\tilde{t}\,)\int_{\bR}g_n(t)y(z+it) dt d\tilde{t},\\
v_n(z)= & \int_{\bR} g_n(\tilde{t}\,)\tau_{\tilde{t}\,}\int_{\bR}g_n(t)v(z+it) dt d\tilde{t},
\end{align*}
where $\displaystyle g_n(t)=\sqrt{\frac{n}{\pi}} e^{-nt^2}$.
We shall regard $x(z)$, $y(z)$ on the enlarged strip $0< \Re z < p_0'/2$ and $0< \Re z < s_0'/2$ respectively.
We have that $\displaystyle \sup_{0< \Re z <1}|F_n|\leq1$.
The interior norm convergence and the contraction of $\mathsf C_{p,s}^{\mathrm I}$ give
\[
 F_n(\vartheta)\longrightarrow
 \langle v,N_r(v)\rangle=1.
\]
By the upper-boundary norm convergence just proved and the bounded bilinear endpoint map,
\[
 F_n(1+it)\longrightarrow
 \bigl\langle \mathsf C_{p_0,s_0}^{\mathrm I}\bigl(x_1(t),\rho_{t/p_0'} y_1(t)\bigr),v_1(t)\bigr\rangle=f(t)
\]
for almost every $|t|<\epsilon$.
The real action $\rho$ is strongly continuous and isometric on
$L^{s_0}$; hence $f$ is continuous at zero.
Lemma~\ref{lem:boundary} gives $f(0)=1$.
This proves the normalized assertion. 
Rescaling gives the formula in the statement.


\end{proof}

\begin{theorem}\label{main:region1}
Suppose that $\bG$ is a unimodular locally compact quantum group.
Let $1<p,s<r\le2$ and $1/r+1=1/p+1/s$. 
For nonzero $x\in L^p(\phi)$ and $y\in L^s(\phi)$, equality in the scaling Young inequality at one $(\alpha,\beta)\in(0,1)^2$ is equivalent to $(x,y)$ being a scaling-invariant bi-shift $(p,s)$-pair of a group-like projection. 
Such a pair attains equality at every point of $[0,1]^2$.
\end{theorem}
\begin{proof}
For $r=2$, this is Theorem~\ref{twYoung:structural}, since the two regional products agree there. 
Let $r<2$. Theorem~\ref{thm:YF-low-outer-equality} gives scaling invariance and canonical Young equality. 
Proposition~\ref{lift:kac} lifts this equality to output exponent $2$ for
\[
 \widetilde{x}=w_x|x|^{p/p_0},\qquad \widetilde{y}=w_y|y|^{s/s_0}.
\]
Both remain scaling-invariant. 
Their forward Fourier images are fixed by covariance, so all their forward twists agree. 
Theorem \ref{twYoung:structural} at output $2$ therefore makes $\widetilde{x},\widetilde{y}$ a bi-shift pair of a group-like projection with matching Fourier supports.

To pass back to $x$, write $\widetilde{x}=\kappa w d_\varphi^{1/p_0}$, with $\kappa>0$ and $w^*w=|w|\in\mathcal C_\varphi$.
Its polar part is $w_x=w$ and
\[
 |x|^{p/p_0}=\kappa |w| d_\varphi^{1/p_0},\qquad
 x=\kappa^{p_0/p}w d_\varphi^{1/p}.
\]
This follows by functional calculus because $|w|$ commutes with $d_\varphi$. 
The same argument applies to $y$. 
This proves the necessity of the theorem.

Conversely, suppose $(x,y)$ is a scaling-invariant bi-shift pair of a group-like projection. 
Write its common size as $S$, its Fourier polar parts as $w_{\hat{x}},w_{\hat{y}}$, and its positive coefficient magnitudes as $\mu_x,\mu_y$. 
As in Theorem~\ref{twYoung:structural}, $w_{\hat{x}},w_{\hat{y}}$ belong to the dual Haar centralizer, $w_{\hat{x}}^*w_{\hat{x}}=w_{\hat{y}}w_{\hat{y}}^*$, and $\hat{w}=w_{\hat{x}}w_{\hat{y}}$ has right support of dual Haar weight $\mathcal{S}(x)^{-1}$. 
The Region I Fourier-product identity \eqref{eq:repair-region-I-Fourier} gives
\[
 \tcF_r(\mathsf C^{\mathrm I}_{p,s}(x,y))
 =\mu_x\mu_y \mathcal{S}(x)^2 \hat{w} \hd^{1/r'}.
\]
Its $L^{r'}$ norm is $\mu_x\mu_y \mathcal{S}(x)^{2-1/r'}=\|x\|_p\|y\|_s$. The forward Hausdorff--Young inequality and Young's upper bound squeeze the output norm to that same value.
Scaling invariance removes both input scaling parameters, proving equality throughout the rectangle.
\end{proof}

\subsection{Young Regions III and IV}

In this section, we shall investigate the equality of scaling Young inequality for unimodular locally compact quantum groups in Young Regions III and IV.

\begin{theorem}\label{thm:outer-structural-parametrization}
Suppose that $\mathbb G$ is a unimodular locally compact quantum group and that $1+1/r=1/p+1/s$.  
Let $0<\alpha,\beta<1$. 
All elements below are assumed nonzero.

In Region IV, suppose $1<p<2<s<r<\infty$ and $\|x\|_p=\|y\|_s=1$.  Equality in \eqref{eq:scalingyoung4} holds if and only if there is a normalized scaling-invariant bi-shift $(p,r')$-pair
$(a,b)\in L^p(\phi)\times L^{r'}(\phi)$ of a group-like projection such that
\[
 x=a^*,\qquad y=N_{s'}(D^{0,0}_{p,r'}(a,b)).
\]
In this case, $\|D^{0,0}_{p,r'}(a,b)\|_{s'}=1$ and $D^{\mathrm{IV};\alpha,\beta}_{p,s}(x,y)=N_{r'}(b)$.

In Region III, suppose $1<s<2<p<r<\infty$ and
$\|x\|_p=\|y\|_s=1$.  Equality in
\eqref{eq:scalingyoung3} holds if and only if there is a normalized
scaling-invariant bi-shift $(r',s)$-pair $(a,b)\in L^{r'}(\phi)\times L^s(\phi)$ such that
\[
 x=N_{p'}(D^{0,0}_{r',s}(a,b)),\qquad y=b^*.
\]
In this case,  $\|D^{0,0}_{r',s}(a,b)\|_{p'}=1$ and  $D^{\mathrm{III};\alpha,\beta}_{p,s}(x,y)=N_{r'}(a)$.

In both cases equality holds for every $(\alpha,\beta)\in[0,1]^2$, all outputs in the corresponding family coincide, and the two inputs and output are invariant under their real scaling actions.
\end{theorem}

\begin{proof}
For Region IV the conjugate exponents satisfy
$1<p,r'<s'<2$ and $1+1/s'=1/p+1/r'$, so the auxiliary
product lies in the strict Region I range.  Suppose equality holds,
and let $v=D^{\mathrm{IV};\alpha,\beta}_{p,s}(x,y)$ and
$b=N_r(v)$.  Then $\|b\|_{r'}=1$ and
$\langle v,b\rangle=1$.  Put
\[
 a=x^*,\qquad v_{\alpha,\beta}=D^{\alpha,\beta}_{p,r'}(a,b).
\]
The transpose identity \eqref{eq:YF-twisted-IV}, Young's inequality
and H\"older's inequality give
\[
 1=\overline{\langle v_{\alpha,\beta},y\rangle}
 \le \|v_{\alpha,\beta}\|_{s'}\le1.
\]
Thus $(a,b)$ is an equality pair for the strict Region I inequality.
Theorem~\ref{main:region1} identifies it as a scaling-invariant bi-shift pair of a group-like projection.  
Theorem~\ref{thm:YF-low-outer-equality} gives $v_{\alpha,\beta}=D^{0,0}_{p,r'}(a,b)$.  
Equality in the norming pairing, including its phase, implies $y=N_{s'}(D^{0,0}_{p,r'}(a,b))$; and $v=N_{r'}(b)$ follows from $b=N_r(v)$.

Conversely, choose such a normalized pair $(a,b)$ and put
$x=a^*$ and $y=N_{s'}(D^{0,0}_{p,r'}(a,b))$.  
Region I structure gives $\|D^{0,0}_{p,r'}(a,b)\|_{s'}=1$ and $D^{\alpha,\beta}_{p,r'}(a,b)=D^{0,0}_{p,r'}(a,b)$ for every parameter pair in
the square. 
Evaluating \eqref{eq:YF-twisted-IV} against $b$ gives
\[
 \big\langle D^{\mathrm{IV};\alpha,\beta}_{p,s}(x,y),b\big\rangle
 =\overline{\langle D^{0,0}_{p,r'}(a,b),N_{s'}(D^{0,0}_{p,r'}(a,b))\rangle}=1.
\]
The scaling Young inequality and the unit norm of $b$ imply equality; 
the uniqueness of the norming vector then gives the output $N_{r'}(b)$.

For Region III, $1<r',s<p'<2$ and
$1+1/p'=1/r'+1/s$.  If equality holds, set
$v=D^{\mathrm{III};\alpha,\beta}_{p,s}(x,y)$,
$a=N_r(v)$, $b=y^*$, and
$v_{\alpha,\beta}=D^{\alpha,\beta}_{r',s}(a,b)$.
Equation \eqref{eq:YF-twisted-III} gives
\[
 1=\overline{\langle v_{\alpha,\beta},x\rangle}
 \le \|v_{\alpha,\beta}\|_{p'}\le1.
\]
Theorem \ref{main:region1} shows the $(a,b)$ is a bishift pair of group-like projection and $v_{\alpha,\beta}=D^{0,0}_{r',s}(a,b)$; uniqueness of norming vectors gives
$x=N_{p'}(D^{0,0}_{r',s}(a,b))$ and $v=N_{r'}(a)$.  
Conversely, the indicated pair and these assignments give unit pairing with $a$, and hence equality and the claimed output for every parameter pair $(\alpha, \beta)$.

Real scaling actions preserve polar decompositions and norming maps $N_p$, and the Region I products are scaling covariant.  
The auxiliary pairs are scaling-invariant, so is $D^{0,0}_{r',s}(a,b)$ and the reconstructed inputs and output are scaling-invariant as well. 
This proves the theorem.
\end{proof}

\subsection{Compatibility on the common edges}\label{sec:edge-compatibility}
The regional subdivision is
\[
 \mathrm I:\ 1/r\ge\tfrac12,\qquad
 \mathrm {II}:\ 1/p,1/s\ge\tfrac12,\ 1/r\le\tfrac12,\qquad
 \mathrm {III}:\ 1/p\le\tfrac12,\qquad
 \mathrm {IV}:\ 1/s \le\tfrac12.
\]
All these conditions are restricted to the exponent triangle.
The nondegenerate shared edges and remaining intersections are
\begin{center}
\begin{tabular}{lll}
\toprule
Regions & Intersection & Original maps \\
\midrule
I--II & $r=2$, $1/p+1/s=3/2$ & Agree \\
II--III & $p=2$, $1\le s\le2$ & Generally differ \\
II--IV & $s=2$, $1\le p\le2$ & Generally differ \\
I--III & $(p,s,r)=(2,1,2)$ & Generally differ \\
I--IV & $(p,s,r)=(1,2,2)$ & Generally differ \\
III--IV & $(p,s,r)=(2,2,\infty)$ & Generally differ \\
\bottomrule
\end{tabular}
\end{center}
Here agreement means equality of bilinear maps on the same inputs.

\begin{proposition}[Exact boundary formulas]\label{prop:edge-formulas}
Suppose that $\bG$ is a locally compact quantum group with scaling constant $\nu=1$. 
On the common I--II edge,
\begin{equation}\label{edge:I-II}
 D^{\alpha,\beta}_{p,s}(x,y)
 =\mathsf C^{\alpha,\beta,\gamma}_{p,s}(x,y),\qquad r=2,
\end{equation}
for every $\gamma\in[0,1]$.
On the II--IV edge the original transpose has the formula
\begin{equation}\label{edge:original-IV}
 D^{\mathrm{IV};\alpha,\beta}_{p,2}(x,y)
 =\mathcal G_{r'}^{(\beta)}
 \left(\bigl(\tcF_p^{(\alpha)}(x^*)\bigr)^*\tcF_2(y)\right).
\end{equation}
On the II--III edge it has the formula
\begin{equation}\label{edge:original-III}
 D^{\mathrm{III};\alpha,\beta}_{2,s}(x,y)
 =\mathcal G_{r'}^{(\alpha)}
 \left(\tcF_2(x)\bigl(\tcF_s^{(\beta)}(y^*)\bigr)^*\right).
\end{equation}
In contrast, the Region II boundary values are
\begin{align}
 \mathsf C^{\theta_1,\theta_2,\theta_3}_{p,2}(x,y)
 &=\mathcal G_{r'}^{(\theta_3)}
   \left(\tcF_p^{(\theta_1)}(x)\tcF_2(y)\right),\label{edge:II-IV}\\
 \mathsf C^{\theta_1,\theta_2,\theta_3}_{2,s}(x,y)
 &=\mathcal G_{r'}^{(\theta_3)}
   \left(\tcF_2(x)\tcF_s^{(\theta_2)}(y)\right).
       \label{edge:II-III}
\end{align}
The coordinate attached to an input exponent $2$ is ineffective.
\end{proposition}
\begin{proof}
Equation~\eqref{edge:I-II} is the output-$2$ identification in Theorem~\ref{thm:YF-low-output-family}; the inverse transform is Plancherel and $\kappa_2=0$.
To prove \eqref{edge:original-IV}, test against $h\in L^{r'}(\phi)$ and use the defining transpose identity. The auxiliary Region I product has output $2$, hence its Fourier transform is
\[
 \tcF_p^{(\alpha)}(x^*)\tcF_{r'}^{(\beta)}(h).
\]
Taking the conjugate of its pairing with $\tcF_2(y)$, and cyclically rearranging the spatial integral, gives
\[
 \left\langle
   \bigl(\tcF_p^{(\alpha)}(x^*)\bigr)^*\tcF_2(y),
   \tcF_{r'}^{(\beta)}(h)\right\rangle.
\]
The adjoint definition of $\mathcal G$ gives the asserted formula. 
The spatial H\"older inequality ensures that the product lies in $L^{r'}$. 
The Region III formula follows in the same way, preserving the order of the two factors. 
Density gives the identities for all inputs. 
Equations~\eqref{edge:II-IV}--\eqref{edge:II-III} are immediate from the Region II definition.
\end{proof}

In the following, we shall unify the Young regions.
\begin{theorem}\label{thm:compatible-outer-products}
Suppose that $\bG$ is a locally compact quantum group with scaling constant $\nu=1$ and recall $J_p$ from Proposition~\ref{prop:symmetry}. We define the convolution $\mathsf E$ as follows:
\begin{align}
 \langle\mathsf E^{\mathrm{IV};\alpha,\beta}_{p,s}(x,y),h\rangle
 &=\overline{\left\langle
   D^{1-\alpha,\beta}_{p,r'}(J_p x,h),y\right\rangle},
       \label{edge:new-IV}\\
 \langle\mathsf E^{\mathrm{III};\alpha,\beta}_{p,s}(x,y),h\rangle
 &=\overline{\left\langle
   D^{\alpha,1-\beta}_{r',s}(h,J_s y),x\right\rangle}.
       \label{edge:new-III}
\end{align}
They are bounded bilinear contractions on their respective regions. 
Their shared-edge formulas are
\begin{align}
 \mathsf E^{\mathrm{IV};\alpha,\beta}_{p,2}(x,y)
 &=\mathcal G_{r'}^{(\beta)}
       \left(\tcF_p^{(\alpha)}(x)\tcF_2(y)\right),
       \label{edge:repaired-IV}\\
 \mathsf E^{\mathrm{III};\alpha,\beta}_{2,s}(x,y)
 &=\mathcal G_{r'}^{(\alpha)}
       \left(\tcF_2(x)\tcF_s^{(\beta)}(y)\right).
       \label{edge:repaired-III}
\end{align}
Consequently the following piecewise family is well defined on the entire closed Young triangle for every $\boldsymbol\theta\in[0,1]^3$:
\begin{equation}\label{edge:global-family}
 \mathsf E^{\boldsymbol\theta}_{p,s}=
 \begin{cases}
 D^{\theta_1,\theta_2}_{p,s},&\mathrm I,\\
 \mathsf C^{\theta_1,\theta_2,\theta_3}_{p,s},&\mathrm {II},\\
 \mathsf E^{\mathrm{III};\theta_3,\theta_2}_{p,s},&\mathrm {III},\\
 \mathsf E^{\mathrm{IV};\theta_1,\theta_3}_{p,s},&\mathrm {IV}.
 \end{cases}
\end{equation}
Its norm from $L^p(\phi)\times L^s(\phi)$ to $L^r(\phi)$ is at most one.
\end{theorem}
\begin{proof}
The maps $J_p$ are conjugate-linear surjective isometries. 
The same duality argument as for the original transposes therefore defines bounded bilinear contractions in \eqref{edge:new-IV}--\eqref{edge:new-III}, including the endpoint $r=\infty$.
On the II--IV edge, repeat the proof of \eqref{edge:original-IV}. 
The first factor is now
\[
 \bigl(\tcF_p^{(1-\alpha)}(J_p x)\bigr)^*
 =\tcF_p^{(\alpha)}(x),
\]
by reflection and $\nu=1$. This proves \eqref{edge:repaired-IV}; 
the III formula follows from the analogous identity with $s$ and $\beta$. 
They identify the II--IV parameters as $(\theta_1,\theta_3)=(\alpha,\beta)$, and the II--III parameters as $(\theta_2,\theta_3)=(\beta,\alpha)$.

Equation~\eqref{edge:I-II} handles the I--II edge. 
At $(1,2,2)$ only $\theta_1$ is effective and all three incident definitions agree.
At $(2,1,2)$ only $\theta_2$ is effective. 
At $(2,2,\infty)$ only $\theta_3$ is effective: both new outer formulas equal
\[
 \mathcal G_1^{(\theta_3)}\bigl(\tcF_2(x)\tcF_2(y)\bigr).
\]
These exhaust all intersections, proving that \eqref{edge:global-family} is independent of the regional choice. 
Each piece is contractive, so the global indexed family is contractive.
\end{proof}

\begin{remark}[Transport of equality statements]
The original products and their equality theorems remain the ones defined by \eqref{eq:YF-twisted-IV}--\eqref{eq:YF-twisted-III}. 
The compatible products satisfy the exact identities
\begin{align*}
 \mathsf E^{\mathrm{IV};\alpha,\beta}_{p,s}(x,y)
 &=D^{\mathrm{IV};1-\alpha,\beta}_{p,s}((J_p x)^*,y),\\
 \mathsf E^{\mathrm{III};\alpha,\beta}_{p,s}(x,y)
 &=D^{\mathrm{III};\alpha,1-\beta}_{p,s}(x,(J_s y)^*).
\end{align*}
The input changes here are linear surjective isometries. 
Thus sharp norms and equality statements transport exactly with these input and parameter changes. 
For example, in the normalized unimodular structural parameterization, the new IV equality pairs have $x=J_p a$, $y=N_{s'}(c)$, and the new III pairs have $x=N_{p'}(c)$, $y=J_s b$, with the same Region I bi-shift data as before. 
Edge agreement implies agreement of the norm and equality sets for the same inputs after the parameter identifications in \eqref{edge:global-family}. 
The structural theorems restricted to strict interiors are not, by themselves, endpoint classification theorems.
\end{remark}

\section{Nikolskii inequalities}\label{sec:nikolskii}

The Nikolskii inequality was first studied by S. M. Nikolskii in \cite{Nik51}.
In \cite{AkyRuz20}, Akylzhanov and Ruzhansky study the Nikolskii inequality for locally compact groups.
In this section, we generalize the Nikolskii inequality to locally compact quantum groups.

The right support of a spatial $L^p$ element can depend on its interpolation exponent. 
We therefore specify the representative used in every Fourier support below. 
If $\omega$ belongs to the compatible intersection $L^p(\bG)\cap L^2(\bG)$, we write $\xi_p(\omega)$ and
$\xi(\omega)$ for its two representatives. 
This extends the notation already used for $\omega\in L^1(\bG)\cap L^2(\bG)$.
For $\xi\in L^2(\bG)$ put
\[
 Y_\xi=\tcF_2(\Phi_2^{-1}(\xi)),\qquad
 q_\xi=\cR(|\tcF_2(\Phi_2^{-1}(\xi))|),\qquad
 \hcS(\xi)=\hat\varphi(q_\xi).
\]
Thus $\hcS$ in this section always uses the $L^2$ Fourier
representative. In particular, no identification with
$\cS(\tcF_p(\Phi_p^{-1}(\xi_p(\omega))))$ is implicit.
All products of spatial $L^p$ elements below are their closed products.

\begin{proposition}\label{prop:nik-inverse-endpoint}
Let $\xi\in L^2(\bG)$ and suppose that
$K=\hcS(\xi)<\infty$. 
Then there is a unique $x\in\mathfrak N_\varphi$ representing $\xi$, and
\[
 \xi=\Lamv(x),\qquad
 x=\hat\lambda(\hat\omega)\quad\hbox{for some }
 \hat\omega\in L^1(\hbG)\cap L^2(\hbG),
 \qquad
 \|x\|_\infty\le \hcS(\xi)^{1/2}\|\xi\|_2.
\]
\end{proposition}

\begin{proof}
Write $\widehat{\xi}=\tcF_2(\Phi_2^{-1}(\xi))$. 
Since $\hat\varphi(q_\xi)<\infty$, $q_\xi\hd^{1/2}\in L^2(\hat\phi)$ with norm $\hat\varphi(q_\xi)^{1/2}$. 
The product $\widehat{\xi} q_\xi \hd^{1/2}\in L^1(\hat\phi)$ defines a normal functional $\hat\omega$ on $\hcM$ and satisfies
\[
 \|\hat\omega\|_1=\|\widehat{\xi} q_\xi \hd^{1/2}\|_{1,\hat\phi}
 \le\|\widehat{\xi}\|_{2,\hat\phi}\hat\varphi(q_\xi)^{1/2}.
\]
For $a\in\mathfrak N_{\hat\varphi}$, the spatial pairing and $ \widehat{\xi} q_\xi=\widehat{\xi}$
give
\[
 \hat\omega(a^*)
 =\int\widehat{\xi} q_\xi\hd^{1/2}a^*\,d\hat\phi
 =\big\langle\hat\Phi_2(\widehat{\xi}),\hLam(a)\big\rangle.
\]
Consequently $\hat\omega\in\mathcal I_{\hat\varphi}$ and $\xi_{\hat\varphi}(\hat\omega)=\hat\Phi_2(\widehat{\xi})$ by the defining GNS criterion for $\mathcal I_{\hat\varphi}$.
Inverse Plancherel and the defining formula for the Fourier transform on $L^1\cap L^2$ now yield
\[
 \Lamv(\hat\lambda(\hat\omega))
 =\widehat{\cF}_2(\xi_{\hat\varphi}(\hat\omega))
 =\xi.
\]
Set $x=\hat\lambda(\hat\omega)$. 
The Fourier representation is contractive from $L^1$ into $L^\infty$, so
\[
 \|x\|_\infty\le\|\hat\omega\|_1  \le \hcS(\xi)^{1/2}\|\xi\|_2.
\]
The inclusion map is injective because $\varphi$ is faithful, proving uniqueness. 
\end{proof}

\begin{theorem}[Nikolskii inequality]\label{thm:nikolskii}
Let $1\le p\le2$, and let $\omega$ be a nonzero common element of
$L^p(\bG)\cap L^2(\bG)$. 
Suppose that $\hcS(\xi(\omega))<\infty$.
Then $\omega$ has compatible representatives in $L^s(\bG)$ for every
$p\le s\le\infty$, and
\begin{equation}\label{eq:nikolskii-correct}
 \|\xi_s(\omega)\|_s
 \le \hcS(\xi(\omega))^{1/p-1/s}\|\xi_p(\omega)\|_p,
 \qquad p\le s\le\infty.
\end{equation}
Here $\xi_\infty(\omega)$ denotes the bounded representative obtained
in Proposition~\ref{prop:nik-inverse-endpoint}.
\end{theorem}

\begin{proof}

Proposition~\ref{prop:nik-inverse-endpoint} gives $\|\xi_\infty(\omega)\|_\infty \le \hcS(\xi(\omega))^{1/2}\|\xi(\omega)\|_2$.
If $p<2$, the log-convexity of compatible interpolation norms gives
\[
 \|\xi(\omega)\|_2\le \|\xi_\infty(\omega)\|_\infty^{1-p/2}\|\xi_p(\omega)\|_p^{p/2}.
\]
Combining these two estimates and dividing by $\|\xi(\omega)\|_2^{1-p/2}$ (which is nonzero) yields
\begin{equation}\label{eq:nik-intermediate}
 \|\xi(\omega)\|_2\le \hcS(\xi(\omega))^{1/p-1/2}\|\xi_p(\omega)\|_p,
 \qquad \|\xi_\infty(\omega)\|_\infty \le \hcS(\xi(\omega))^{1/p}\|\xi_p(\omega)\|_p.
\end{equation}
For $p=2$ the second estimate follows directly from the endpoint proposition and the first one is an identity.
Interpolation between the compatible $L^p$ and $L^\infty$ representatives gives, for every finite $s\ge p$,
\[
 \|\xi_s(\omega)\|_s
 \le \|\xi_\infty(\omega)\|_\infty^{1-p/s}\|\xi_p(\omega)\|_p^{p/s}
 \le \hcS(\xi(\omega))^{1/p-1/s}\|\xi_p(\omega)\|_p.
\]
The estimate for $s=\infty$ is already contained in \eqref{eq:nik-intermediate}. 
\end{proof}

\section{Quantum $az+b$}

In this section, we study the scaling Fourier constants for $1<p<2$ and $0<\theta<1$.
This shows that the quantum $az+b$ has ordinary Hausdorff-Young inequality but has nontrivial scaling Hausdorff-Young inequality.
Let $0<q<1$ be the deformation parameter, and write
$1<p<2$ and $r=p'$. 
In this section $q$ is not the conjugate exponent. 
Put
\[
 \Gamma_q=\{q^\ell z:\ell\in\bZ,\ z\in\T\},\qquad
 \cH=L^2(\Gamma_q)\cong\ell^2(\bZ^2).
\]
The radial measure is counting measure and the angular measure is the probability Haar measure. 
The orthonormal basis $\{e_{k, \ell}\}_{k, \ell}$ is defined as follows:
\[
 e_{k,\ell}(q^jz)=\mathbf1_{\{j=\ell\}}z^k.
\]
The affiliated normal generators $a, b$ act on $L^2(\Gamma_q)$ by
\begin{equation}\label{az:generators}
 a e_{k,\ell}=q^k e_{k,\ell-1},\qquad
 b e_{k,\ell}=q^\ell e_{k+1,\ell}.
\end{equation}
They form a regular \(q^2\)-pair in the sense of
\cite[Definition 2.1]{PuszSoltan}; in particular, on the algebraic span
of this basis,
\[
 ab=q^2ba,\qquad ab^*=b^*a.
\]
The comultiplication $\Delta$ is defined by
\begin{equation}\label{az:coproduct}
 \Delta(a)=a\otimes a,\qquad
 \Delta(b)=a\otimes b\dot{+}b\otimes1.
\end{equation}

We have that $\cM=L^\infty(\bG)=\cB(\cH)$ with a tracial weight $\Tr$.
Choose the right Haar weight $\psi$ and the corresponding left Haar weight $\varphi$ as
\begin{equation}\label{az:weights}
 \psi(x)=\Tr(|b|x|b|),\qquad
 \varphi(x)=\Tr(D_\varphi^{1/2}xD_\varphi^{1/2}),\qquad
 D_\varphi=q^2|a^{-1}b|^2.
\end{equation}
Here and below, positive products involving unbounded densities are
interpreted as normal semifinite weights. The right Haar formula and
the unitary antipode give the left formula; see
\cite[Section 7.2, Proposition 7.7]{Krajczok}. In these conventions the
unitary antipode sends \(b\) to \(-qa^{-1}b\), so applying it to the
right density produces the left density in \eqref{az:weights}.
Explicitly,
\[
 qa^{-1}b\,e_{k,\ell}=q^{\ell-k}e_{k+1,\ell+1},
 \qquad
 D_\varphi e_{k,\ell}=q^{2(\ell-k)}e_{k,\ell}.
\]
Thus the source Haar weight is nontracial. 

Finite matrix units lie in a common Tomita core. 
Indeed \(\sigma_t^\varphi(x)=D_\varphi^{it}xD_\varphi^{-it}\), and every matrix unit is an entire analytic eigenvector of this action with finite left and right Haar square norms. 
The same holds after taking its adjoint. 
Consequently the formulas below are valid on a dense Fourier core and extend by the general Hausdorff--Young inequality.

The self-duality of \(\Gamma_q\) is described by the bicharacter $\chi: \Gamma_q\times \Gamma_q\to \T$ by $\chi(q^kz,q^\ell w)=z^\ell w^k.$

The quantum exponential and its angular Fourier coefficients are
\begin{align}
 F_q(tz)&=\prod_{k=0}^\infty
       \frac{1+q^{2k}tz^{-1}}{1+q^{2k}tz},
       && t\in q^\bZ,\quad z\in\T,\label{az:quantum-exponential}\\
 c_n(t)&=\int_\T F_q(tz)z^{-n}\,dz,
 &&F_q(tz)=\sum_{n\in\bZ}c_n(t)z^n.\label{az:coefficients}
\end{align}
At removable points the continuous extension is understood. 
The function in \eqref{az:quantum-exponential} is unimodular, and
\begin{equation}\label{az:parseval}
 |c_n(t)|\le1,\qquad \sum_{n\in\bZ}|c_n(t)|^2=1.
\end{equation}
Its coefficients are real, as follows by replacing \(z\) with
\(\bar z\) and conjugating.

Let \((\widehat a,\widehat b)\) be the copy of \eqref{az:generators} acting on a second copy of \(\cH\).
The representation theorem of Pusz--So{\l}tan gives the unitary
\begin{equation}\label{az:U}
 U=F_q(\widehat b\otimes b)
   \chi(\widehat a\otimes1,1\otimes a),
 \qquad
 (\id\otimes\Delta)U=U_{12}U_{13}.
\end{equation}
The second leg is the original quantum group. The canonical pairing is
\begin{equation}\label{az:W-orientation}
 W=\Sigma U\Sigma\in \cM\overline\otimes\widehat \cM,
 \qquad
 \lambda(\omega)=(\omega\otimes\id)(W).
\end{equation}
There is no adjoint in \eqref{az:W-orientation}. 
This orientation follows from the universal-dual interpretation of \cite[Theorem 4.2]{PuszSoltan}. 
Krajczok's type I Plancherel description identifies the dual von Neumann algebra with \(B(H)\); the dual quantum group is the opposite quantum \(az+b\) group \cite[Section 7.2, footnote 8]{Krajczok}.
Thus, in the displayed dual generators,
\begin{equation}\label{az:dual-coproduct}
 \widehat\Delta(\widehat a)=\widehat a\otimes\widehat a,\qquad
 \widehat\Delta(\widehat b)
   = \widehat b\otimes\widehat a\dotplus 1\otimes\widehat b.
\end{equation}

\begin{lemma}\label{az:kernel}
The multiplicative unitary \(W\) acts by
\begin{equation}\label{az:W-kernel}
 W(e_{k,\ell}\otimes e_{r,s})
   =\sum_{n\in\bZ}c_n(q^{s+\ell-k-r})
       e_{k+n,\ell-r}\otimes e_{r+n,s-k}.
\end{equation}
\end{lemma}
\begin{proof}
The bicharacter in \eqref{az:U} acts as
\[
 \chi(\widehat a,a)(e_{r,s}\otimes e_{k,\ell})
   =e_{r,s-k}\otimes e_{k,\ell-r}.
\]
For example, the radial part of \(\widehat a\) selects the angular power in the second factor, and its phase translates the second radial coordinate; 
applying the definition of \(\chi\) gives exactly the two shifts displayed. 
The modulus of \(\widehat b\otimes b\) on this vector has eigenvalue \(q^{s+\ell-k-r}\), 
while its phase increments the first index in both tensor factors. 
Expanding the function \(F_q\) by \eqref{az:coefficients} and flipping the legs gives \eqref{az:W-kernel}. 
For the fixed input, the output vectors indexed by \(n\) are orthonormal, and \eqref{az:parseval} proves convergence.
\end{proof}

Putting \(\mu=q^{-1}>1\), and using Baaj's coefficients \cite[Appendix A, Proposition A.3(b)]{Baaj}, we obtain that
\begin{equation}\label{az:Baaj-change}
 g_\mu(\mu^m z)=F_q(q^{1-m}z^{-1}),\quad A(m,k)=\int_\T g_\mu(\mu^m z)z^{-k}\,dz, \quad 
 c_n(q^m)= A(1-m,-n).
\end{equation}
We have that for every \(n\in\bZ\),
\begin{equation}\label{az:weighted-sum}
 \sum_{m\in\bZ}q^{2m}|c_n(q^m)|^2=q^{2(n+1)}.
\end{equation}

Under the normalization \eqref{az:weights}, the dual left
Haar weight density is
\begin{equation}\label{az:Dhat}
 \widehat D=q^{-2}|\widehat b|^2,\qquad
       \widehat D e_{r,s}=q^{2(s-1)}e_{r,s}. 
\end{equation}

We use precisely the Haar normalizations of Section~4 of \cite{WY26b} \emph{Hausdorff--Young inequalities}:
\[
 D_\varphi e_{k,\ell}=q^{2(\ell-k)}e_{k,\ell},\qquad
 \widehat D e_{d,s}=q^{2(s-1)}e_{d,s}.
\]
Thus the spatial realizations are the ordinary Schatten operators
$X=aD_\varphi^{1/p}$ and
$\tcF_p(X)=\lambda(a\varphi)\widehat D^{1/r}$. The ordinary
sharp constant is $H_p=1$, by \cite[Theorem~4.3]{WY26b}.

The source scaling group is $\tau_t=\operatorname{Ad}|a|^{-2it}$, where
$|a|e_{k,\ell}=q^ke_{k,\ell}$, which is periodic with period $\pi/|\log q|$. 
Its fixed spatial subspace is therefore
\begin{equation}\label{azgap:fixed-source}
 \mathcal X_p=
 \bigoplus_{k\in\mathbb Z}^{\ell^p}
 S_p\bigl(\ell^2(\mathbb Z)\bigr),\qquad
 \|X\|_p^p=\sum_k\|X_k\|_{S_p}^p.
\end{equation}
Here $X_k$ is the block on the vectors $e_{k,\ell}$, indexed by
$\ell$. This is the entire fixed subspace, including all
off-diagonal entries inside each block.

Put
\begin{equation}\label{azgap:kernel}
 h_r(j)=q^{2(j-1)/r}c_0(q^j),\qquad j\in\mathbb Z,
\end{equation}
where $c_n$ is defined in \eqref{az:coefficients}. 
For a finite matrix $X$ in \eqref{azgap:fixed-source}, the Fourier image is
block diagonal in the first target coordinate $d$ and satisfies
\begin{equation}\label{azgap:fixed-transform}
 (\tcF_p X)_m(j,j')
    =\sum_{\tilde\ell\in\mathbb Z} h_r(j+\tilde\ell)(X_{j'-j})_{\tilde\ell+m,\tilde\ell}.
\end{equation}
To check every weight in this formula, take the physical input
$X=|e_{k,\ell}\rangle\langle e_{k,\ell'}|$. 
Its bounded coefficient is $q^{-2(\ell'-k)/p}X$, so its normal functional is the vector slice with coefficient $q^{2(\ell'-k)/r}$. 
In \cite[equation~(4.10)]{WY26b}, slicing the first input $e_{k,\ell}$ against the first output $e_{k,\ell'}$ forces $n=0$ and $m=\ell-\ell'$. 
Multiplication on the right by $\widehat D^{1/r}$ gives
\[
 (\tcF_pX)e_{\tilde{j},\tilde{k}}
 =1_{\{\tilde{j}=d\}}q^{2(\tilde{k}+\ell'-k-1)/r} c_0(q^{\tilde{k}+\ell'-k})e_{d,\tilde{k}-k}.
\]
Set $j=\tilde{k}-k$, $j'=\tilde{k}$, and $\tilde\ell=\ell'$ to obtain \eqref{azgap:fixed-transform}. 
In particular, $\|\tcF_p (X)\|_{r}=\displaystyle \left(\sum_m \|(\tcF_pX)_m\|_{S_r}^r\right)^{1/r}$.

\begin{lemma}\label{azgap:bilinear}
Define
\[
 S_r(q)=\sum_{j\in\mathbb Z}q^{j-1}|c_0(q^j)|^{r/2}.
\]
Whenever this sum is finite, we have 
\begin{equation}\label{azgap:fixed-bound}
 H_p^{\mathrm{fix}}(\bG)\le \|h_r\|_{\ell^{r/2}}=S_r(q)^{2/r}.
\end{equation}
In particular the sum is finite when $r\ge4$.
\end{lemma}
\begin{proof}
For an arbitrary finite scalar sequence $h$ in $\ell^1$ and finite block
matrix $X$, we define the bilinear map
\[
 \mathcal B(h,X)_m(j,j')
       =\sum_\ell h(j+\ell)(X_{j'-j})_{\ell+m,\ell}.
\]
There are two endpoint bounds, both with constant one:
\begin{align}
 \|\mathcal B(h,X)\|_{\ell^\infty(S_\infty)}
    &\le \|h\|_{\ell^\infty}\|X\|_{\ell^1(S_1)},
       \label{azgap:endpoint-one}\\
 \|\mathcal B(h,X)\|_{\ell^2(S_2)}
    &\le \|h\|_{\ell^1}\|X\|_{\ell^2(S_2)}.
       \label{azgap:endpoint-two}
\end{align}
For \eqref{azgap:endpoint-one}, hold $k=j'-j$ fixed. 
Each target block of this contribution is a weighted shift, of norm at most
\[
 \|h\|_\infty\sum_\ell|(X_k)_{\ell+m,\ell}|
           \le\|h\|_\infty\|X_k\|_{S_1}.
\]
The last inequality is trace duality with the shift whose
diagonal phases match these entries. Sum over $k$ and take the
supremum over $d$. For \eqref{azgap:endpoint-two}, put $v=u+k$
and apply the scalar convolution estimate
$\ell^1*\ell^2\to\ell^2$ to the $u$ variable:
\[
 \sum_{m,k,u}\left|\sum_\ell h(u+\ell)
                    (X_k)_{\ell+m,\ell}\right|^2
 \le\|h\|_1^2\sum_{m,k,\ell}|(X_k)_{\ell+m,\ell}|^2.
\]

Bilinear complex interpolation between these two bounds at
$\eta=2/r$, first for finite coordinate compressions with the
same endpoint constants and then by density, gives
\begin{equation}\label{azgap:interpolation}
 \|\mathcal B(h,X)\|_{\ell^r(S_r)}
      \le\|h\|_{\ell^{r/2}}\|X\|_{\ell^p(S_p)},
 \qquad \frac1p=1-\frac1r.
\end{equation}
Indeed the three interpolation couples are
$(\ell^\infty,\ell^1)$,
$(\ell^1(S_1),\ell^2(S_2))$, and
$(\ell^\infty(S_\infty),\ell^2(S_2))$.
These are ordinary sequence and tracial Schatten interpolation
couples.
Apply \eqref{azgap:interpolation} to truncations of $h_r$. 
Their limits have the entries in \eqref{azgap:fixed-transform}, so that formula and density of finite matrices prove \eqref{azgap:fixed-bound}.

For finiteness, $|c_0(q^j)|\le1$, while \eqref{az:weighted-sum} with $n=0$ says
\begin{equation}\label{azgap:weighted-Plancherel}
 \sum_jq^{2(j-1)}|c_0(q^j)|^2=1.
\end{equation}
For $j\ge1$, bound the summands of $S_r$ by $q^{j-1}$.
For $j\le0$ and $r\ge4$, bound them by
$q^{2(j-1)}|c_0(q^j)|^2$. Both majorants are summable.
\end{proof}

\begin{proposition} \label{azgap:near-one}
Suppose that $\bG$ is quantum $az+b$ and $1<p<2$.
For every deformation parameter $0<q<1$, one has
$H_p^{\mathrm{fix}}(\bG)<1$.
Consequently, for $0<\theta<1$,
\begin{equation}\label{azgap:all-exponents}
 H_{p,\theta}(\bG)<1.
\end{equation}
\end{proposition}
\begin{proof}
For every integer $j$, the unimodular function
$z\mapsto F_q(q^jz)$ is nonconstant, hence
$|c_0(q^j)|<1$. For $j\ge1$, nonconstancy follows, for example,
from the convergent logarithmic expansion of its product for
$q^j<1$. 
For $j\le0$, put $n=1-j\ge1$ and use
\[
 F_q(q^jz)=z^{-n}F_q(q^{n+1}z).
\]
The second factor has winding number zero, so the product has
winding number $-n$ and is nonconstant. The summable majorants
in the proof of Lemma~\ref{azgap:bilinear} now give
$S_r(q)\to0$ as $r\to\infty$, by dominated convergence.
In particular $S_r(q)<1$ for all sufficiently large $r$.
Lemma~\ref{azgap:bilinear} proves strictness at some exponent
$p_0>1$, and Proposition~\ref{azgap:fixed-interpolation} propagates it to every $1<p<2$.

The scaling actions have common period $P=\pi/|\log q|$.
Apply Theorem~\ref{near:periodic-criterion}, using $H_p=1$.
The final assertion follows from reflection symmetry and log convexity in the scaling parameter.
\end{proof}

\begin{remark}
We see that $H_{p, \theta}(\bG)$ is not a constant and its minimum at $\theta=1/2$ is strictly less than one.
\end{remark}

\begin{proposition}\label{prop:Young-gap}
Suppose that $\bG$ is a locally compact quantum group.
We have that
\begin{equation}\label{eq:Young-product-bound}
 \|\mathsf C_{p,s}^{\boldsymbol\theta}\|
 \le H_{p,\theta_1} H_{s,\theta_2} H_{r',\theta_3}
\end{equation}
At $r=2$ the last factor is one. 
In particular, the optimal bilinear norm is strictly less than one if at least one of the three factors on the right-hand side of \eqref{eq:Young-product-bound} is strictly less than one.
\end{proposition}
\begin{proof}
It follows from the H\"{o}lder inequality and scaling Hausdorff-Young inequality.    
\end{proof}

\begin{theorem}\label{thm:azb-all-vertices-one}
For the real quantum $az+b$, let $1<p,s<2$, $2\le r<\infty$, and $1+1/r=1/p+1/s$.
The bilinear norm of $\mathsf C_{p,s}^{\boldsymbol\theta}$ is one at each of the six vertices of $\Omega_{p,s}$ when $r>2$, and at each of the four vertices when $r=2$. In particular $\|\mathsf C_{p,s}^{\boldsymbol0}\|=1$.
\end{theorem}
\begin{proof}
For $1<t<2$ and $n<0$, let $j-\tilde{j}=n$ and choose
\[
 Y_{j,\tilde{j}}=|e_{0,1+n+\tilde{j}}\rangle\langle e_{n,1+n}|.
\]
We have that $\|Y_{j,\tilde{j}}\|_{S_t}=1$.
We claim that
\begin{equation}\label{az:peak}
 \bigl\|\tcF_t(Y_{j,\tilde{j}})-
 |e_{j,0}\rangle\langle e_{\tilde{j},0}|\bigr\|_{t'}\longrightarrow0
 \quad\text{as }n\to-\infty.
\end{equation}
Indeed, slicing \eqref{az:W-kernel} gives, for a general matrix unit $Y_{k,k', \ell, \ell'}=|e_{k,\ell}\rangle\langle e_{k',\ell'}|$, with
$n=k'-k$,
\[
 \tcF_t(Y_{k,k', \ell, \ell'})e_{\ell-\ell',v}
 =q^{2(v+\ell'-k'-1)/t'}c_n(q^{v+\ell'-k})e_{\ell-\ell'+n,v-k},
\]
and the transform vanishes on the other first-coordinate blocks.
At $v=0$ this gives the coefficient $\gamma_{-n}=c_{n}(q^{1+n})$ for $Y_{j,\tilde{j}}$. 
The product formula for $F_q$ implies
\[
 F_q(q^{1+n}z)=z^{n}F_q(q^{-n+1}z),\qquad \gamma_{-n}=c_0(q^{-n+1})\longrightarrow1.
\]
The transform is a weighted shift, and its $S_{t'}$ norm is at most
one. 
The sum of the $t'$-powers of all its other coefficients tends to zero. 
This proves \eqref{az:peak}.

The scaling constant of this quantum group is $1$. 
By reflection, $\tcF_t^{(1)}(J_t Y_{k,k', \ell, \ell'})=\tcF_t(Y_{k,k', \ell, \ell'})^*$. 
Thus for $j>\tilde{j}$, applying \eqref{az:peak} to $Y_{\tilde{j},j}$ shows that the rank-one matrix unit $|e_{j,0}\rangle\langle e_{\tilde{j},0}|$ is approximated in $S_{t'}$ by the image under $\tcF_t^{(1)}$ of a norm-one input. 
For $j<\tilde{j}$ the same is true at parameter zero.

Let $j_m,k_m,\ell_m$ be three distinct integers and choose norm-one inputs
$X_m,Y_m,Z_m$ so that their scaling Fourier images converge,
respectively, to
\[
 A_m=|e_{j_m,0}\rangle\langle e_{k_m,0}|,\quad
 B_m=|e_{k_m,0}\rangle\langle e_{\ell_m,0}|,\quad
 C_m=|e_{j_m,0}\rangle\langle e_{\ell_m,0}|,
\]
in $S_{p'},S_{s'},S_r$, where every pairwise distance between $j_m,k_m,\ell_m$ tends to infinity with $m$.
At parameter zero use the construction when the row index is smaller than the column index; 
at parameter one use reflection when it is larger. 
Since $A_mB_m=C_m$, the spatial H\"older inequality and Fourier duality give
\[
 \bigl\langle\mathsf C_{p,s}^{\boldsymbol\theta}(X_m,Y_m),Z_m\bigr\rangle
 \longrightarrow1.
\]
The six orderings of $j_m,k_m,\ell_m$ realize precisely the six effective vertices represented by
\[
 (0,0,0),\ (1,0,0),\ (1,0,1),\
 (1,1,1),\ (0,1,1),\ (0,1,0).
\]
Hence each corresponding bilinear norm is at least one; Young's
inequality gives the reverse bound.

For $r=2$, the third twist is ineffective. 
Choose $j_m,k_m,\ell_m$ so that the two input index inequalities realize any prescribed pair $(\theta_1,\theta_2)\in\{0,1\}^2$, and take $Z_m$ to be the exact inverse Plancherel transform of $C_m$. 
The same argument gives norm one at all four vertices of the rectangle.
\end{proof}

\begin{corollary}\label{az:exact-one-locus}
Under the strict input assumptions of Theorem~\ref{thm:azb-all-vertices-one}, 
the bilinear norm is one precisely at the effective vertices. 
At every other point of $\Omega_{p,s}$ it is strictly smaller than one.
\end{corollary}
\begin{proof}
For $r>2$, each Fourier exponent $p,s,r'$ is strict. Proposition~\ref{prop:Young-gap} and Proposition~\ref{azgap:near-one} show that a cube representation with any coordinate in $(0,1)$ has norm strictly smaller than one.
A nonvertex effective point has such a representation: otherwise it is the image of a cube corner. 
The only cube corners whose images are not effective vertices are $(0,0,1)$ and $(1,1,0)$. 
Their images are $(\kappa_{r'},\kappa_{r'})$ and $(\kappa_p,\kappa_s)$.
For $0<\varepsilon<\min(\kappa_p,\kappa_s,\kappa_{r'})$ the decompositions
\[
 (\kappa_{r'},\kappa_{r'})=(\varepsilon,\varepsilon)+(\kappa_{r'}-\varepsilon,\kappa_{r'}-\varepsilon),\qquad
 (\kappa_p,\kappa_s)=(\kappa_p-\varepsilon,\kappa_s-\varepsilon)+(\varepsilon,\varepsilon)
\]
represent them with all three effective cube times strictly interior. 
If $r=2$, every nonvertex point of the rectangle has an interior input coordinate. 
The same product bound proves strictness there. 
The theorem supplies equality of the optimal norm to one at all vertices.
\end{proof}

\begin{theorem}
\label{thm:low-quantitative-resolved}
Assume $\nu=1$ and use the analytic family and boundary identification proved in Theorem~\ref{thm:YF-low-output-family}. 
Let $1<p,s<r\le2$, $1+1/r=1/p+1/s$, and put
\[
 t=\frac{2}{r'},\qquad \tilde{p}'=\frac{2p'}{r'},\qquad
 \tilde{s}'=\frac{2s'}{r'}.
\]
Then $1<\tilde{p},\tilde{s}<2$, $1/\tilde{p}+1/\tilde{s}=3/2$, and
\begin{equation}\label{eq:low-quantitative-resolved}
 \|D_{p,s}^{\alpha,\beta}\|
 \le\bigl(H_{\tilde{p},\alpha}H_{\tilde{s},\beta}\bigr)^{2/r'}.
\end{equation}
For $1<p<2<s<r<\infty$ with $1/r+1=1/p+1/s$,
\begin{equation}\label{eq:IV-quantitative-resolved}
 \|D_{p,s}^{\mathrm{IV};\alpha,\beta}\|
 \le \left(H_{(2p'/s)',\alpha}H_{(2r/s)',\beta}\right)^{2/s}.
\end{equation}
For $1<s<2<p<r<\infty$ with $1/r+1=1/p+1/s$,
\begin{equation}\label{eq:III-quantitative-resolved}
 \|D_{p,s}^{\mathrm{III};\alpha,\beta}\|
 \le \left(H_{(2r/p)',\alpha}H_{(2s'/p)',\beta}\right)^{2/p}.
\end{equation}
All norms are the corresponding bilinear operator norms. 
In particular, for real quantum $az+b$, the estimates are strictly less than one if either $\alpha$ or $\beta$ is in $(0,1)$.
\end{theorem}

\begin{proof}
The exponent relation is equivalent to $1/p'+1/s'=1/r'$.
Since $p,s<r$, we have $p',s'>r'$, hence $\tilde{p}',\tilde{s}'>2$ and
\[
 \frac1{\tilde{p}}+\frac1{\tilde{s}}
 =2-\frac{r'}2\left(\frac1{p'}+\frac1{s'}\right)=\frac32.
\]
Moreover
\[
 \frac1p=1-t+\frac t{\tilde{p}},\qquad
 \frac1s=1-t+\frac t{\tilde{s}},\qquad
 \frac1r=1-\frac t2.
\]
The case $r=2$ has $t=1$ and follows directly from the  Plancherel identity and H\"older inequality. 
Suppose $r<2$. 
We apply bilinear Stein interpolation to the same analytic family as in \eqref{eq:YF-analytic-bilinear-family}, now with $p_0=\tilde{p},s_0=\tilde{s}$.
At the lower edge the predual convolution bound is $M_0=1$.
At the upper edge, the boundary identification from that proof and the factorization through the two scaling Fourier transforms give the stronger bound
\[
 M_1\le H_{\tilde{p},\alpha}H_{\tilde{s},\beta}.
\]
The real modular and scaling actions occurring on the boundary are
isometries on the pertinent endpoint spaces, so this bound is uniform
in the imaginary part of the strip parameter. 
Thus the interpolation gives $M_0^{1-t}M_1^t$ at $t$, and the dense core identity
\eqref{eq:YF-low-core-formula} identifies the interpolated map with
$D_{p,s}^{\alpha,\beta}$. This proves
\eqref{eq:low-quantitative-resolved}.

For Region IV, the defining duality identity gives the equality of
bilinear norms
\[
 \|D_{p,s}^{\mathrm{IV};\alpha,\beta}\|
 =\|D_{p,r'}^{\alpha,\beta}:L^p\times L^{r'}\to L^{s'}\|.
\]
Indeed both directions follow by taking suprema over the unit balls
in $L^p,L^s,L^{r'}$, and the involution is isometric. Since
$p,r'<s'<2$, apply the Region I estimate with $r$ replaced by
$s'$. 
This gives \eqref{eq:IV-quantitative-resolved}. 
Similarly,
\[
 \|D_{p,s}^{\mathrm{III};\alpha,\beta}\|
 =\|D_{r',s}^{\alpha,\beta}:L^{r'}\times L^s\to L^{p'}\|,
\]
and the Region I estimate with output $p'$ yields
\eqref{eq:III-quantitative-resolved}. 
For real quantum $az+b$,
Proposition~\ref{azgap:near-one} gives $H_{\mathfrak{p},\gamma}<1$ for $1<\mathfrak{p}<2$ and $0<\gamma<1$; 
the other factor is at most one. 
All Fourier exponents above lie strictly between one and two. 
This proves the asserted strictness.
\end{proof}

\section*{Declaration of generative AI use}

During the preparation of this work, the authors used generative AI
tools to improve the language and readability of the manuscript and
to assist with developing and checking mathematical arguments.
The authors reviewed and revised the AI-assisted content and verified
the mathematical statements and proofs.
The authors take full responsibility for the accuracy and integrity
of the manuscript.

\bibliographystyle{plain}
\bibliography{scaling}

@article{Bab61,
  author  = {Babenko, K. I.},
  title   = {An inequality in the theory of {F}ourier integrals},
  journal = {Izv. Akad. Nauk SSSR Ser. Mat.},
  volume  = {25},
  year    = {1961},
  pages   = {531--542}
}

@article{Baaj,
  author  = {Baaj, S.},
  title   = {Repr\'esentation r\'eguli\`ere du groupe quantique des d\'eplacements de {W}oronowicz},
  journal = {Ast\'erisque},
  volume  = {232},
  year    = {1995},
  pages   = {11--48},
  url     = {https://www.numdam.org/item/AST_1995__232__11_0.pdf},
}

@article{Bec75,
  author  = {Beckner, W.},
  title   = {Inequalities in {F}ourier analysis},
  journal = {Annals of Mathematics},
  volume  = {102},
  number  = {1},
  year    = {1975},
  pages   = {159--182}
}

@article{CMMP2019,
  author  = {Cowling, M. G. and Martini, A. and M{\"u}ller, D. and Parcet, J.},
  title   = {The {H}ausdorff--{Y}oung inequality on {L}ie groups},
  journal = {Math. Ann.},
  volume  = {375},
  year    = {2019},
  pages   = {93--131}
}

@article{Cas13,
  author  = {Caspers, M.},
  title   = {The {$L^p$}-{F}ourier transform on locally compact quantum groups},
  journal = {Journal of Operator Theory},
  volume  = {69},
  year    = {2013},
  pages   = {161--193}
}

@article{Con80,
  author  = {Connes, A.},
  title   = {On the spatial theory of von {N}eumann algebras},
  journal = {Journal of Functional Analysis},
  volume  = {35},
  year    = {1980},
  pages   = {153--164}
}

@article{Coo10,
  author  = {Cooney, T.},
  title   = {A {H}ausdorff--{Y}oung inequality for locally compact quantum groups},
  journal = {International Journal of Mathematics},
  volume  = {21},
  number  = {12},
  year    = {2010},
  pages   = {1619--1632}
}

@article{Clarkson36,
  author  = {Clarkson, James A.},
  title   = {Uniformly Convex Spaces},
  journal = {Transactions of the American Mathematical Society},
  volume  = {40},
  number  = {3},
  pages   = {396--414},
  year    = {1936},
  doi     = {10.2307/1989630}
}

@book{EnoSch92,
  author    = {Enock, M. and Schwartz, J.-M.},
  title     = {Kac algebras and duality of locally compact groups},
  publisher = {Springer-Verlag},
  year      = {1992}
}

@article{Hil81,
  author  = {Hilsum, M.},
  title   = {Les espaces {$L^p$} d'une alg{\`e}bre de von {N}eumann d{\'e}finies par la d{\'e}riv{\'e}e spatiale},
  journal = {Journal of Functional Analysis},
  volume  = {40},
  year    = {1981},
  pages   = {151--169}
}

@article{HJX,
  author  = {Haagerup, U. and Junge, M. and Xu, Q.},
  title   = {A reduction method for noncommutative {$L_p$}-spaces and applications},
  journal = {Trans. Amer. Math. Soc.},
  volume  = {362},
  year    = {2010},
  pages   = {2125--2165},
  url     = {https://arxiv.org/abs/0806.3635}
}

@article{Izu97,
  author  = {Izumi, H.},
  title   = {Constructions of non-commutative {$L^p$}-spaces with a complex parameter arising from modular actions},
  journal = {International Journal of Mathematics},
  volume  = {8},
  year    = {1997},
  pages   = {1029--1066}
}

@article{JJL20,
  author  = {Jaffe, A. and Jiang, C. and Liu, Z. and Ren, Y. and Wu, J.},
  title   = {Quantum {F}ourier analysis},
  journal = {PNAS},
  volume  = {117},
  number  = {20},
  year    = {2020},
  pages   = {10715--10720}
}

@article{JLW16,
  author  = {Jiang, C. and Liu, Z. and Wu, J.},
  title   = {Noncommutative uncertainty principles},
  journal = {Journal of Functional Analysis},
  volume  = {270},
  year    = {2016},
  pages   = {264--311}
}

@article{JLW18,
  author  = {Jiang, C. and Liu, Z. and Wu, J.},
  title   = {Uncertainty principles for locally compact quantum groups},
  journal = {Journal of Functional Analysis},
  volume  = {274},
  year    = {2018},
  pages   = {2399--2445}
}

@article{JLW19,
  author  = {Jiang, C. and Liu, Z. and Wu, J.},
  title   = {Block maps and {F}ourier analysis},
  journal = {Science China Mathematics},
  volume  = {62},
  year    = {2019},
  pages   = {1585--1614}
}

@article{Kahng,
    author ={B.J. Kahng} ,
    title ={A \((2n+1)\)-dimensional quantum group constructed from a skew-symmetric matrix} ,
    journal = {J. Geom. Phys.},
    year = {2011},
    volume={61},
    pages={2081--2097},
}

@article{KR,
    author = {A. Klein and B. Russo},
    title = {Sharp inequalities for {W}eyl operators and {H}eisenberg groups},
    journal = {Math. Ann.},
    year = {1978},
    volume={235},
    pages={175--194},
}

@article{Krajczok,
  author = {Krajczok, J.},
  title  = {Modular properties of type {I} locally compact quantum groups},
  journal={J. Operator Theory},
  volume={87},
  number={2},
  pages={319-354},
  url    = {10.7900/jot.2020sep23.2312},
  year={2022}
}

@article{KusVae00,
  author  = {Kustermans, J. and Vaes, S.},
  title   = {Locally compact quantum groups},
  journal = {Annales Scientifiques de l'{\'E}cole Normale Sup{\'e}rieure},
  volume  = {33},
  number  = {6},
  year    = {2000},
  pages   = {837--934}
}

@article{KusVae03,
  author  = {Kustermans, J. and Vaes, S.},
  title   = {Locally compact quantum groups in the von {N}eumann algebraic setting},
  journal = {Mathematica Scandinavica},
  volume  = {92},
  number  = {1},
  year    = {2003},
  pages   = {68--92}
}

@article{Kunze1958,
  author  = {Kunze, R. A.},
  title   = {{$L^{p}$} {F}ourier transforms on locally compact unimodular groups},
  journal = {Trans. Amer. Math. Soc.},
  volume  = {89},
  year    = {1958},
  pages   = {519--540}
}

@article{Lie90,
  author  = {Lieb, E.},
  title   = {Gaussian kernels have only {G}aussian maximizers},
  journal = {Inventiones Mathematicae},
  volume  = {102},
  year    = {1990},
  pages   = {179--208}
}

@misc{LPW19,
  author        = {Liu, Z. and Palcoux, S. and Wu, J.},
  title         = {Fusion bialgebras and {F}ourier analysis},
  year          = {2019},
  eprint        = {1910.12059},
  archivePrefix = {arXiv},
  note          = {arXiv:1910.12059}
}

@article{LWW17,
  author  = {Liu, Z. and Wang, S. and Wu, J.},
  title   = {Young's inequality for locally compact quantum groups},
  journal = {Journal of Operator Theory},
  volume  = {77},
  year    = {2017},
  pages   = {109--131}
}

@article{LiuWu18,
  author  = {Liu, Z. and Wu, J.},
  title   = {Extremal pairs of {Y}oung's inequality for {K}ac algebras},
  journal = {Pacific Journal of Mathematics},
  volume  = {295},
  number  = {1},
  year    = {2018},
  pages   = {103--121}
}

@article{LiuWu17,
  author  = {Liu, Z. and Wu, J.},
  title   = {Uncertainty principles for {K}ac algebras},
  journal = {Journal of Mathematical Physics},
  volume  = {58},
  year    = {2017},
  pages   = {052102}
}

@article{Mc67,
  author  = {McCarthy, Charles A.},
  title   = {\(C_p\)},
  journal = {Israel Journal of Mathematics},
  volume  = {5},
  number  = {4},
  pages   = {249--271},
  year    = {1967},
  doi     = {10.1007/BF02771613}
}

@article{Nei63,
  author  = {O'Neil, R.},
  title   = {Convolution operators and {$L(p,q)$} spaces},
  journal = {Duke Math Journal},
  volume  = {30},
  year    = {1963},
  pages   = {129--142}
}

@article{PuszSoltan,
  author = {Pusz, W. and So{\l}tan, P. M.},
  title  = {Functional form of unitary representations of the quantum {$az+b$} group},
  journal={Rep. Math. Phys.},
  year={2003},
  pages={309-319},
  volume={52},
  issue={2},
  url    = {https://arxiv.org/abs/math/0604613},
}

@book{Str81,
  author    = {Stratila, S.},
  title     = {Modular Theory in Operator Algebras},
  publisher = {Abacus Press},
  year      = {1981}
}

@article{Rus74,
  author  = {Russo, B.},
  title   = {The norm of {$L_p$}-{F}ourier transform on unimodular groups},
  journal = {Transactions of the American Mathematical Society},
  volume  = {192},
  number  = {4},
  year    = {1974},
  pages   = {293--305}
}

@article{Stein56,
  title={Interpolation of linear operators},
  author={Stein, Elias M.},
  journal={Transactions of the American Mathematical Society},
  volume={83},
  number={2},
  pages={482--492},
  year={1956},
  doi={10.1090/S0002-9947-1956-0082586-0}
}

@article{Ter81,
  author       = {Terp, M.},
  title        = {{$L^p$} spaces associated with von {N}eumann algebras},
  journal = {Notes, Report No.~3a+3b, K{\o}benhavns Universitets Matematiske Institut},
  year         = {1981}
}

@article{Ter82,
  author  = {Terp, M.},
  title   = {Interpolation spaces between a von {N}eumann algebra and its predual},
  journal = {Journal of Operator Theory},
  volume  = {8},
  year    = {1982},
  pages   = {327--360}
}

@article{Terp2017,
  author  = {Terp, M.},
  title   = {{$L^{p}$} {F}ourier transformation on non-unimodular locally compact groups},
  journal = {Adv. Oper. Theory},
  volume  = {2},
  year    = {2017},
  pages   = {547--583}
}

@article{You26,
  author        = {Youn, S.-G.},
  title         = {Twisted {F}ourier transforms on non-{K}ac compact quantum groups},
  year          = {2026},
  eprint        = {2503.23316},
  archivePrefix = {arXiv},
  journal          = {arXiv:2503.23316v2},
  url           = {https://arxiv.org/abs/2503.23316v2}
}

@article{VaiKac74,
  author  = {Vainerman, L. I. and Kac, G. I.},
  title   = {Nonunimodular Ring Groups and {H}opf--von {N}eumann Algebras},
  journal = {Sbornik: Mathematics},
  volume  = {23},
  number  = {2},
  year    = {1974},
  pages   = {185--214}
}

@article{WY26a,
    author = {J. Wu and F. Yang},
    title ={Sharp {H}ausdorff–{Y}oung inequalities for solvable and semisimple {L}ie groups} ,
    journal ={ArXiv:2611.02597} ,
    year = {2026},
}

@article{WY26b,
    author = {J. Wu and F. Yang},
    title ={Sharp {H}ausdorff–{Y}oung inequalities for some locally compact quantum groups} ,
    journal ={ArXiv:2611.05863v1} ,
    year = {2026},
}

@article{Yap69,
  author  = {Yap, L.},
  title   = {Some remarks on convolution operators and {$L(p,q)$} spaces},
  journal = {Duke Math. Journal},
  volume  = {36},
  year    = {1969},
  pages   = {647--658}
}

@article{You12,
  author  = {Young, W. H.},
  title   = {On the multiplication of successions of {F}ourier constants},
  journal = {Proceedings of the Royal Society Series A},
  volume  = {87},
  year    = {1912},
  pages   = {331--339}
}

@article{Nik51,
 author = {Nikol'skii, S. M.},
 title = {Inequalities for entire functions of finite degree and their application to the theory of differentiable functions of several variables},
 journal = {Trudy Matematicheskogo Instituta imeni V. A. Steklova},
 volume = {38}, pages = {244--278}, year = {1951},
 note = {In Russian}
}

@article{AkyRuz20,
 author = {Akylzhanov, Rauan and Ruzhansky, Michael},
 title = {{$L^p$--$L^q$ multipliers on locally compact groups}},
 journal = {Journal of Functional Analysis}, volume = {278}, number = {3},
 pages = {108324}, year = {2020}
}

@article{QX,
    author = {H. Qiu and X. Xiong},
    title = {Noncommutative sharp {H}ausdorff--{Y}oung inequality},
    journal = {ArXiv:2609.10424v2},
    year = {2026}
}
\end{document}